\documentclass[preprint,3p,sort&compress,final,times]{elsarticle}
\usepackage[overload]{empheq}
\usepackage{amssymb}
\usepackage{mathtools}
\usepackage{bm}
\usepackage[usenames,dvipsnames]{color}

\usepackage{tikz}
\usetikzlibrary{decorations.pathreplacing,decorations.markings}

\newcommand{\innerprod}[2]{\left\langle #2\,,\, #1\right\rangle}

\newcommand{\ex}{\boldsymbol{e}_{\lambda}}
\newcommand{\ey}{\boldsymbol{e}_{\varphi}}
\newcommand{\ez}{\boldsymbol{e}_{z}}

\newcommand{\ux}{u}
\newcommand{\uy}{v}
\newcommand{\uz}{w}
\newcommand{\upar}{\boldsymbol{u}_{\parallel}}
\newcommand{\uperp}{\boldsymbol{u}_{\perp}}

\newcommand{\omegapar}{\boldsymbol{\omega}_{\parallel}}
\newcommand{\omegaperp}{\boldsymbol{\omega}_{\perp}}
\newcommand{\omegaparpar}{\boldsymbol{\omega}_{\parallel,\parallel}}
\newcommand{\omegaparperp}{\boldsymbol{\omega}_{\parallel,\perp}}

\newcommand{\gradpar}{\nabla_{\parallel}}
\newcommand{\gradperp}{\nabla_{\perp}}
\newcommand{\curlpar}{\nabla_{\parallel}\times}
\newcommand{\curlperp}{\nabla_{\perp}\times}

\newcommand{\pspace}{\mathcal{P}_{h}}
\newcommand{\wspace}{\mathcal{W}_{h}}
\newcommand{\wspacepar}{\mathcal{W}_{\parallel,h}}
\newcommand{\wspaceperp}{\mathcal{W}_{\perp,h}}
\newcommand{\uspace}{\mathcal{U}_{h}}
\newcommand{\uspacepar}{\mathcal{U}_{\parallel,h}}
\newcommand{\uspaceperp}{\mathcal{U}_{\perp,h}}
\newcommand{\qspace}{\mathcal{Q}_{h}}

\newcommand{\pcardinal}{d_{\mathcal{P}}}
\newcommand{\wcardinal}{d_{\mathcal{W}}}
\newcommand{\wcardinalpar}{d_{\mathcal{W}_{\parallel}}}
\newcommand{\wcardinalperp}{d_{\mathcal{W}_{\perp}}}
\newcommand{\ucardinal}{d_{\mathcal{U}}}

\newcommand{\qcardinal}{d_{\mathcal{Q}}}

\newcommand{\pbasis}{\epsilon^{\mathcal{P}}}
\newcommand{\wbasis}{\boldsymbol{\epsilon}^{\mathcal{W}}}
\newcommand{\wbasispar}{\boldsymbol{\epsilon}^{\mathcal{W}_{\parallel}}}
\newcommand{\wbasisperp}{\boldsymbol{\epsilon}^{\mathcal{W}_{\perp}}}
\newcommand{\ubasis}{\boldsymbol{\epsilon}^{\mathcal{U}}}
\newcommand{\ubasispar}{\boldsymbol{\epsilon}^{\mathcal{U}_{\parallel}}}
\newcommand{\ubasisperp}{\boldsymbol{\epsilon}^{\mathcal{U}_{\perp}}}
\newcommand{\qbasis}{\epsilon^{\mathcal{Q}}}

\newcommand{\ualgebraicpar}{\boldsymbol{\mathsf{u}}^{\parallel}}
\newcommand{\ualgebraicperp}{\boldsymbol{\mathsf{u}}^{\perp}}

\newcommand{\massu}{\boldsymbol{\mathsf{M}}^{\mathcal{U}}}
\newcommand{\massupar}{\boldsymbol{\mathsf{M}}^{\mathcal{U}_{\parallel}}}
\newcommand{\nupar}{\boldsymbol{\mathsf{N}}^{\mathcal{U}_{\parallel}}}
\newcommand{\massuperp}{\boldsymbol{\mathsf{M}}^{\mathcal{U}_{\perp}}}
\newcommand{\nuperp}{\boldsymbol{\mathsf{N}}^{\mathcal{U}_{\perp}}}

\newcommand{\masswpar}{\boldsymbol{\mathsf{M}}^{\mathcal{W}_{\parallel}}}
\newcommand{\masswperp}{\boldsymbol{\mathsf{M}}^{\mathcal{W}_{\perp}}}
\newcommand{\massq}{\boldsymbol{\mathsf{M}}^{\mathcal{Q}}}
\newcommand{\luq}{\boldsymbol{\mathsf{L}}^{\mathcal{U}_{\perp},\mathcal{Q}}}

\newcommand{\Ualgebraicperp}[1]{\boldsymbol{\mathsf{u}}^{\perp, #1}}
\newcommand{\Nuperp}[1]{\boldsymbol{\mathsf{N}}^{\mathcal{U}_{\perp, #1}}}

\begin{document}
\begin{frontmatter}

\title{An energetically balanced, quasi-Newton integrator for non-hydrostatic vertical atmospheric dynamics}
\author[MON]{David Lee\corref{cor}}
\ead{davelee2804@gmail.com}

\address[MON]{Department of Mechanical and Aerospace Engineering, Monash University, Melbourne 3800, Australia}
\cortext[cor]{Corresponding author. Tel. +61 452 262 804.}

\begin{abstract}
An energetically balanced, implicit integrator for non-hydrostatic vertical atmospheric dynamics on the 
sphere is presented. The integrator allows for the exact balance of energy exchanges in space and time 
for vertical atmospheric motions by preserving the skew-symmetry of the non-canonical Hamiltonian 
formulation of the compressible Euler equations. The performance of the integrator is accelerated by a 
preconditioning strategy that reduces the dimensionality of the inner linear system.
Here we reduce the four component velocity, density, density weighted potential temperature
and Exner pressure system into a single equation for the density weighted potential temperature via
repeated Schur complement decomposition and the careful selection of coupling terms. 
As currently formulated, the integrator is based on a horizontal-vertical spatial splitting
that does not permit bottom topography.
The integrator is
validated for standard test cases for baroclinic instability and a non-hydrostatic gravity wave on the 
sphere and a rising bubble in a high-resolution planar geometry, and shows robust convergence 
across all of these regimes.
\end{abstract}

\begin{keyword}
Poisson integrator\sep
Non-hydrostatic\sep
Euler equations\sep
Cubed sphere\sep
Horizontally explicit/vertically implicit\sep
Quasi-Newton
\end{keyword}

\end{frontmatter}

\section{Introduction}

Energy conserving, semi-implicit time integrators for various geophysical systems in Hamiltonian form have been 
previously introduced for the shallow water \cite{Bauer18,Wimmer20} and thermal shallow water \cite{Eldred18,WCB20}
equations on the sphere, and the compressible Euler equations in vertical slice geometry \cite{WCB20}. 
They have also been implemented on other areas, such as plasma physics \cite{KS19}.
These integrators are based on discrete gradient methods \cite{McLachlan99,Cohen11} for non-canonical 
Hamiltonian systems which preserve the skew-symmetric property of the Poisson bracket in the spatial 
discretisation and the exact integration of the variational derivatives of the Hamiltonian in the temporal discretisation. 
In doing so these integrators satisfy the exact balance of energy exchanges, 
as well as the orthogonality of the vorticity dynamics to these energetic exchanges. By mitigating against 
internal biases in the representation of dynamical processes these integrators may help to improve the 
statistical quality of long time integrations.

In the present article a skew-symmetric mimetic spectral element spatial discretisation of the 3D compressible 
Euler equations on the sphere \cite{LP20} is extended to incorporate such an integrator for non-hydrostatic 
vertical dynamics. The vertical integrator is coupled to the explicit horizontal dynamics using a trapazoidal 
horizontally explicit/vertically implicit (HEVI) splitting scheme \cite{Weller13,Lock14}, which 
complements the centered time integration of the vertical scheme. Skew-symmetric formulations have 
previously been used as the basis for energy conserving spatial discretisations of vertical dynamics in 
non-hydrostatic atmospheric modelling \cite{Taylor19}. In order to accelerate the performance of 
such an integrator, Schur complement reduction may be applied in order to transform the coupled nonlinear system into a 
single equation at each linear solve \cite{Melvin19,MMM20}. In the present formulation the coupled velocity, density, 
density weighted potential temperature, Exner pressure system is reduced to a single equation for the density 
weighted potential temperature, for which the implicit Helmholtz operator incorporates corrections to both the 
gradient and divergence operators in order to account for the thermodynamics. This operator is derived from the 
reduction of an approximate Jacobian for a deliberate choice of coupling terms.

The motivation for such an integrator is twofold. Firstly, almost all of the energy of the atmosphere exists as
either potential or internal energy, which are dependent on the density and potential temperature (or equivalent 
thermodynamic variables). These quantities are strongly stratified, such that their gradients are much sharper 
in the vertical dimension than the horizontal. As such it is hoped that by balancing energetic exchanges 
associated with transient vertical processes this will improve the overall representation of dynamics within the 
model. Secondly, by developing such an integrator and demonstrating its robust performance across different 
scales and physical regimes, this may inform the development of a fully three-dimensional energetically balanced
integrator.

In section 2 the 3D compressible Euler equations are introduced in continuous form, including a brief 
discussion of their energetic properties. Section 3 briefly describes the mixed mimetic spectral
element spatial discretisation of the 3D compressible Euler equations \cite{LP20}. For more extensive 
discussions of this discretisation \cite{Gerritsma11,Kreeft13} and its properties for geophysical flow 
simulation \cite{LPG18,LP18} the reader is referred to the aforementioned articles. 
The energetic properties of the spatially and temporally discrete system are discussed in section 4.
Section 5 describes the formulation of the quasi-Newton vertical integrator, including the nonlinear
preconditioning strategy. Results for the reproduction of standard test cases are presented in section 6, 
and finally conclusions and discussion of future work are presented in section 7.

\section{The rotating 3D compressible Euler equations}

The three dimensional compressible Euler equations for a shallow atmosphere may be expressed as \cite{Nevir09,Gassmann13,LP20}
\begin{subequations}\label{eq:compressible_euler}
\begin{align}
\frac{\partial\bm u}{\partial t} + (\bm\omega + \bm f)\times\bm u + \nabla\Bigg(\frac{1}{2}\|u\|^2 + gz\Bigg) + \theta\nabla\Pi &= 0,\label{eq::mom}\\
\frac{\partial\rho}{\partial t} + \nabla\cdot(\rho\bm u) &= 0,\label{eq::dens}\\
\frac{\partial\rho\theta}{\partial t} + \nabla\cdot(\rho\theta \bm u) &= 0,\label{eq::temp}
\end{align}
\end{subequations}
where $\bm u = u\ex + v\ey + w\ez$ are the zonal, meridional and vertical velocity components respectively,
$\rho$ is the density, $\bm f=f\ez$ is the Coriolis term, $g$ is the acceleration due to gravity, $\theta$
is the potential temperature, and $\Pi$ is the Exner pressure (including the specific heat at constant pressure).
The potential temperature and Exner pressure are defined with respect to the standard thermodynamic variables 
of temperature, $T$, and pressure, $p$, as
\begin{subequations}
	\label{eq:definition_exner_potential_temperature}
	\begin{align}
		\Pi &:= c_p\Bigg(\frac{p}{p_0}\Bigg)^{\frac{R}{c_p}}, \label{eq:exener_definition} \\
	 	\theta &:= \frac{c_{p}T}{\Pi}, \label{eq:potential_temperature_definition}
	\end{align}
\end{subequations}
where $c_p$ is the specific heat at constant pressure, $p_0$ is the reference pressure, and $R = c_p-c_v$ 
is the ideal gas constant (and $c_v$ is the specific heat at constant volume). Using \eqref{eq:definition_exner_potential_temperature}
the ideal gas law, $p=R\rho T$ may be reformulated as
\begin{equation}\label{eq::eos}
\Pi = c_p\Bigg(\frac{R\rho\theta}{p_0}\Bigg)^\frac{R}{c_v}.
\end{equation}

To obtain a closed system for the solution of the compressible Euler equations, \eqref{eq:compressible_euler} and
\eqref{eq::eos} must be supplemented by appropriate boundary conditions. We impose homogeneous Dirichlet boundary 
conditions on the $z$-component of velocity, $w$ as
\begin{equation}
	\left.w\right|_{z = 0} = \left.w\right|_{z = z^{\mathrm{top}}} = 0, \label{eq:dirichlet_bc_velocity}
\end{equation}
where $z^{\mathrm{top}}$ corresponds to the $z$-coordinates of the top boundary of the domain. We also apply
homogeneous Neumann boundary conditions on the Exner pressure as
\begin{equation}\label{eq::bcs_neumann}
	\frac{\partial\Pi}{\partial z}\Bigg|_{z=0} = \frac{\partial\Pi}{\partial z}\Bigg|_{z=z^{top}} = 0.
\end{equation}

Note that in this formulation we have invoked the shallow atmosphere approximation, for which gravity is 
constant throughout the fluid column and the height of the fluid column is negligible with respect to the 
earth's radius.
We further assume that the horizontal components of the Coriolis term are small and omit these 
from the formulation \cite{White05}.

\subsection{Energetics}

In this section the energetics of the 3D compressible Euler equations will be briefly analysed in
continuous form. The energetics of the discrete formulation will be discussed in subsequent sections.

\subsubsection{Kinetic, potential, and internal energy}
	The kinetic energy, $K$, is defined as
	\begin{equation}
		K := \frac{1}{2}\langle\boldsymbol{u}, \rho\boldsymbol{u} \rangle = \frac{1}{2}\int_{\Omega}\rho \|u\|^{2}, \label{eq:kinetic_energy_definition}
	\end{equation}
	where $\|u\| := \langle \boldsymbol{u}, \boldsymbol{u} \rangle$, and $\langle \cdot, \cdot \rangle$ is the $L^{2}$ inner product given 
	for scalar fields as
	\begin{equation}
		\langle f, g \rangle := \int_{\Omega} f g\,\mathrm{d}\Omega\,,
	\end{equation}
and for vector fields as
	\begin{equation}
		\langle \boldsymbol{u}, \boldsymbol{v} \rangle := \int_{\Omega} \boldsymbol{u} \cdot \boldsymbol{v}\,\mathrm{d}\Omega.
	\end{equation}

The time variation of kinetic energy is obtained by summing the $L^{2}$ inner product, between the momentum equation, \eqref{eq::mom},
and $\rho\bm u$, and between the continuity equation, \eqref{eq::dens}, and $\frac{1}{2}\|u\|^2$
\begin{equation}
	\frac{\partial K}{\partial t}  = -\langle g ,\rho \uz\rangle - \langle\rho\bm u,\theta\nabla\Pi\rangle\,,\label{eq::ke_evol}
\end{equation}
where again $\uz$ is the $z$-component of the velocity field, $\boldsymbol{u}$.

	The potential energy, $P$, is given by
	\begin{equation}
		P := \innerprod{gz}{\rho} = \int_{\Omega}\rho gz \,\mathrm{d}\Omega, \label{eq:potential_energy_definition}
	\end{equation}
	and its time derivative follows directly
\begin{equation}
\frac{\partial P}{\partial t} = \innerprod{\frac{\partial\rho}{\partial t}}{gz} \stackrel{\eqref{eq::dens}}{=}
-\langle gz,\nabla\cdot(\rho\bm u)\rangle
= \langle g,\rho \uz\rangle\,,\label{eq:time_variation_potential_energy}
\end{equation}
where we have used integration by parts on the last identity and assumed periodic boundary conditions in the horizontal directions, together with homogeneous boundary conditions for the vertical component of the velocity field, \eqref{eq:dirichlet_bc_velocity}.

	The internal energy, $I$,  is defined as
	\begin{equation}
I := \int_{\Omega}c_{v}\rho T \,\mathrm{d}\Omega\stackrel{\eqref{eq:potential_temperature_definition}}{=} \int_{\Omega}\frac{c_{v}}{c_{p}}\rho\theta\Pi \,\mathrm{d}\Omega = \int_{\Omega}c_{v}\rho\theta\left(\frac{R\rho\theta}{p_{0}}\right)^{\frac{R}{c_{v}}} \,\mathrm{d}\Omega= \int_{\Omega} c_{v} \left(\frac{R}{p_{0}}\right)^{\frac{R}{c_{v}}}\left(\rho\theta\right)^{\frac{c_{p}}{c_{v}}}\,\mathrm{d}\Omega\,. \label{eq:internal_energy_definition}
	\end{equation}
After some manipulation, the time variation of internal energy is given by
	\begin{equation}
		\frac{\partial I}{\partial t} = - \innerprod{\Pi}{\nabla\cdot(\rho\theta \bm u)} = \innerprod{\theta\nabla\Pi}{\rho\bm u} \,. \label{eq:time_variation_internal_energy}
	\end{equation}
	where integration by parts was used on the last identity, together with homogeneous boundary conditions for $\boldsymbol{u}$ and periodic boundary conditions on the horizontal directions.
Note that the substitution of the boundary conditions as defined in \eqref{eq:dirichlet_bc_velocity} and \eqref{eq::bcs_neumann}
into \eqref{eq::ke_evol}, \eqref{eq:time_variation_potential_energy} and \eqref{eq:time_variation_internal_energy}
ensures that there is no temporal change in energy at the boundary. Mass conservation is also ensured by the substitution of
\eqref{eq:dirichlet_bc_velocity} into \eqref{eq::dens}.

\subsubsection{Conservation of total energy}
Following \cite{Gassmann13}, the total energy, $\mathcal{H}$, is given as the sum of kinetic, $K$, potential, $P$, and internal, $I$, energy as
 \begin{equation}
\mathcal{H} := K + P + I = \int_{\Omega}\frac{1}{2}\rho u^2\,\mathrm{d}\Omega +  \int_{\Omega}\rho gz\,\mathrm{d}\Omega +  \int_{\Omega}\frac{c_v}{c_p}\Theta\Pi\,\mathrm{d}\Omega\,. \label{eq::tot_en}
 \end{equation}
where $\Theta := \rho\theta$. The variational derivatives of $\mathcal{H}$ with respect to the prognostic variables 
$\boldsymbol{u}$, $\rho$, and $\Theta$, are given as
\begin{equation}\label{eq::variations}
	\frac{\delta\mathcal{H}}{\delta\bm u} = \rho\bm u = \bm U,\qquad
	\frac{\delta\mathcal{H}}{\delta\rho} = \frac{1}{2}u^2 + gz = \Phi,\qquad
	\frac{\delta\mathcal{H}}{\delta\Theta} = \Pi.
\end{equation}
Writing the prognostic variables as a column vector of the form, $\boldsymbol{\mathsf{a}} = [\boldsymbol{u}\quad\rho\quad\Theta]^{\top}$,
we may then express the original system \eqref{eq:compressible_euler} as \cite{Hairer06}
\begin{equation}\label{eq:da_dt_B_H}
\frac{\partial\boldsymbol{\mathsf{a}}}{\partial t} = \boldsymbol{\mathsf{B}} \,\frac{\delta\mathcal{H}}{\delta\boldsymbol{\mathsf{a}}}, 
\end{equation}
where $\boldsymbol{\mathsf{B}}$ is a skew-symmetric operator of the form
\begin{equation}\label{eq:en_con_3}
	\boldsymbol{\mathsf{B}} :=
	\left[
		\begin{array}{ccc}
			-\bm q\times(\cdot) & -\nabla(\cdot) & -\theta\nabla(\cdot)\\
             -\nabla\cdot(\cdot) & 0 & 0 \\
             -\nabla\cdot(\theta\cdot) & 0 & 0\\
		\end{array}
	\right],
\end{equation}
and $\boldsymbol{q} = (\boldsymbol{\omega} + \boldsymbol{f})/\rho$ is the potential vorticity.
Energy conservation is then assured since
\begin{equation}\label{eq::conservation}
\frac{\partial\mathcal{H}}{\partial t} =
\frac{\delta\mathcal{H}}{\delta\boldsymbol{\mathsf{a}}}
\cdot\frac{\partial\boldsymbol{\mathsf{a}}}{\partial t} \stackrel{\eqref{eq:da_dt_B_H}}{=} \frac{\delta\mathcal{H}}{\delta\boldsymbol{\mathsf{a}}}\cdot\left(
 \boldsymbol{\mathsf{B}}\frac{\delta\mathcal{H}}{\delta\boldsymbol{\mathsf{a}}}\right) = 0,
\end{equation}
where the last identity follows from the skew-symmetry of $\boldsymbol{\mathsf{B}}$.

\section{Spatial discretisation}

In this section we briefly describe some of the important features of the mixed mimetic spectral 
element spatial discretisation used in this work. For a complete description of this discretisation
as it applies to the 3D compressible Euler equations on the cubed sphere see \cite{LP20}.

The mixed mimetic spectral element method \cite{Gerritsma11,Kreeft13} is a compatible family of
finite element spaces. For finite dimensional polynomial spaces of the form
$\pspace(\Omega) \subset H^{1}(\Omega)$, $\wspace(\Omega)\subset H(\mathrm{curl}, \Omega)$, 
$\uspace(\Omega) \subset H(\mathrm{div}, \Omega)$, and $\qspace(\Omega) \subset L^{2}(\Omega)$ 
the discrete analogues of the differential operators satisfy a series of compatibility relations
of the form
\begin{equation}\label{eq:de_rham_sequence_discrete}
\mathbb{R}\longrightarrow\pspace(\Omega) \stackrel{\nabla}{\longrightarrow} \wspace(\Omega) 
\stackrel{\nabla\times}{\longrightarrow} \uspace(\Omega)  \stackrel{\nabla\cdot}{\longrightarrow} \qspace(\Omega) \longrightarrow 0, 
\end{equation}
where the space $H^{1}(\Omega)$ represents square integrable functions over $\Omega$ whose gradient 
is also square integrable, $H(\mathrm{curl}, \Omega)$ and $H(\mathrm{div}, \Omega)$ contain square 
integrable vector fields over $\Omega$ with square integrable curl and divergence respectively, and 
the function space $L^{2}(\Omega)$ contains square integrable functions. Each of these polynomial 
function spaces has an associated finite set of basis functions $\pbasis_{i}$, $\wbasis_{i}$, 
$\ubasis_{i}$, and $\qbasis_{i}$, such that
\begin{equation}\label{eq:function_spaces_span}
\pspace = \mathrm{span}\{\pbasis_{1}, \dots, \pbasis_{\pcardinal}\}, \quad 
\wspace = \mathrm{span}\{\wbasis_{1}, \dots, \wbasis_{\wcardinal}\}, \quad 
\uspace = \mathrm{span}\{\ubasis_{1}, \dots, \ubasis_{\ucardinal}\}, \quad\mathrm{and}
\quad\qspace = \mathrm{span}\{\qbasis_{1}, \dots, \qbasis_{\qcardinal}\}. 
\end{equation}
The discrete mappings given in \eqref{eq:de_rham_sequence_discrete} are represented by
strong form topological relations as \cite{Gerritsma11}:
\begin{equation}
\nabla\pbasis_{j} = \sum_{k=0}^{\wcardinal}\mathsf{E}_{k,j}^{1,0}\wbasis_{k},\qquad\nabla\times\wbasis_{j} = \sum_{k=0}^{\ucardinal}\mathsf{E}_{k,j}^{2,1}\ubasis_{k},\qquad \mathrm{and} \qquad \nabla\cdot\ubasis_{j} = \sum_{k=0}^{\qcardinal}\mathsf{E}_{k,j}^{3,2}\qbasis_{k}, \label{eq_incidence_matrices}
\end{equation}
where $\boldsymbol{\mathsf{E}}^{1,0}$, $\boldsymbol{\mathsf{E}}^{2,1}$, and $\boldsymbol{\mathsf{E}}^{3,2}$ 
are so called \emph{incidence} matrices corresponding to the discrete versions of the differential operators 
grad, curl and div. Conversely, reverse mappings may be applied through weak form adjoints to these operators
\cite{Kreeft13}.

\subsection{Directional splitting}

Here we briefly describe the splitting procedure used to decompose the three dimensional problem into
horizontal (two dimensional) and vertical sub-problems. Due to the nature of the various function spaces,
this splitting implies no additional approximation to the spatial discretisation. However this construction
does involve a temporal approximation, by which the horizontal and vertical sub problems are solved 
explicitly and implicitly respectively \cite{Weller13,Lock14}.

We defined the horizontal, $\upar$, and vertical, $\uperp$, components of the velocity field $\boldsymbol{u} = u\ex + v\ey + w\ez$ as
\begin{equation}
\upar := \ux\ex + \uy\ey\,, \qquad \uperp := \uz\ez\,. \label{eq:vector_splitting}
\end{equation}
Moreover, let $\gradpar$ and $\gradperp$ represent the horizontal and vertical components of the gradient operator of a scalar field $\rho$
\begin{equation}
\gradpar \rho := \frac{1}{a\cos(\phi)}\frac{\partial\rho}{\partial\lambda} \ex + \frac{1}{a}\frac{\partial\rho}{\partial\phi} \ey\,, \qquad 
\gradperp\rho := \frac{\partial\rho}{\partial z} \ez\,. \label{eq:grad_splitting}
\end{equation}
In a similar way, $\curlpar$ and $\curlperp$ represent, respectively, the horizontal and vertical components of the curl of a vector field $\boldsymbol{u} = u\ex + v\ey + w\ez$
\begin{equation}
\curlpar \boldsymbol{u} := \left(\frac{1}{a}\frac{\partial w}{\partial\phi} - \frac{\partial v}{\partial z}\right)\ex + \left(\frac{\partial u}{\partial z} - \frac{1}{a\cos(\phi)}\frac{\partial w}{\partial\lambda}\right)\ey\,, \qquad  
\curlperp \boldsymbol{u} := \frac{1}{a\cos(\phi)}\left(\frac{\partial v}{\partial\lambda} - \frac{\partial(\cos(\phi) u)}{\partial\phi}\right)\ez\,. \label{eq:curl_splitting}
\end{equation}
Note that we have assumed the shallow atmosphere approximation of constant radius, $a$, in the above expressions. 
In practice, these spherical transformations are absorbed into the definition of the Jacobian \cite{Guba14} (including the
shallow atmosphere approximation) and the 
associated Piola transformations for the different function spaces \cite{LP18,LP20}.

From \eqref{eq:grad_splitting} and \eqref{eq:curl_splitting} it follows directly that
\[
\nabla\rho = \gradpar\rho + \gradperp\rho\,, \quad\mathrm{and}\quad \nabla\times\boldsymbol{u} = \curlpar\boldsymbol{u} +\curlperp\boldsymbol{u}\,.
\]
We may then reformulate the vorticity, $\boldsymbol{\omega} := \nabla\times\boldsymbol{u}$, as
\begin{equation}
\boldsymbol{\omega} = \underbrace{\overbrace{\curlpar\upar}^{\omegaparpar}+ \overbrace{\curlpar\uperp}^{\omegaparperp}}_{\omegapar} + \underbrace{\curlperp\upar}_{\omegaperp}\,. \label{eq:vorticity_splitting}
\end{equation}
While the spatial discretisation outlined above is valid for basis functions of any polynomial degree, $p$, in practice
we use polynomials of arbitrary order in the horizontal only, and use a lowest order discretisation in the vertical for
which the $\mathcal{P}_h$ bases are piecewise linear and the $\mathcal{Q}_h$ bases are piecewise constant \cite{LP20}.

Using \eqref{eq:vector_splitting}, \eqref{eq:grad_splitting}, \eqref{eq:curl_splitting}, and \eqref{eq:vorticity_splitting}, 
we may split the compressible Euler equations, \eqref{eq:compressible_euler}, into horizontal and vertical components
\begin{subequations}\label{eq:compressible_euler_splitting}
\begin{align}
\frac{\partial\upar}{\partial t} +(\omegaperp + \boldsymbol{f}_{\perp})\times\upar + \omegaparpar\times\uperp + \frac{1}{2}\gradpar\|\upar\|^{2} + \theta\gradpar\Pi &= 0, \label{eq:mom_par_split}\\
\frac{\partial\uperp}{\partial t} + \omegaparperp\times\upar + \gradperp\left(\frac{1}{2}\|\uperp\|^{2} + gz\right) + \theta\gradperp\Pi &= 0, \label{eq:mom_perp_split} \\
\frac{\partial\rho}{\partial t} + \nabla\cdot\left(\rho\upar\right) + \nabla\cdot\left(\rho\uperp\right) &= 0, \label{eq:dens_split} \\
\frac{\partial (\rho\theta)}{\partial t} + \nabla\cdot\left(\rho\theta\upar\right) + \nabla\cdot\left(\rho\theta\uperp\right) &= 0, \label{eq:temp_split}
\end{align}
\end{subequations}
where, as in \eqref{eq:vorticity_splitting}, $\omegaparpar := \curlpar\upar$ and $\omegaparperp := \curlpar\uperp$.
The same splitting into horizontal and vertical components may be applied for the basis functions, \eqref{eq:function_spaces_span},
\begin{equation}
\wbasis_{j} = \left(\wbasis_{j}\right)_{\parallel} +  \left(\wbasis_{j}\right)_{\perp} := \wbasispar_{j} + \wbasisperp_{j}, \label{eq:wbasis_split}
\end{equation}
and
\begin{equation}
\ubasis_{j} = \left(\ubasis_{j}\right)_{\parallel} +  \left(\ubasis_{j}\right)_{\perp} := \ubasispar_{j} + \ubasisperp_{j}. \label{eq:ubasis_split}
\end{equation}
Note that by splitting the vector spaces at the level of the basis functions and excluding horizontal-vertical cross terms
from the Jacobian \cite{LP20} the current dimensional splitting does not permit variations in topography. In order to incorporate 
topography the cross terms will need to be introduced into the Jacobian, such that the horizontal forcings project onto the vertical
momentum equation (and vice versa) as source terms. These source terms would not be energetically balanced according to the HEVI
splitting described below.

Multiplying each of the equations in \eqref{eq:compressible_euler_splitting} by the appropriate finite dimensional test functions, 
and discretising the solution variables accordingly, the horizontal discrete equations are given as \cite{LP20}
\begin{subequations}\label{eq:euler_wf_discrete_algebraic_par_split_incidence}
			\begin{align}
	\begin{split}
\innerprod{\ubasispar_{i}}{\ubasispar_{j}}_{\Omega}\frac{\mathrm{d}\mathsf{u}_{i,\parallel}}{\mathrm{d} t} + \innerprod{(\boldsymbol{\omega}_{h,\perp} + \boldsymbol{f}_{h,\perp})\times\ubasispar_{i}}{\ubasispar_{j}}_{\Omega}\mathsf{u}_{i,\parallel} + & \\
\innerprod{\boldsymbol{\omega}_{h,\parallel,\parallel}\times\ubasisperp_{i}}{\ubasispar_{j}}_{\Omega}\mathsf{u}_{i,\perp} - & \\
\left(\mathsf{E}^{3,2}_{\parallel}\right)_{j,k}^{\top}\innerprod{\frac{1}{2}\boldsymbol{u}_{h,\parallel} \cdot \ubasispar_{i}}{\qbasis_{k}}_{\Omega}\mathsf{u}_{i,\parallel}  + 
\innerprod{\theta_{h} \ubasispar_{i}}{\ubasispar_{j}}_{\Omega}\mathsf{P}_{i,\parallel} &= 0, \qquad \forall \ubasispar_{j}\in\uspacepar,
	\end{split} \label{eq::mom_weak_discrete_algebraic_par_split_incidence}\\
-\left(\mathsf{E}^{3,2}_{\parallel}\right)_{j,k}^{\top}\innerprod{\qbasis_{i}}{\qbasis_{k}}_{\Omega} \mathsf{\Pi}_{i} - \innerprod{\ubasispar_{i}}{\ubasispar_{j}}_{\Omega}\mathsf{P}_{i,\parallel} &= 0, \qquad\forall \ubasispar_{j}\in\uspacepar, \label{eq:Pi_P_discrete_algebraic_par_split_incidence}\\
\innerprod{\rho_{h}\ubasispar_{i}}{\ubasispar_{j}}_{\Omega}\mathsf{u}_{i,\parallel} - \innerprod{\ubasispar_{i}}{\ubasispar_{j}}_{\Omega}\mathsf{U}_{i,\parallel} &= 0, \qquad \forall \ubasispar_{j}\in\uspacepar, \label{eq:rho_u_V_discrete_algebraic_par_split_incidence}\\
\innerprod{\theta_{h}\ubasispar_{i}}{\ubasispar_{j}}_{\Omega}\mathsf{U}_{i,\parallel} - \innerprod{\ubasispar_{i}}{\ubasispar_{j}}_{\Omega}\mathsf{F}_{i,\parallel} &= 0, \qquad \forall \ubasispar_{j}\in\uspacepar, \label{eq:theta_V_F_discrete_algebraic_par_split_incidence} \\
-\left(\mathsf{E}^{2,1}_{\parallel,\parallel}\right)_{j,k}^{\top}\innerprod{\ubasispar_{i}}{\ubasispar_{k}}_{\Omega}\mathsf{u}_{i,\parallel} - \innerprod{\wbasispar_{i}}{\wbasispar_{j}}_{\Omega}\mathsf{\omega}_{i,\parallel,\parallel} & = 0, 
\qquad\forall\wbasispar_{j}\in\wspacepar, \label{eq:curl_u_omega_discrete_algebraic_par_par_split_incidence} \\
-\left(\mathsf{E}^{2,1}_{\parallel,\perp}\right)_{j,k}^{\top}\innerprod{\ubasispar_{i}}{\ubasispar_{k}}_{\Omega}\mathsf{u}_{i,\parallel} - \innerprod{\wbasisperp_{i}}{\wbasisperp_{j}}_{\Omega}\mathsf{\omega}_{i,\perp} & = 0, 
\qquad\forall\wbasisperp_{j}\in\wspaceperp, \label{eq:curl_u_omega_discrete_algebraic_perp_split_incidence}
			\end{align}
		\end{subequations}
where we have introduced $\omega_{h,\parallel,\parallel} := \sum_{i=0}^{\wcardinalpar-1}\mathsf{\omega}_{i,\parallel,\parallel} \wbasispar_{i}$, $\omega_{h,\parallel,\perp} := \sum_{i=0}^{\wcardinalpar-1} \mathsf{\omega}_{i,\parallel,\perp}\wbasispar_{i}$. Note that we have introduced additional diagnostic equations for the pressure gradients \eqref{eq:Pi_P_discrete_algebraic_par_split_incidence}, 
mass and temperature fluxes, \eqref{eq:rho_u_V_discrete_algebraic_par_split_incidence}, \eqref{eq:theta_V_F_discrete_algebraic_par_split_incidence} and
vorticity terms, \eqref{eq:curl_u_omega_discrete_algebraic_par_par_split_incidence}, \eqref{eq:curl_u_omega_discrete_algebraic_perp_split_incidence}.
In the same way, the vertical discrete equations are
\begin{subequations}\label{eq:euler_wf_discrete_algebraic_perp_split_incidence}
\begin{align}
\begin{split}
\innerprod{\ubasisperp_{i}}{\ubasisperp_{j}}_{\Omega}\frac{\mathrm{d}\mathsf{u}_{i,\perp}}{\mathrm{d} t} + \innerprod{\boldsymbol{\omega}_{h,\parallel,\perp}\times\ubasispar_{i}}{\ubasisperp_{j}}_{\Omega}\mathsf{u}_{i,\parallel} - & \\
\left(\mathsf{E}^{3,2}_{\perp}\right)_{j,k}^{\top}\innerprod{\frac{1}{2}{\boldsymbol{u}}_{h,\perp}\cdot\ubasisperp_{i}}{\qbasis_{k}}_{\Omega}\mathsf{u}_{i,\perp} - &\\
g\left(\mathsf{E}^{3,2}_{\perp}\right)_{j,k}^{\top}\innerprod{\qbasis_i}{\qbasis_{k}}_{\Omega}\mathsf{z}_{i} +
\innerprod{\theta_{h} \ubasisperp_{i}}{\ubasisperp_{j}}_{\Omega}\mathsf{P}_{i,\perp} &= 0, \qquad\forall\ubasisperp_{j}\in\uspaceperp,
					\end{split} \label{eq::mom_weak_discrete_algebraic_perp_split_incidence}\\
-\left(\mathsf{E}^{3,2}_{\perp}\right)_{j,k}^{\top}\innerprod{\qbasis_{i}}{\qbasis_{j}}_{\Omega} \mathsf{\Pi}_{i} - \innerprod{\ubasisperp_{i}}{\ubasisperp_{j}}_{\Omega}\mathsf{P}_{i,\perp} &= 0, \qquad\forall\ubasisperp_{j}\in\uspaceperp,\label{eq:Pi_P_discrete_algebraic_perp_split_incidence}\\
\innerprod{\rho_{h}\ubasisperp_{i}}{\ubasisperp_{j}}_{\Omega}\mathsf{u}_{i,\perp} - \innerprod{\ubasisperp_{i}}{\ubasisperp_{j}}_{\Omega}\mathsf{U}_{i,\perp} &= 0, \qquad\forall\ubasisperp_{j}\in\uspaceperp, \label{eq:rho_u_V_discrete_algebraic_perp_split_incidence}\\
\innerprod{\theta_{h}\ubasisperp_{i}}{\ubasisperp_{j}}_{\Omega}\mathsf{U}_{i,\perp} - \innerprod{\ubasisperp_{i}}{\ubasisperp_{j}}_{\Omega}\mathsf{F}_{i,\perp} &= 0, \qquad\forall\ubasisperp_{j}\in\uspaceperp, \label{eq:theta_V_F_discrete_algebraic_perp_split_incidence} \\
-\left(\mathsf{E}^{2,1}_{\perp}\right)_{j,k}^{\top}\innerprod{\ubasisperp_{i}}{\ubasisperp_{k}}_{\Omega}\mathsf{u}_{i,\perp} - \innerprod{\wbasispar_{i}}{\wbasispar_{j}}_{\Omega}\mathsf{\omega}_{i,\parallel,\perp} & = 0, \qquad\forall\wbasispar_{j}\in\wspacepar, \label{eq:curl_u_omega_discrete_algebraic_par_perp_split_incidence}
			\end{align}
		\end{subequations}
where we have introduced $\omega_{h,\perp} := \sum_{i=0}^{\wcardinalperp-1}\mathsf{\omega}_{i,\perp} \wbasisperp_{i}$.

Additionally, we also have the flux form equations for density and density weighted potential temperature transport that
contain both vertical and horizontal components. While we have not included these in the split systems described in
\eqref{eq:euler_wf_discrete_algebraic_par_split_incidence} and
\eqref{eq:euler_wf_discrete_algebraic_perp_split_incidence}, since doing so incurs a temporal splitting
error, in practice these equations are also split between their horizontal and vertical components. 
These equations are given as
\begin{subequations}
\label{eq:euler_wf_discrete_algebraic_no_par_no_perp_split_incidence}
\begin{align}
\begin{split}
\innerprod{\qbasis_{i}}{\qbasis_{j}}_{\Omega}\frac{\mathrm{d}\mathsf{\rho}_{i}}{\mathrm{d}t} + \innerprod{\qbasis_{k}}{\qbasis_{j}}_{\Omega}\left(\mathsf{E}^{3,2}_{\parallel}\right)_{k,i}\mathsf{U}_{i,\parallel} + 
\innerprod{\qbasis_{k}}{\qbasis_{j}}_{\Omega}\left(\mathsf{E}^{3,2}_{\perp}\right)_{k,i}\mathsf{U}_{i,\perp} &= 0, \qquad\forall\qbasis_{j}\in\qspace,
\end{split} \label{eq::dens_weak_discrete_algebraic_split_incidence} \\
\begin{split}
\innerprod{\qbasis_{i}}{\qbasis_{j}}_{\Omega}\frac{\mathrm{d}\mathsf{\Theta}_{i}}{\mathrm{d}t} + \innerprod{\qbasis_{k}}{\qbasis_{j}}_{\Omega}\left(\mathsf{E}^{3,2}_{\parallel}\right)_{k,i}\mathsf{F}_{i,\parallel} +
\innerprod{\qbasis_{k}}{\qbasis_{j}}_{\Omega}\left(\mathsf{E}^{3,2}_{\perp}\right)_{k,i}\mathsf{F}_{i,\perp} &= 0, \qquad\forall\qbasis_{j}\in\qspace,
\end{split} \label{eq::temp_weak_discrete_algebraic_split_incidence}
\end{align}
\end{subequations}
We also have two additional diagnostic equations for the potential temperature and Exner pressure. These are
given respectively as
\begin{equation}
\label{eq:rho_theta_Theta_discrete_algebraic} 
\innerprod{\qbasis_{i}}{\epsilon^{\mathcal{U}_{\perp}}_{j}}_{\Omega}\mathsf{\Theta}_{i}
- \innerprod{\rho_{h}\epsilon^{\mathcal{U}_{\perp}}_{i}}{\epsilon^{\mathcal{U}_{\perp}}_{j}}_{\Omega}\mathsf{\theta}_{i}
= 0, \qquad\forall\epsilon^{\mathcal{U}_{\perp}}_{j}\in\uspaceperp,
\end{equation}
\begin{equation}
\label{eq:Pi_Theta_discrete_algebraic}
c_p\Big(\frac{R}{p_0}\Big)^{R/c_v}\innerprod{(\qbasis_{i}\mathsf{\Theta}_{i})^{R/c_v}}{\qbasis_{j}}_{\Omega} -
\innerprod{\qbasis_{i}}{\qbasis_{j}}_{\Omega}\mathsf{\Pi}_{i} = 0, \qquad\forall\qbasis_{j}\in\qspace.
\end{equation}

For the sake of brevity, we use compact matrix notation for the remainder of this article. Using this notation,
the horizontal system, \eqref{eq:euler_wf_discrete_algebraic_par_split_incidence} is expressed as
\begin{subequations}\label{eq:euler_wf_discrete_algebraic_par_split_incidence_matrix}
			\begin{align}
\massupar\frac{\mathrm{d}\ualgebraicpar}{\mathrm{d}t} + \boldsymbol{\mathsf{R}}^{\parallel,\parallel}\ualgebraicpar +
\boldsymbol{\mathsf{R}}^{\parallel, \perp}\ualgebraicperp -
\left(\boldsymbol{\mathsf{E}}^{3,2}_{\parallel}\right)^{\top}\boldsymbol{\mathsf{T}}^{\mathcal{U}_{\parallel}}\ualgebraicpar +
\boldsymbol{\mathsf{S}}^{\mathcal{U}_{\parallel}}\boldsymbol{\mathsf{P}}^{\parallel} &= \boldsymbol{\mathsf{0}},\label{eq::mom_weak_discrete_algebraic_par_split_incidence_matrix}\\
-\left(\boldsymbol{\mathsf{E}}^{3,2}_{\parallel}\right)^{\top}\massq\boldsymbol{\mathsf{\Pi}} - \massupar\boldsymbol{\mathsf{P}}^{\parallel} &= \boldsymbol{\mathsf{0}}, \label{eq:Pi_P_discrete_algebraic_par_split_incidence_matrix}\\
				\nupar\ualgebraicpar - \massupar\boldsymbol{\mathsf{U}}^{\parallel} & = \boldsymbol{\mathsf{0}}, \label{eq:rho_u_V_discrete_algebraic_par_split_incidence_matrix}\\
\boldsymbol{\mathsf{S}}^{\mathcal{U}_{\parallel}}\boldsymbol{\mathsf{U}}^{\parallel} - \massupar\boldsymbol{\mathsf{F}}^{\parallel} &= \boldsymbol{\mathsf{0}}, \label{eq:theta_V_F_discrete_algebraic_par_split_incidence_matrix}\\
-\left(\boldsymbol{\mathsf{E}}^{2,1}_{\parallel,\parallel}\right)^{\top}\massupar\ualgebraicpar - \masswpar\boldsymbol{\mathsf{\omega}}^{\parallel,\parallel} &= \boldsymbol{\mathsf{0}}, \label{eq:curl_u_omega_discrete_algebraic_par_split_incidence_matrix}\\
-\left(\boldsymbol{\mathsf{E}}^{2,1}_{\perp}\right)^{\top}\massuperp\ualgebraicperp - \masswpar\boldsymbol{\mathsf{\omega}}^{\parallel,\perp} &= \boldsymbol{\mathsf{0}}.
			\end{align}
		\end{subequations}
In a similar fashion, the discrete vertical equations, \eqref{eq:euler_wf_discrete_algebraic_perp_split_incidence}, 
are expressed as
\begin{subequations}\label{eq:euler_wf_discrete_algebraic_perp_split_incidence_matrix}
\begin{align}
\massuperp\frac{\mathrm{d}\ualgebraicperp}{\mathrm{d}t} + \boldsymbol{\mathsf{R}}^{\perp,\parallel}\ualgebraicpar -
\left(\boldsymbol{\mathsf{E}}^{3,2}_{\perp}\right)^{\top}\boldsymbol{\mathsf{T}}^{\mathcal{U}_{\perp}}\ualgebraicperp -
g\left(\boldsymbol{\mathsf{E}}^{3,2}_{\perp}\right)^{\top}\massq\boldsymbol{\mathsf{z}} +
\boldsymbol{\mathsf{S}}^{\mathcal{U}_{\perp}}\boldsymbol{\mathsf{P}}^{\perp} &= \boldsymbol{\mathsf{0}},\label{eq::mom_weak_discrete_algebraic_perp_split_incidence_matrix}\\
-\left(\boldsymbol{\mathsf{E}}_{\perp}^{3,2}\right)^{\top}\massq\boldsymbol{\mathsf{\Pi}} - \massuperp\boldsymbol{\mathsf{P}}^{\perp} &= \boldsymbol{\mathsf{0}},\label{eq:Pi_P_discrete_algebraic_perp_split_incidence_matrix}\\
\nuperp\ualgebraicperp - \massuperp\boldsymbol{\mathsf{U}}^{\perp} &= \boldsymbol{\mathsf{0}},\label{eq:rho_u_V_discrete_algebraic_perp_split_incidence_matrix}\\
\boldsymbol{\mathsf{S}}^{\mathcal{U}_{\perp}}\boldsymbol{\mathsf{U}}^{\perp} - \massuperp\boldsymbol{\mathsf{F}}^{\perp} &= \boldsymbol{\mathsf{0}}, \label{eq:theta_V_F_discrete_algebraic_perp_split_incidence_matrix} \\
-\left(\boldsymbol{\mathsf{E}}^{2,1}_{\parallel,\perp}\right)^{\top}\massuperp\ualgebraicperp - \masswperp\boldsymbol{\mathsf{\omega}}^{\perp} &= \boldsymbol{\mathsf{0}}.\label{eq:curl_u_omega_discrete_algebraic_perp_split_incidence_matrix}
\end{align}
\end{subequations}
Finally, 
\eqref{eq:euler_wf_discrete_algebraic_no_par_no_perp_split_incidence}-\eqref{eq:Pi_Theta_discrete_algebraic} 
may be written in compact matrix notation as
\begin{subequations}\label{eq:euler_wf_discrete_algebraic_no_par_no_perp_split_incidence_matrix}
\begin{align}
\massq\frac{\mathrm{d}\boldsymbol{\mathsf{\rho}}}{\mathrm{d}t} + \massq\boldsymbol{\mathsf{E}}^{3,2}_{\parallel}\boldsymbol{\mathsf{U}}^{\parallel} + \massq\boldsymbol{\mathsf{E}}^{3,2}_{\perp}\boldsymbol{\mathsf{U}}^{\perp} &= \boldsymbol{\mathsf{0}},\label{eq::dens_weak_discrete_algebraic_split_incidence_matrix}\\
\massq\frac{\mathrm{d}\boldsymbol{\mathsf{\Theta}}}{\mathrm{d}t} + \massq\boldsymbol{\mathsf{E}}^{3,2}_{\parallel}\boldsymbol{\mathsf{F}}^{\parallel} + \massq\boldsymbol{\mathsf{E}}^{3,2}_{\perp}\boldsymbol{\mathsf{F}}^{\perp} &= \boldsymbol{\mathsf{0}},\label{eq::temp_weak_discrete_algebraic_split_incidence_matrix} \\
\luq\boldsymbol{\mathsf{\Theta}}
- \nuperp\boldsymbol{\mathsf{\theta}}
&= \boldsymbol{\mathsf{0}},\label{eq:rho_theta_Theta_discrete_algebraic_split_incidence_matrix}\\
c_p\Big(\frac{R}{p_0}\Big)^{R/c_v}\innerprod{(\qbasis_{i}\mathsf{\Theta}_{i})^{R/c_v}}{\qbasis_{j}}_{\Omega}
- \massq \boldsymbol{\mathsf{\Pi}} &= \boldsymbol{\mathsf{0}}.\label{eq:discrete_eos}
\end{align}
\end{subequations}

\section{Discrete energetics}

In \cite{LP20} the energetic properties of the discrete form of the compressible Euler equations were analysed
for a spatial discretisation using mixed mimetic spectral elements. Here that analysis is extended to 
incorporate the temporal discretisation.
 
The discrete Hamiltonian $\mathcal{H}_{h} := \mathcal{H}[\boldsymbol{u}_{h}, \rho_{h}, \Theta_{h}]$ is given as \cite{LP20}
\begin{equation}\label{eq::hamiltonian}
	\mathcal{H}[\boldsymbol{u}_{h}, \rho_{h}, \Theta_{h}] = \int_{\Omega}\frac{1}{2}\rho_{h}\|\boldsymbol{u}_{h}\|^{2}\,\mathrm{d}\Omega + \int_{\Omega}\rho_{h}gz_{h}\,\mathrm{d}\Omega +  \int_{\Omega} c_{v} \left(\frac{R}{p_{0}}\right)^{\frac{R}{c_{v}}}\Theta_{h}^{\frac{c_{p}}{c_{v}}}\,\mathrm{d}\Omega\,.
\end{equation}
Invoking the definition of the variational derivative \cite{Celledoni12}, these are given for the 
Hamiltonian \eqref{eq::hamiltonian} as
\begin{subequations}\label{eq::variational_derivs}
\begin{align}
\innerprod{\frac{\delta\mathcal{H}}{\delta\boldsymbol{u}_{h}}}{\ubasis_{j}}_{\Omega} &=
\innerprod{\rho_{h}\boldsymbol{u}_{h}}{\ubasis_{j}}_{\Omega}=
\innerprod{\boldsymbol{U}_{h}}{\ubasis_{j}}_{\Omega}, \qquad\forall\ubasis_{j}\in\uspace,\\
\innerprod{\frac{\delta\mathcal{H}}{\delta\rho_{h}}}{\qbasis_{j}}_{\Omega} &=
\innerprod{\frac{1}{2}\|\boldsymbol{u}_{h}\|^{2} + gz_h}{\qbasis_{j}}_{\Omega} =
\innerprod{\Phi_{h}}{\qbasis_{j}}_{\Omega}, \qquad\forall\qbasis_{j}\in\qspace,\\
\innerprod{\frac{\delta\mathcal{H}}{\delta\Theta_{h}}}{\qbasis_{j}}_{\Omega} &=
c_p\Bigg(\frac{R}{p_0}\Bigg)^{R/c_v}\innerprod{\Theta_h^{R/c_v}}{\qbasis_{j}}_{\Omega}=
\innerprod{\Pi_h}{\qbasis_{j}}_{\Omega}
,\qquad\forall\qbasis_{j}\in\qspace.
\label{eq::var_div_pi}
\end{align}
\end{subequations}
The semi-discrete form of the compressible Euler equations may then be formulated as a skew-symmetric system as
\begin{equation}\label{eq:euler_eqns_discrete_skew_symmetric}
\left[
\begin{array}{c}
	\massu\boldsymbol{\mathsf{u}}_{,t} \\
	\massq\mathsf{\rho}_{,t} \\
	\massq\mathsf{\Theta}_{,t}
\end{array}
\right] =
\begin{bmatrix}
-\boldsymbol{\mathsf{R}}_q & \left(\boldsymbol{\mathsf{E}}^{3,2}\right)^{\top}\massq &
\boldsymbol{\mathsf{S}}^{\mathcal{U}}\left(\massu\right)^{-1} \left(\boldsymbol{\mathsf{E}}^{3,2}\right)^{\top}\massq\\
-\massq\boldsymbol{\mathsf{E}}^{3,2} & \boldsymbol{\mathsf{0}} & \boldsymbol{\mathsf{0}} \\
-\massq\boldsymbol{\mathsf{E}}^{3,2}\left(\massu\right)^{-1} \boldsymbol{\mathsf{S}}^{\mathcal{U}} &
\boldsymbol{\mathsf{0}} & \boldsymbol{\mathsf{0}}\\
\end{bmatrix}
\left[
\begin{array}{c}
	\boldsymbol{\mathsf{U}} \\
	\mathsf{\Phi} \\
	\mathsf{\Pi}
\end{array}
\right].
\end{equation}
This formulation is consistent with both the horizontal 
\eqref{eq:euler_wf_discrete_algebraic_par_split_incidence_matrix} and vertical 
\eqref{eq:euler_wf_discrete_algebraic_perp_split_incidence_matrix} discretisations
detailed above.
The fully discrete, second order in time \cite{Cohen11} form of 
\eqref{eq:euler_eqns_discrete_skew_symmetric} is then given as
\begin{equation}\label{eq:euler_eqns_discrete_skew_symmetric_2}
\left[
\begin{array}{c}
\massu\boldsymbol{\mathsf{u}}^{n+1} \\
\massq\mathsf{\rho}^{n+1} \\
\massq\mathsf{\Theta}^{n+1}
\end{array}
\right] =
\left[
\begin{array}{c}
\massu\boldsymbol{\mathsf{u}}^{n} \\
\massq\mathsf{\rho}^{n} \\
\massq\mathsf{\Theta}^{n}
\end{array}
\right] + 
\Delta t
\begin{bmatrix}
-\widehat{\boldsymbol{\mathsf{R}}}_q & \left(\boldsymbol{\mathsf{E}}^{3,2}\right)^{\top}\massq &
\widehat{\boldsymbol{\mathsf{S}}^{\mathcal{U}}}\left(\massu\right)^{-1} \left(\boldsymbol{\mathsf{E}}^{3,2}\right)^{\top}\massq\\
-\massq\boldsymbol{\mathsf{E}}^{3,2} & \boldsymbol{\mathsf{0}} & \boldsymbol{\mathsf{0}} \\
-\massq\boldsymbol{\mathsf{E}}^{3,2}\left(\massu\right)^{-1} \widehat{\boldsymbol{\mathsf{S}}^{\mathcal{U}}} &
\boldsymbol{\mathsf{0}} & \boldsymbol{\mathsf{0}}\\
\end{bmatrix}
\left[
\begin{array}{c}
	\overline{\boldsymbol{\mathsf{U}}} \\
	\overline{\mathsf{\Phi}} \\
	\overline{\mathsf{\Pi}}
\end{array}
\right],
\end{equation}
where $\widehat{\boldsymbol{\mathsf{R}}}_q$ and $\widehat{\boldsymbol{\mathsf{S}}^{\mathcal{U}}}$ 
are time centered operators and 
$\overline{\boldsymbol{\mathsf{U}}}$, $\overline{\mathsf{\Phi}}$ and $\overline{\mathsf{\Pi}}$ 
are second order in time versions of the variational derivatives of the energy.
These are given as \cite{Bauer18,Wimmer20}
\begin{subequations}\label{eq::variational_derivs_discrete}
\begin{align}
\massuperp\overline{\boldsymbol{\mathsf{U}}^{\perp}} &= 
\frac{1}{3}\Nuperp{n}\Ualgebraicperp{n} + 
\frac{1}{6}\Nuperp{n}\Ualgebraicperp{k} + 
\frac{1}{6}\Nuperp{k}\Ualgebraicperp{n} + 
\frac{1}{3}\Nuperp{k}\Ualgebraicperp{k} \\
\massq\overline{\boldsymbol{\mathsf{\Phi}}} &=
\frac{1}{3}\boldsymbol{\mathsf{T}}^{\mathcal{U}_{\perp},n}\Ualgebraicperp{n} + 
\frac{1}{3}\boldsymbol{\mathsf{T}}^{\mathcal{U}_{\perp},k}\Ualgebraicperp{n} + 
\frac{1}{3}\boldsymbol{\mathsf{T}}^{\mathcal{U}_{\perp},k}\Ualgebraicperp{k} + g\massq\boldsymbol{\mathsf{z}} \\
\massq\overline{\boldsymbol{\mathsf{\Pi}}} &= \frac{1}{2}\massq\boldsymbol{\mathsf{\Pi}}^n + \frac{1}{2}\massq\boldsymbol{\mathsf{\Pi}}^k.
\end{align}
\end{subequations}
Note that these variational derivatives are exactly integrated in time, as can be seen by approximating
the dependent variables by piecewise linear bases for second order temporal accuracy as
$[\boldsymbol{\mathsf{u}}^{\perp},\boldsymbol{\rho},\boldsymbol{\mathsf{\Pi}}] \approx 
[\Ualgebraicperp{n},\boldsymbol{\rho}^n,\boldsymbol{\mathsf{\Pi}}^n](1-s)+[\Ualgebraicperp{k},\boldsymbol{\rho}^k,\boldsymbol{\mathsf{\Pi}}^k]s$,
and integrating each of the above expressions between dimensionless times $s=0$ and $s=1$.
This ensures that the chain rule as given in \eqref{eq::dHda_dadt_2} is preserved in the 
discrete form. We further note that unlike the temporal discretisation the spatial integrals 
are inexact due to both the inclusion of transcendental functions within the metric terms \cite{LP18,LP20}, 
and the non-integer polynomial within the equation of state \eqref{eq::var_div_pi}. 
It is perhaps possible that this use of inexact integration may lead to a loss of energy 
conservation due to aliasing errors. However any loss of energy conservation is likely to
be extremely small, since in each case the under-integrated functions are slowly varying.
This is borne out by anecdotal experience, since in all cases the linear systems converge
to a small tolerance ($1.0^{-16}$) with no difficulty, and as shown in subsequent sections
the energy conservation errors are small when only the vertical dynamics are considered.

Multiplying both sides of \eqref{eq:euler_eqns_discrete_skew_symmetric_2} by 
$\begin{bmatrix}\overline{\boldsymbol{\mathsf{U}}}^{\top} & \overline{\mathsf{\Phi}}^{\top} & \overline{\mathsf{\Pi}}^{\top}\end{bmatrix}$, 
gives
\begin{equation}\label{eq::dHda_dadt}
\overline{\boldsymbol{\mathsf{U}}}^{\top}\massu(\boldsymbol{\mathsf{u}}^{n+1} - \boldsymbol{\mathsf{u}}^{n}) +
\overline{\boldsymbol{\mathsf{u}}}^{\top}(\boldsymbol{\mathsf{T}}^{\mathcal{U}})^{\top}
(\boldsymbol{\mathsf{\rho}}^{n+1} - \boldsymbol{\mathsf{\rho}}^{n}) + 
g\boldsymbol{\mathsf{z}}^{\top}\massq(\boldsymbol{\mathsf{\rho}}^{n+1} - \boldsymbol{\mathsf{\rho}}^{n}) +
\overline{\boldsymbol{\mathsf{\Pi}}}^{\top}\massq(\boldsymbol{\mathsf{\Theta}}^{n+1} - \boldsymbol{\mathsf{\Theta}}^{n}) = 0,
\end{equation}
which is the discrete form of the continuous relation
\begin{equation}\label{eq::dHda_dadt_2}
\frac{\delta K}{\delta\boldsymbol{u}}\cdot\frac{\partial\boldsymbol{u}}{\partial t} + 
\frac{\delta K}{\delta \rho}\cdot\frac{\partial \rho}{\partial t} + 
\frac{\delta P}{\delta \rho}\cdot\frac{\partial \rho}{\partial t} + 
\frac{\delta I}{\delta \Theta}\cdot\frac{\partial \Theta}{\partial t} = 0.
\end{equation}
Note that the left hand side row vectors in \eqref{eq::dHda_dadt} are integrated exactly for the 
second order temporal approximation. This ensures that the chain rule as expressed in \eqref{eq::dHda_dadt_2}
is preserved in the discrete form. 
Note also that for
$\mathsf{R}_{q,ij} := \innerprod{\boldsymbol{q}_{h}\times\ubasis_{j}}{\ubasis_{i}}_{\Omega}$,
where $\boldsymbol{q}_h$ is the potential vorticity \cite{LPG18},
this is itself a skew-symmetric operator such that
$\boldsymbol{\mathsf{U}}^{\top}\boldsymbol{\mathsf{R}}_q\boldsymbol{\mathsf{U}} =
\boldsymbol{\mathsf{U}}^{\top}\boldsymbol{\mathsf{R}}\boldsymbol{\mathsf{u}} = \boldsymbol{\mathsf{0}}$.
As such neither $\boldsymbol{\mathsf{R}}_q$ nor $\boldsymbol{\mathsf{R}}$ projects onto the
energy in the discrete form.

Equation \eqref{eq::dHda_dadt} is the discrete equivalent of \eqref{eq::conservation}, and ensures
that total energy is conserved between time levels $n$ and $n+1$. This is achieved by the balanced
exchanges of kinetic, potential and internal energy, which are given as

\begin{align}
\frac{\partial K_h}{\partial t} &= g\overline{\boldsymbol{\mathsf{U}}}^{\top}(\boldsymbol{\mathsf{E}}^{3,2})^{\top}\massq\boldsymbol{\mathsf{z}} +
\overline{\boldsymbol{\mathsf{U}}}^{\top}\widehat{\boldsymbol{\mathsf{S}}^{\mathcal{U}}}\left(\massu\right)^{-1}
(\boldsymbol{\mathsf{E}}^{3,2})^{\top}\massq\overline{\boldsymbol{\mathsf{\Pi}}},\label{eq::dKdt}\\
\frac{\partial P_h}{\partial t} &= -g\boldsymbol{\mathsf{z}}^{\top}\massq\boldsymbol{\mathsf{E}}^{3,2}\overline{\boldsymbol{\mathsf{U}}},\label{eq::dPdt}\\
\frac{\partial I_h}{\partial t} &= -\overline{\boldsymbol{\mathsf{\Pi}}}^{\top}\massq\boldsymbol{\mathsf{E}}^{3,2}
\left(\massu\right)^{-1}\widehat{\boldsymbol{\mathsf{S}}^{\mathcal{U}}}\overline{\boldsymbol{\mathsf{U}}}.\label{eq::dIdt}
\end{align}
The right hand side terms of \eqref{eq::dKdt} exactly balance those of \eqref{eq::dPdt} and \eqref{eq::dIdt},
thus allowing for the exact balances of kinetic to potential and kinetic to internal energy respectively.

\section{Implicit vertical solver}

The horizontal and vertical spatial discretisations described above 
\eqref{eq:euler_wf_discrete_algebraic_par_split_incidence_matrix}-\eqref{eq:euler_wf_discrete_algebraic_no_par_no_perp_split_incidence_matrix}
are time split using a horizontally explicit/vertically implicit (HEVI) scheme.
While implicit vertical solvers are traditionally used in order to time step over the CFL 
limit of the vertical acoustic modes, here we have the additional motivation of conserving 
energetic exchanges associated with vertical atmospheric processes.
Various flavors of HEVI schemes have previously been employed in non-hydrostatic atmospheric models
\cite{Ullrich12,Giraldo13,Bao15,Gardner18}.
Here we base our splitting scheme on the horizontally third order, vertically second order,
trapazoidal TRAP(2,3,2) scheme \cite{Weller13,Lock14}, which has an overall accuracy of second order in time. 
The TRAP(2,3,2) scheme involves two time centered, Crank-Nicolson, vertical substeps. These substeps
are unconditionally stable, implying that that the energy conservation errors remain bounded. The energetically
balanced integrator is in some regards a natural extension of such a scheme, involving additional cross terms
\eqref{eq::variational_derivs_discrete} so as to preserve the chain rule as given in \eqref{eq::dHda_dadt_2} 
in the discrete form, as well as the 
skew-symmetric property of the Poisson bracket \eqref{eq:euler_eqns_discrete_skew_symmetric}.

Given the horizontal and vertical state vectors 
$\boldsymbol{a}^{\parallel}=[\boldsymbol{u}^{\parallel}\quad\rho\quad\Theta]^{\top}$ and
$\boldsymbol{a}^{\perp}=[\boldsymbol{u}^{\perp}\quad\rho\quad\Theta]^{\top}$, the 
full state vector $\boldsymbol{a}=[\boldsymbol{u}^{\parallel}\quad\boldsymbol{u}^{\perp}\quad\rho\quad\Theta]^{\top}$, 
and the explicit horizontal tendencies, $H(\boldsymbol{a^{\parallel}})$
this scheme integrates the
equations of motion over a time step $\Delta t$ between time levels $n$ and $n+1$ as:
\begin{subequations}\label{eq::hevi}
\begin{align}
\boldsymbol{a}^{\parallel,1} &= \boldsymbol{a}^{\parallel,n} - \Delta t H(\boldsymbol{a}^{\parallel,n})\\
\boldsymbol{a}^{2} + \Delta tV(\boldsymbol{a}^{\perp,n},\boldsymbol{a}^{\perp,2}) &= \boldsymbol{a}^{n} -
\frac{\Delta t}{2}H(\boldsymbol{a}^{\parallel,n}) - \frac{\Delta t}{2}H(\boldsymbol{a}^{\parallel,1})\\
\boldsymbol{a}^{n+1} + \Delta tV(\boldsymbol{a}^{\perp,n},\boldsymbol{a}^{\perp,n+1}) &= \boldsymbol{a}^{n} -
\frac{\Delta t}{2}H(\boldsymbol{a}^{\parallel,n}) - \frac{\Delta t}{2}H(\boldsymbol{a}^{\parallel,2}),
\end{align}
\end{subequations}
where $V(\boldsymbol{a}^{\perp,n},\boldsymbol{a}^{\perp,n+1})$ represents the action of
the skew-symmetric operator on the variational derivatives of the energy as given in
\eqref{eq:euler_eqns_discrete_skew_symmetric_2}.
For a linear system, there are no cross terms in the variational derivatives, and the scheme will 
revert to the original TRAP(2,3,2) formulation and exhibit the same linear stability properties, which are
detailed for acoustic modes in \cite{Lock14}. As for 
the Crank-Nicholson substeps in the TRAP(2,3,2) scheme, the energetically balanced vertical integrator is 
unconditionally stable, by virtue of the fact that it conserves energy.

\subsection{Implicit vertical solve}\label{sec:vertical_implicit_half_step}

The errors at a given Newton iteration for the implicit vertical solver in \eqref{eq::hevi} at
nonlinear iteration, $k$ may be expressed as residuals, drawing from
\eqref{eq:euler_wf_discrete_algebraic_perp_split_incidence_matrix} and
\eqref{eq:euler_wf_discrete_algebraic_no_par_no_perp_split_incidence_matrix} as
\begin{subequations}\label{eq::residuals}
\begin{align}
\boldsymbol{F}_u = &
\massuperp\Ualgebraicperp{k} - \massuperp\Ualgebraicperp{n} + 
\Delta t\widehat{\boldsymbol{\mathsf{R}}^{\perp,\parallel}}\overline{\ualgebraicpar} +
\Delta t\left(\boldsymbol{\mathsf{E}}^{3,2}_{\perp}\right)^{\top}\massq\overline{\boldsymbol{\mathsf{\Phi}}} +
\Delta t\widehat{\boldsymbol{\mathsf{S}}^{\mathcal{U}_{\perp}}}\left(\massuperp\right)^{-1}
\left(\boldsymbol{\mathsf{E}}_{\perp}^{3,2}\right)^{\top}\massq\overline{\boldsymbol{\mathsf{\Pi}}} \\
\boldsymbol{F}_{\rho} = &
\massq\boldsymbol{\mathsf{\rho}}^{k} + \massq\boldsymbol{\mathsf{\rho}}^{n} + 
\Delta t\massq\boldsymbol{\mathsf{E}}^{3,2}_{\perp}\overline{\boldsymbol{\mathsf{U}}^{\perp}} \\
\boldsymbol{F}_{\Theta} = &
\massq\boldsymbol{\mathsf{\Theta}}^{k} + \massq\boldsymbol{\mathsf{\Theta}}^{n} + 
\Delta t\massq\boldsymbol{\mathsf{E}}^{3,2}_{\perp}\left(\massuperp\right)^{-1}
\widehat{\boldsymbol{\mathsf{S}}^{\mathcal{U}_{\perp}}}\overline{\boldsymbol{\mathsf{U}}^{\perp}}
\\
\boldsymbol{F}_{\Pi} = &
\innerprod{\log(\qbasis_{i}\mathsf{\Pi}_i^k) - 
\log\left(\qbasis_{i}\frac{R}{c_v}\mathsf{\Theta}_{i}^k\right) - \qbasis_{i}\log(c_p) - 
\qbasis_{i}\left(\frac{R}{c_v}\log\left(\frac{R}{p_0}\right)\right)}{\qbasis_{j}}_{\Omega},\label{eq::F_pi}
\end{align}
\end{subequations}
where \eqref{eq::F_pi} is the natural logarithm of the discrete equation of state \eqref{eq:discrete_eos}.

At each Newton iteration, the updates to the solutions are computed from the approximate Jacobian, given as a
$4\times 4$ block system as
\begin{equation}\label{eq::jacobian}
\begin{bmatrix}
\boldsymbol{\mathsf{M}}_u & \boldsymbol{\mathsf{0}} & \boldsymbol{\mathsf{G}}_{\Theta} & \boldsymbol{\mathsf{G}}_{\Pi} \\
\boldsymbol{\mathsf{D}}_{\rho}  & \boldsymbol{\mathsf{M}}_{\rho} & \boldsymbol{\mathsf{0}}   & \boldsymbol{\mathsf{0}} \\
\boldsymbol{\mathsf{D}}_{\Theta} & \boldsymbol{\mathsf{Q}}_{\Theta,\rho} & \boldsymbol{\mathsf{M}}_{\Theta} & \boldsymbol{\mathsf{0}} \\
\boldsymbol{\mathsf{0}}   & \boldsymbol{\mathsf{0}}   & \boldsymbol{\mathsf{C}}_{\Theta} & \boldsymbol{\mathsf{C}}_{\Pi} \\
\end{bmatrix} 
\begin{bmatrix} 
\delta\ualgebraicperp \\ \delta\boldsymbol{\rho} \\ \delta\boldsymbol{\mathsf{\Theta}} \\ \delta\boldsymbol{\mathsf{\Pi}} \\ 
\end{bmatrix} \approx -
\begin{bmatrix} \boldsymbol{F}_{u} \\ \boldsymbol{F}_{\rho} \\ \boldsymbol{F}_{\Theta} \\ \boldsymbol{F}_{\Pi} \\ \end{bmatrix},
\end{equation}
where the blocks are defined as
\begin{subequations}\label{eq::jacobian_blocks}
\begin{align}
\boldsymbol{\mathsf{M}}_u &= \massuperp,\\
\boldsymbol{\mathsf{G}}_{\Theta} &= 
\frac{\Delta t}{2}\innerprod{\qbasis_j}{\ubasisperp_i\cdot\widehat{\mathsf{P}}_h}_{\Omega}\innerprod{\widehat{\rho}_h\qbasis_j}{\qbasis_l}_{\Omega}^{-1}\massq,\\
\boldsymbol{\mathsf{G}}_{\Pi} &= \frac{\Delta t}{2}\widehat{\boldsymbol{\mathsf{S}}^{\mathcal{U}_{\perp}}}\left(\massuperp\right)^{-1}
\left(\boldsymbol{\mathsf{E}}_{\perp}^{3,2}\right)^{\top}\massq,\\
\boldsymbol{\mathsf{D}}_{\rho} &= 
\frac{\Delta t}{2}\massq\boldsymbol{\mathsf{E}}^{3,2}_{\perp}\left(\massuperp\right)^{-1}\widehat{\nuperp},\\
\boldsymbol{\mathsf{M}}_{\rho} &= \massq,\\
\boldsymbol{\mathsf{D}}_{\Theta} &= \frac{\Delta t}{2}\innerprod{\widehat{\Theta}_h\qbasis_j}{\qbasis_i}_{\Omega}\boldsymbol{\mathsf{E}}_{\perp}^{3,2},\\
\boldsymbol{\mathsf{Q}}_{\Theta,\rho} &= 
\frac{\Delta t}{2}\innerprod{\widehat{u_{\perp h}}\cdot\ubasisperp_j}{\qbasis_i}_{\Omega}
\left(\massuperp\right)^{-1}
\left(\boldsymbol{\mathsf{E}}^{3,2}_{\perp}\right)^{\top}
\innerprod{\widehat{\theta}_h\qbasis_m}{\qbasis_l}_{\Omega},\\
\boldsymbol{\mathsf{M}}_{\Theta} &= \massq,\\
\boldsymbol{\mathsf{C}}_{\Theta} &= -\frac{R}{c_v}\massq\innerprod{\Theta_h\qbasis_j}{\qbasis_i}_{\Omega}^{-1}\massq,\label{eq::C_Theta}\\
\boldsymbol{\mathsf{C}}_{\Pi} &= \massq\innerprod{\Pi_h\qbasis_j}{\qbasis_i}_{\Omega}^{-1}\massq.
\end{align}
\end{subequations}
Note that unlike the original system of residual vectors \eqref{eq::residuals}, in order to preserve 
energetic balances this approximate Jacobian need not be itself a skew-symmetric system.
All that is required is that this operator facilitate the convergence of such a system.
The operators $\boldsymbol{\mathsf{G}}_{\Theta}$ and $\boldsymbol{\mathsf{Q}}_{\Theta,\rho}$ are of 
particular significance, since as will be seen these allow for the Helmholtz structure of the reduced
system. The operator $\boldsymbol{\mathsf{Q}}_{\Theta,\rho}$ has the structure of the material advection
term, which accounts for the fact that while the density advection equation \eqref{eq::dens} is strictly
a flux form equation, the density weighted potential temperature equation \eqref{eq::temp} may 
equivalently be expressed as a material advection equation for the potential temperature 
\cite{Natale16,Melvin19,WCB20}.

In order to reduce \eqref{eq::jacobian} into a more computationally tractable form we begin by 
eliminating the density update as
\begin{equation}\label{eq::delta_rho}
\delta\boldsymbol{\rho} = -\boldsymbol{\mathsf{M}}_{\rho}^{-1}(\boldsymbol{F}_{\rho} + 
\boldsymbol{\mathsf{D}}_{\rho}\delta\boldsymbol{\mathsf{u}}^{\perp})
\end{equation}
Such that
\begin{equation}
\boldsymbol{\mathsf{D}}_{\Theta}\delta\boldsymbol{\mathsf{u}}^{\perp} - \boldsymbol{\mathsf{Q}}_{\Theta,\rho}\boldsymbol{\mathsf{M}}_{\rho}^{-1}
(\boldsymbol{F}_{\rho} + \boldsymbol{\mathsf{D}}_{\rho}\delta\boldsymbol{u}) + 
\boldsymbol{\mathsf{M}}_{\Theta}\delta\boldsymbol{\mathsf{\Theta}} = -\boldsymbol{F}_{\Theta}
\end{equation}

\begin{equation}
\begin{bmatrix}
\boldsymbol{\mathsf{M}}_u & \boldsymbol{\mathsf{G}}_{\Theta} & \boldsymbol{\mathsf{G}}_{\Pi} \\
\boldsymbol{\mathsf{D}}_{\Theta} -  
\boldsymbol{\mathsf{Q}}_{\Theta,\rho}\boldsymbol{\mathsf{M}}_{\rho}^{-1}\boldsymbol{\mathsf{D}}_{\rho}
& \boldsymbol{\mathsf{M}}_{\Theta} & \boldsymbol{\mathsf{0}} \\
\boldsymbol{\mathsf{0}}   & \boldsymbol{\mathsf{C}}_{\Theta} & \boldsymbol{\mathsf{C}}_{\Pi} \\
\end{bmatrix} 
\begin{bmatrix} \delta\boldsymbol{\mathsf{u}}^{\perp} \\ \delta\boldsymbol{\mathsf{\Theta}} \\ \delta\boldsymbol{\mathsf{\Pi}} \\ \end{bmatrix} = 
-\begin{bmatrix} \boldsymbol{F}_{u} \\ 
\boldsymbol{F}_{\Theta} - \boldsymbol{\mathsf{Q}}_{\Theta,\rho}\boldsymbol{\mathsf{M}}_{\rho}^{-1}\boldsymbol{F}_{\rho} 
\\ \boldsymbol{F}_{\Pi} \\ \end{bmatrix}
\end{equation}
Eliminating the Exner pressure update then gives
\begin{equation}\label{eq::delta_pi}
\delta\boldsymbol{\mathsf{\Pi}} = -\boldsymbol{\mathsf{C}}_{\Pi}^{-1}(\boldsymbol{F}_{\Pi} + 
\boldsymbol{\mathsf{C}}_{\Theta}\delta\boldsymbol{\mathsf{\Theta}})
\end{equation}
Such that
\begin{equation}\label{eq::delta_u}
\boldsymbol{\mathsf{M}}_u\delta\boldsymbol{\mathsf{u}}^{\perp} + \boldsymbol{\mathsf{G}}_{\Theta}\delta\boldsymbol{\mathsf{\Theta}} - 
\boldsymbol{\mathsf{G}}_{\Pi}\boldsymbol{\mathsf{C}}_{\Pi}^{-1}(\boldsymbol{F}_{\Pi} + 
\boldsymbol{\mathsf{C}}_{\Theta}\delta\boldsymbol{\mathsf{\Theta}}) =
-\boldsymbol{F}_{u}
\end{equation}

\begin{equation}
\begin{bmatrix}
\boldsymbol{\mathsf{M}}_u & \boldsymbol{\mathsf{G}}_{\Theta} - 
\boldsymbol{\mathsf{G}}_{\Pi}\boldsymbol{\mathsf{C}}_{\Pi}^{-1}\boldsymbol{\mathsf{C}}_{\Theta} \\
\boldsymbol{\mathsf{D}}_{\Theta} -  
\boldsymbol{\mathsf{Q}}_{\Theta,\rho}\boldsymbol{\mathsf{M}}_{\rho}^{-1}\boldsymbol{\mathsf{D}}_{\rho} &
\boldsymbol{\mathsf{M}}_{\Theta} \\
\end{bmatrix} 
\begin{bmatrix} \delta\boldsymbol{\mathsf{u}}^{\perp} \\ \delta\boldsymbol{\mathsf{\Theta}} \\ \end{bmatrix} = 
-\begin{bmatrix} \boldsymbol{F}_{u}  - \boldsymbol{\mathsf{G}}_{\Pi}\boldsymbol{\mathsf{C}}_{\Pi}^{-1}\boldsymbol{F}_{\Pi} \\ 
\boldsymbol{F}_{\Theta} - \boldsymbol{\mathsf{Q}}_{\Theta,\rho}\boldsymbol{\mathsf{M}}_{\rho}^{-1}\boldsymbol{F}_{\rho} \\ \end{bmatrix}
\end{equation}
Finally, Schur complement reduction yields a single equation for the density weighted potential temperature
update as
\begin{multline}\label{eq::delta_rt}
\begin{bmatrix}
\boldsymbol{\mathsf{M}}_{\Theta} - (\boldsymbol{\mathsf{D}}_{\Theta} -  
\boldsymbol{\mathsf{Q}}_{\Theta,\rho}\boldsymbol{\mathsf{M}}_{\rho}^{-1}\boldsymbol{\mathsf{D}}_{\rho})\boldsymbol{\mathsf{M}}_u^{-1}
(\boldsymbol{\mathsf{G}}_{\Theta} - 
\boldsymbol{\mathsf{G}}_{\Pi}\boldsymbol{\mathsf{C}}_{\Pi}^{-1}\boldsymbol{\mathsf{C}}_{\Theta})
\end{bmatrix}\delta\boldsymbol{\mathsf{\Theta}} = \\
-\begin{bmatrix}
\boldsymbol{F}_{\Theta} - \boldsymbol{\mathsf{Q}}_{\Theta,\rho}\boldsymbol{\mathsf{M}}_{\rho}^{-1}\boldsymbol{F}_{\rho} -
(\boldsymbol{\mathsf{D}}_{\Theta} -  
\boldsymbol{\mathsf{Q}}_{\Theta,\rho}\boldsymbol{\mathsf{M}}_{\rho}^{-1}\boldsymbol{\mathsf{D}}_{\rho})\boldsymbol{\mathsf{M}}_u^{-1}
(\boldsymbol{F}_{u}  - \boldsymbol{\mathsf{G}}_{\Pi}\boldsymbol{\mathsf{C}}_{\Pi}^{-1}\boldsymbol{F}_{\Pi})\end{bmatrix}
\end{multline}
Once \eqref{eq::delta_rt} has been solved for the update to the density weighted potential temperature
at nonlinear iteration $k$, the density, Exner pressure and velocity may be updated by \eqref{eq::delta_rho},
\eqref{eq::delta_pi} and \eqref{eq::delta_u} respectively. Equation \eqref{eq::delta_rt} has the structure 
of a Helmholtz equation, for which both the divergence, $\boldsymbol{\mathsf{D}}_{\Theta}$ and gradient, 
$\boldsymbol{\mathsf{G}}_{\Theta}$ terms are \emph{corrected} by secondary terms due to the additional
equations and dynamics. While the gradient correction will always be of the same sign as the original 
gradient due to the sign of \eqref{eq::C_Theta}, it is at least
theoretically possible that $\boldsymbol{\mathsf{Q}}_{\Theta,\rho}$ could change sign. This could potentially
result in multiple extrema and a catastrophic loss of convergence for the nonlinear system as a whole. 
However this has caused no problems in the tests performed here, 
and we are not sure of any physical scenarios in which this would be reversed.

While preconditioning of compressible systems is typically applied so as to solve an inner linear 
equation for the pressure update \cite{Melvin19}, update equations for for density weighted potential 
temperature have previously been derived for high resolution atmospheric models using Jacobian-Free
Newton-Krylov methods \cite{Reisner05}.

Once the updates have been determined, the solutions at Newton iteration $k+1$ are given as
\begin{subequations}
\begin{align}
\Ualgebraicperp{k+1} &= \Ualgebraicperp{k} + \delta\boldsymbol{\mathsf{u}}^{\perp}\\
\boldsymbol\rho^{k+1} &= \boldsymbol\rho^{k} + \delta\boldsymbol\rho\\
\boldsymbol{\mathsf{\Theta}}^{k+1} &= \boldsymbol{\mathsf{\Theta}}^{k} + \delta\boldsymbol{\mathsf{\Theta}}\\
\boldsymbol{\mathsf{\Pi}}^{k+1} &= \boldsymbol{\mathsf{\Pi}}^{k} + \delta\boldsymbol{\mathsf{\Pi}}.
\end{align}
\end{subequations}
Convergence is achieved once $|\delta\boldsymbol{\mathsf{u}}^{\perp}|/|\Ualgebraicperp{k+1}|$, 
$|\delta\boldsymbol\rho|/|\boldsymbol\rho^{k+1}|$,
$|\delta\boldsymbol{\mathsf{\Theta}}|/|\boldsymbol{\mathsf{\Theta}}^{k+1}|$ and 
$|\delta\boldsymbol{\mathsf{\Pi}}|/|\boldsymbol{\mathsf{\Pi}}^{k+1}|$ are all less that a specified tolerance
in all vertical columns (here we use $10.0^{-8}$ as the tolerance for all solution variables unless otherwise 
stated).
Since the vertical kinetic energy is extremely small in comparison to both the potential
and internal energies, (which depend on $\rho$ and $\Theta$ respectively), 
and vertical velocities may approach zero in some regions, it is perhaps
advisable to remove $|\delta\boldsymbol{\mathsf{u}}^{\perp}|/|\Ualgebraicperp{k+1}|$ 
from the convergence criteria. This has been done for the third test case presented here
(3D rising bubble), since the bubble is initially at rest and hydrostatic balance is exact,
so the initial variations in vertical velocity are negligible, and also for conservation tests 
involving only the vertical dynamics.

Note that while the implicit temporal discretisation described above is cast in the context of 
the vertical dynamics, since the elements are high order in the horizontal, each vertical solve
actually applies to $p^2$ degrees of freedom at each vertical level, where $p$ is the polynomial
degree of the horizontal component of the basis functions $\qbasis_i\in\qspace$.

\subsection{Implementation details}

The vertical integrator described in the preceding section was implemented using the PETSc 
linear algebra library \cite{petsc-user-ref,petsc-web-page,petsc-efficient}. Since the solution 
variables for the vertical solve are all discontinuous across horizontal element boundaries, each
horizontal element may be solved independently, which allows for dramatically increased 
computational performance. All matrix inverses detailed in the preceding section may then be 
computed directly using LU decomposition, including the final solve for the density weighted
potential temperature update \eqref{eq::delta_rt}. The horizontal terms are solved in
parallel using GMRES with block Jacobi preconditioning \cite{LP20}. While the PETSc linear solver
is used for the inner vertical solves, the outer Newton iteration is implemented directly without
use of the PETSc Newton solver, or any line search augmentation.

In order to stabilise the model we add a biharmonic viscosity to both the horizontal momentum 
equation \cite{LP18}, and the horizontal density weighted potential temperature equation,
using a value of $0.072\Delta x^{3.2}$ \cite{Guba14}, where $\Delta x$ is the average spacing between spectral 
element nodes. 
These are added as explicit tendencies to the right hand sides of 
\eqref{eq::mom_weak_discrete_algebraic_par_split_incidence_matrix} and 
\eqref{eq::temp_weak_discrete_algebraic_split_incidence_matrix} respectively.
A Rayleigh damping term \cite{Klemp08} has also been applied to the vertical momentum equation in the 
top three levels, with descending values of
$4.0/\Delta t, 2.0/\Delta t, 1.0/\Delta t$, so as to damped oscillations associated with the initial
hydrostatic adjustment process. Note that this term is applied only to the first test 
case described below.

\section{Results}

\subsection{Baroclinic instability}

The model is validated using a $z$-level dry baroclinic instability test case \cite{Ullrich14} with 
the shallow atmosphere approximation. The initial state is one of geostrophic horizontal and hydrostatic
vertical balance, overlaid with a small, $\mathcal{O}(1\mathrm{m/s})$, perturbation to the zonal and
meridional velocity components. The model was run with $24\times 24$ elements of degree $p=3$ on each 
face of the cubed sphere (and piecewise constant/linear elements in the vertical), for an averaged resolution of 
$\Delta x\approx 128\mathrm{km}$ and 30 vertical levels on 96 processors with a time step of 
$\Delta t = 120\mathrm{s}$.

\begin{figure}[!hbtp]
\begin{center}
\includegraphics[width=0.48\textwidth,height=0.36\textwidth]{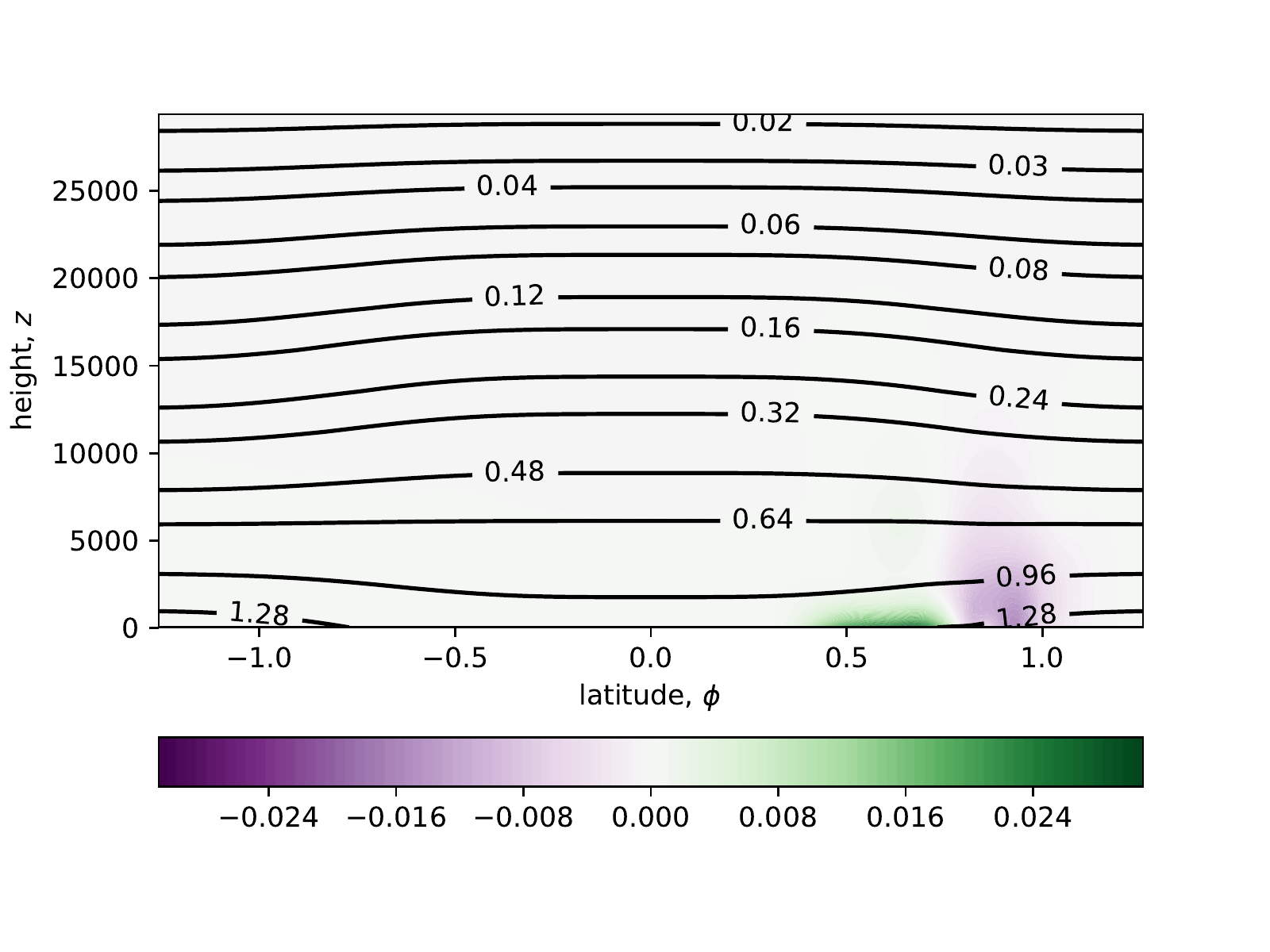}
\includegraphics[width=0.48\textwidth,height=0.36\textwidth]{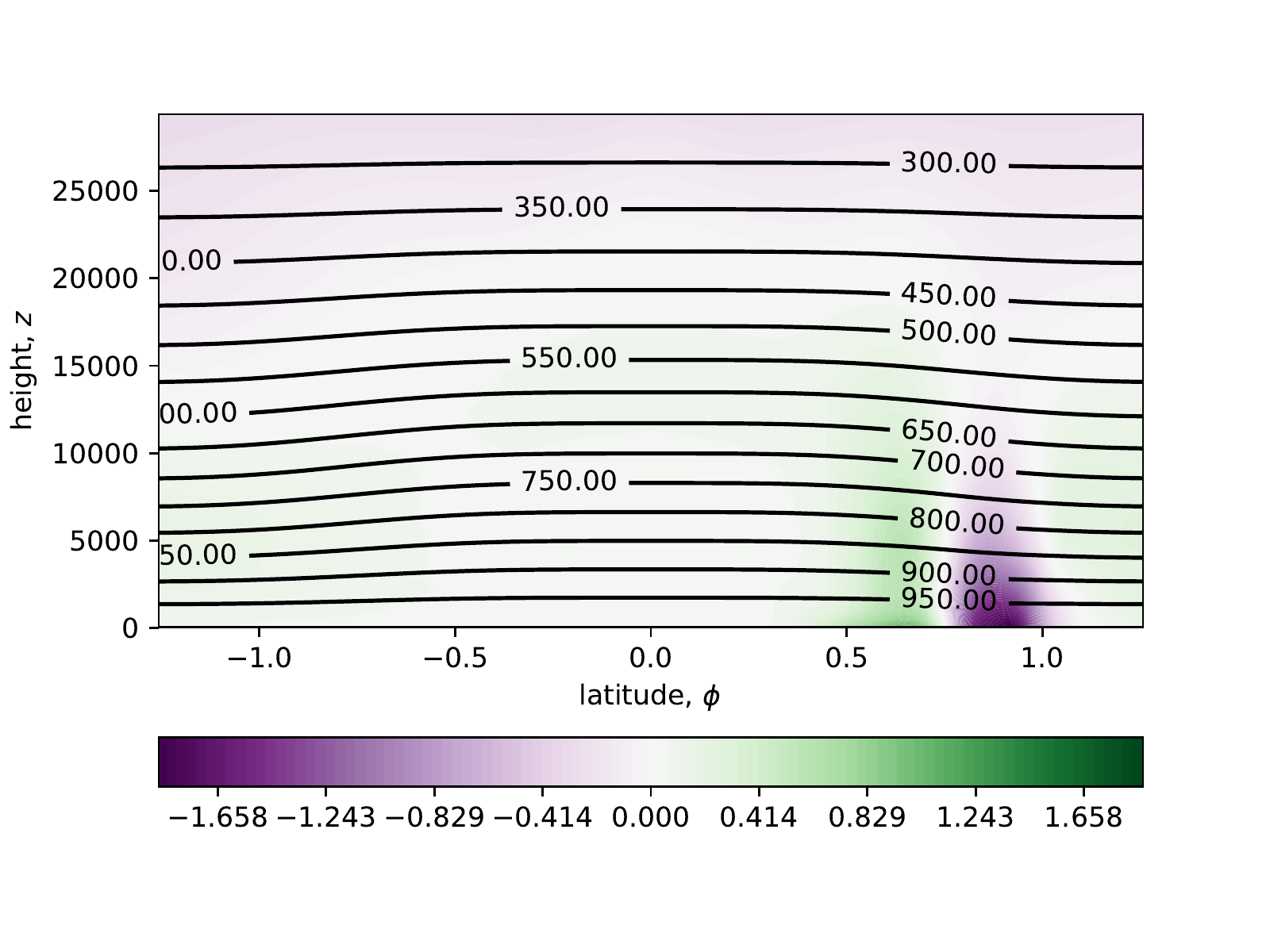}
\caption{
Baroclinic instability test case: 
Zonal averages of density, $\rho_h$ in $\mathrm{kg}\cdot\mathrm{m}^{-3}$ (left) and
Exner pressure, $\Pi_h$ in $\mathrm{m}^{2}\mathrm{s}^{-2}\mathrm{K}^{-1}$ (right) at day 10.
Contours represent absolute values, and shades represent differences from initial values.}
\label{fig::zonal_avg_1}
\end{center}
\end{figure}

\begin{figure}[!hbtp]
\begin{center}
\includegraphics[width=0.48\textwidth,height=0.36\textwidth]{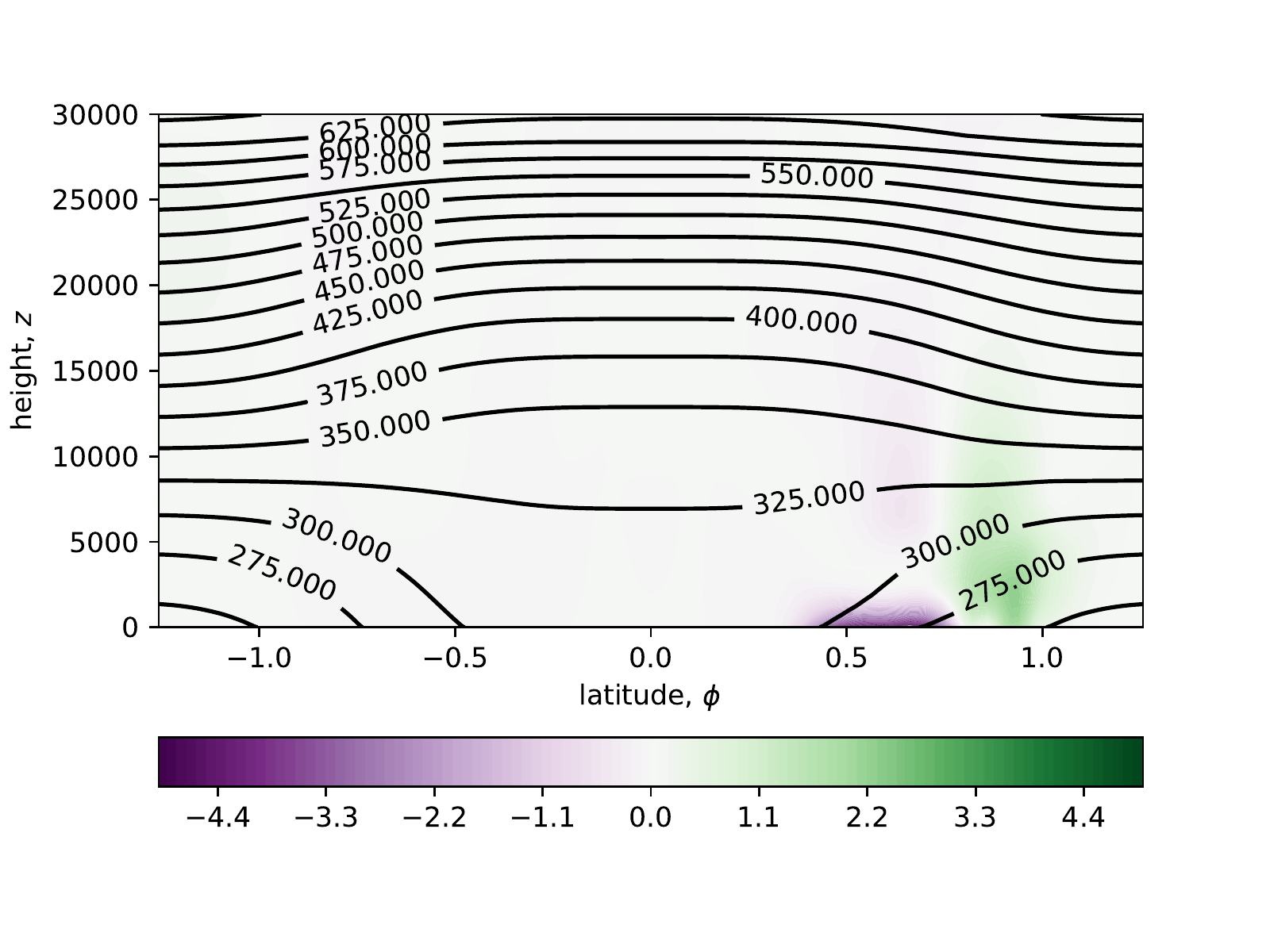}
\includegraphics[width=0.48\textwidth,height=0.36\textwidth]{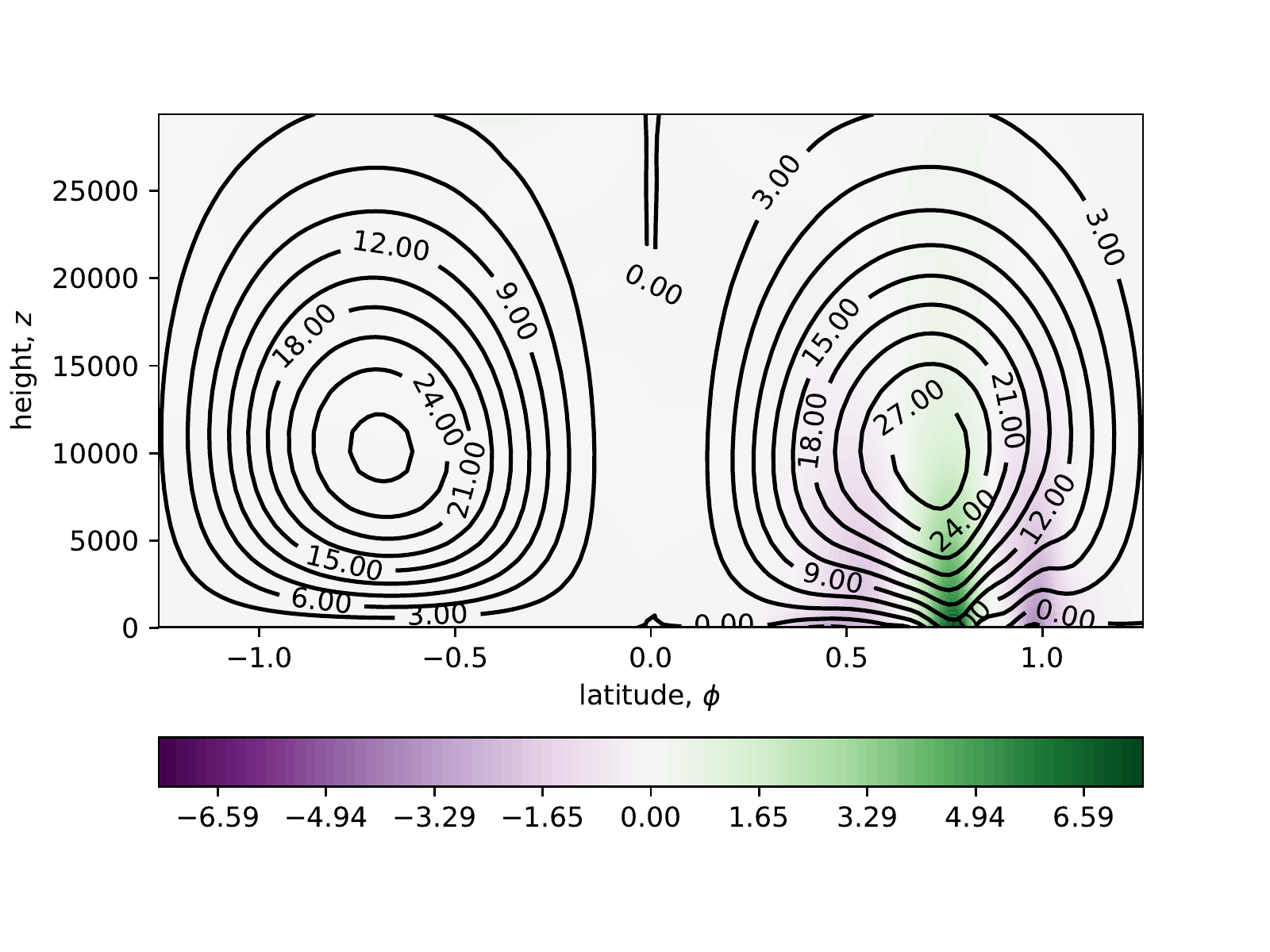}
\caption{
Baroclinic instability test case: 
Zonal averages of potential temperature, $\theta_h$ in $\mathrm{K}$ (left) and
zonal velocity, $u_h$ in $\mathrm{m}\cdot\mathrm{s}^{-1}$ (right) at day 10.
Contours represent absolute values, and shades represent differences from initial values.}
\label{fig::zonal_avg_2}
\end{center}
\end{figure}

Figures \ref{fig::zonal_avg_1} and \ref{fig::zonal_avg_2} show the zonal averages of density $\rho$,
Exner pressure $\Pi$, potential temperature $\theta$ and zonal velocity $u$ at day 10
(solid lines), as well as the differences between the final and initial states. These profiles
show little difference between the initial and final states, with the exception of the zonal velocity,
which exhibits a small kink near the bottom boundary where the baroclinic instability occurs,
demonstrating that the leading order geostrophic and hydrostatic balances in the horizontal and vertical
are well satisfied. 

\begin{figure}[!hbtp]
\begin{center}
\includegraphics[width=0.48\textwidth,height=0.36\textwidth]{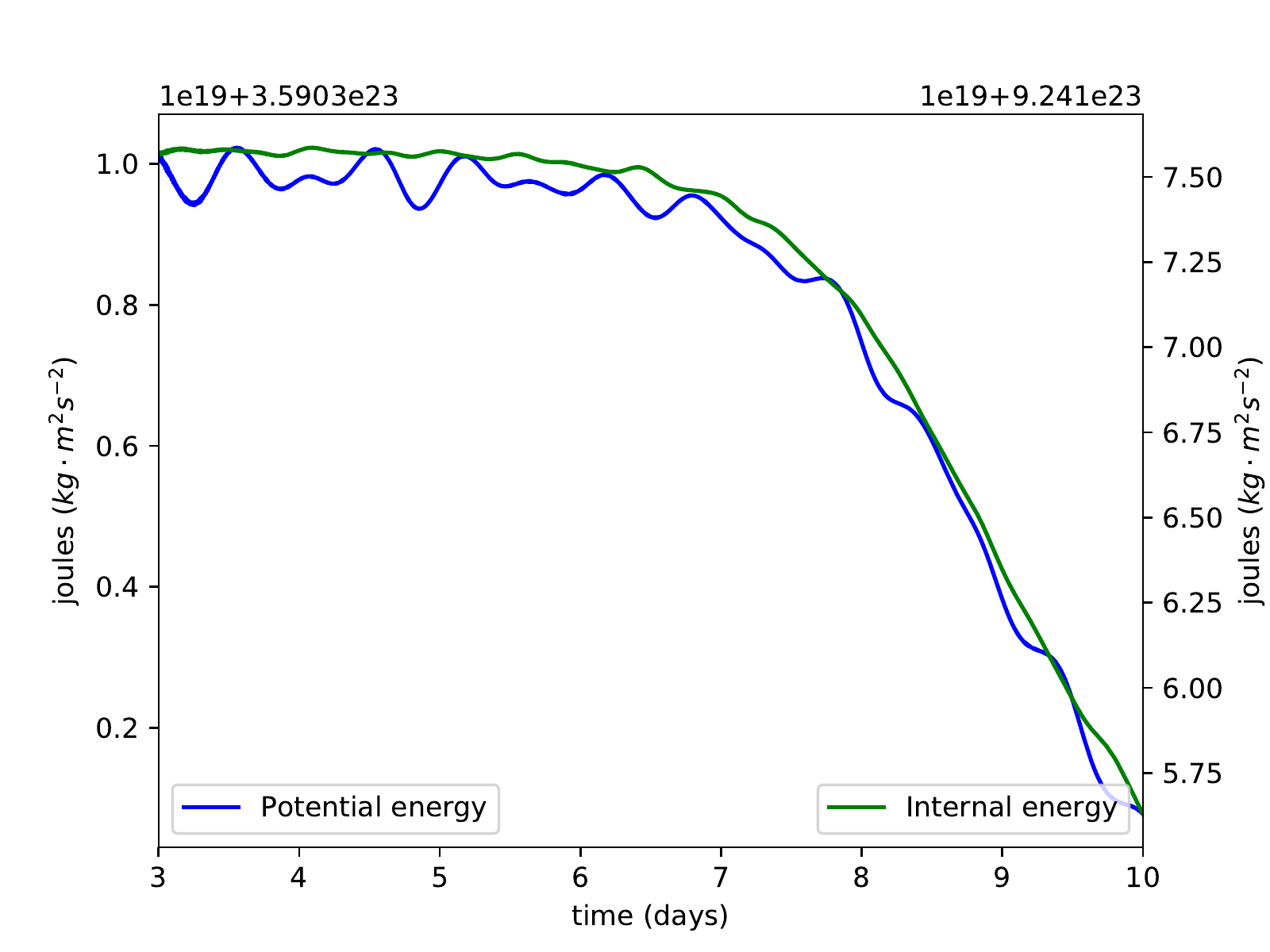}
\includegraphics[width=0.48\textwidth,height=0.36\textwidth]{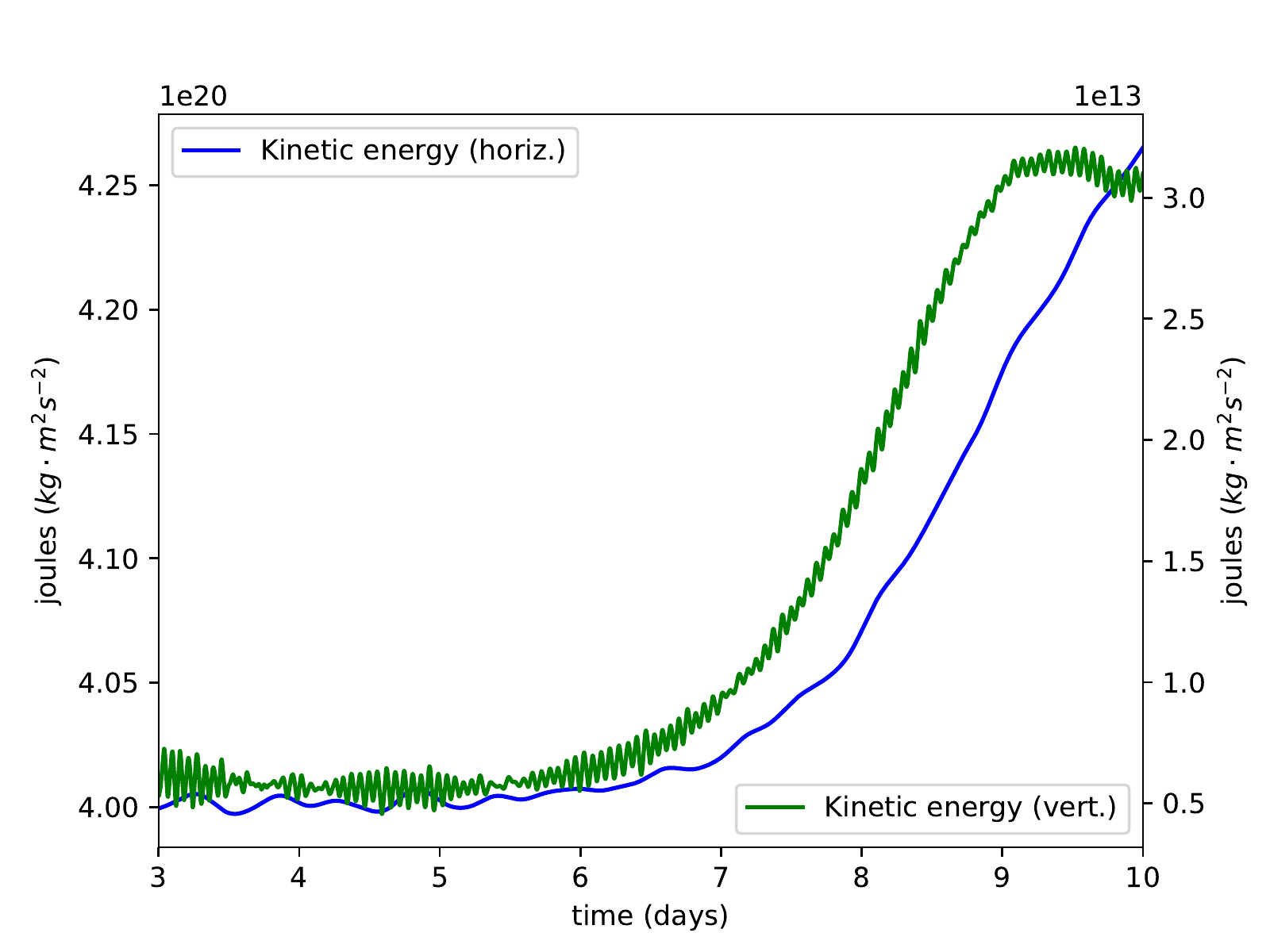}
\caption{
Baroclinic instability test case: 
Left: time evolution of the potential and internal energies. Right: time evolution of
the kinetic energy (horizontal and vertical).}
\label{fig::energetics_1}
\end{center}
\end{figure}

Figure \ref{fig::energetics_1} shows the evolution of the
potential, internal and kinetic (horizontal and vertical) energy with time. As can be seen the
advent of the baroclinic instability coincides with a loss of both potential and internal 
energy, and an increase of kinetic energy due to acceleration in both the horizontal and vertical
directions. The signature of the internal gravity waves is also visible in the smaller scale
oscillation of the potential and horizontal kinetic energies. 
Note that the total amounts of potential and internal energy are approximately
$3.6\times 10^{23}$ and $9.2\times 10^{23}\mathrm{kg\cdot m^2s^{-2}}$ respectively, and so are
several orders of magnitude greater than the amounts of horizontal and vertical kinetic energy
(approximately $4.0\times 10^{20}$ and $3.0\times 10^{13}\mathrm{kg\cdot m^2s^{-2}}$ respectively).
As such the flattening of the density contours from which the baroclinic instability draws energy
are barely evident in Fig. \ref{fig::zonal_avg_1}.

\begin{figure}[!hbtp]
\begin{center}
\includegraphics[width=0.48\textwidth,height=0.36\textwidth]{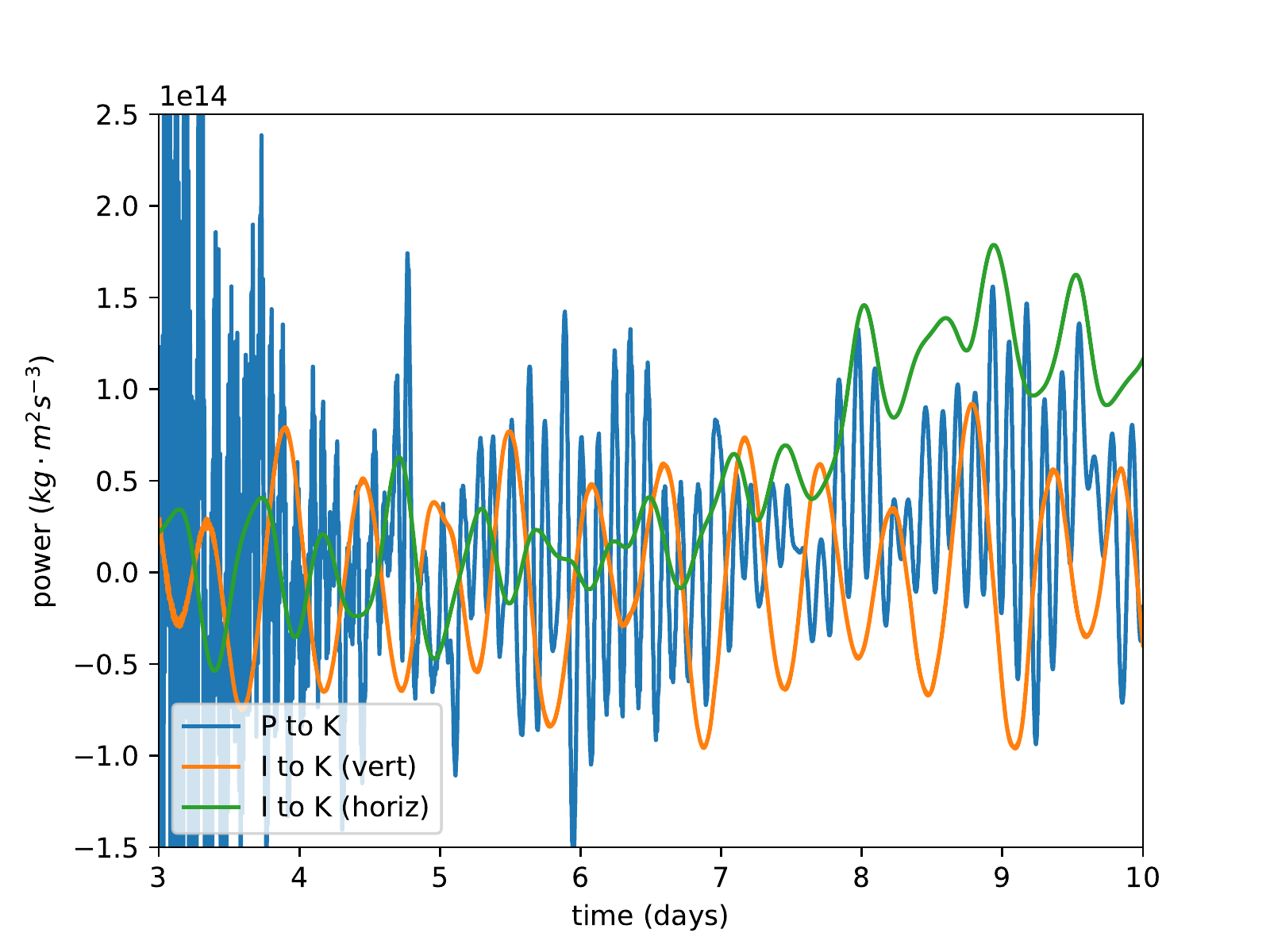}
\includegraphics[width=0.48\textwidth,height=0.36\textwidth]{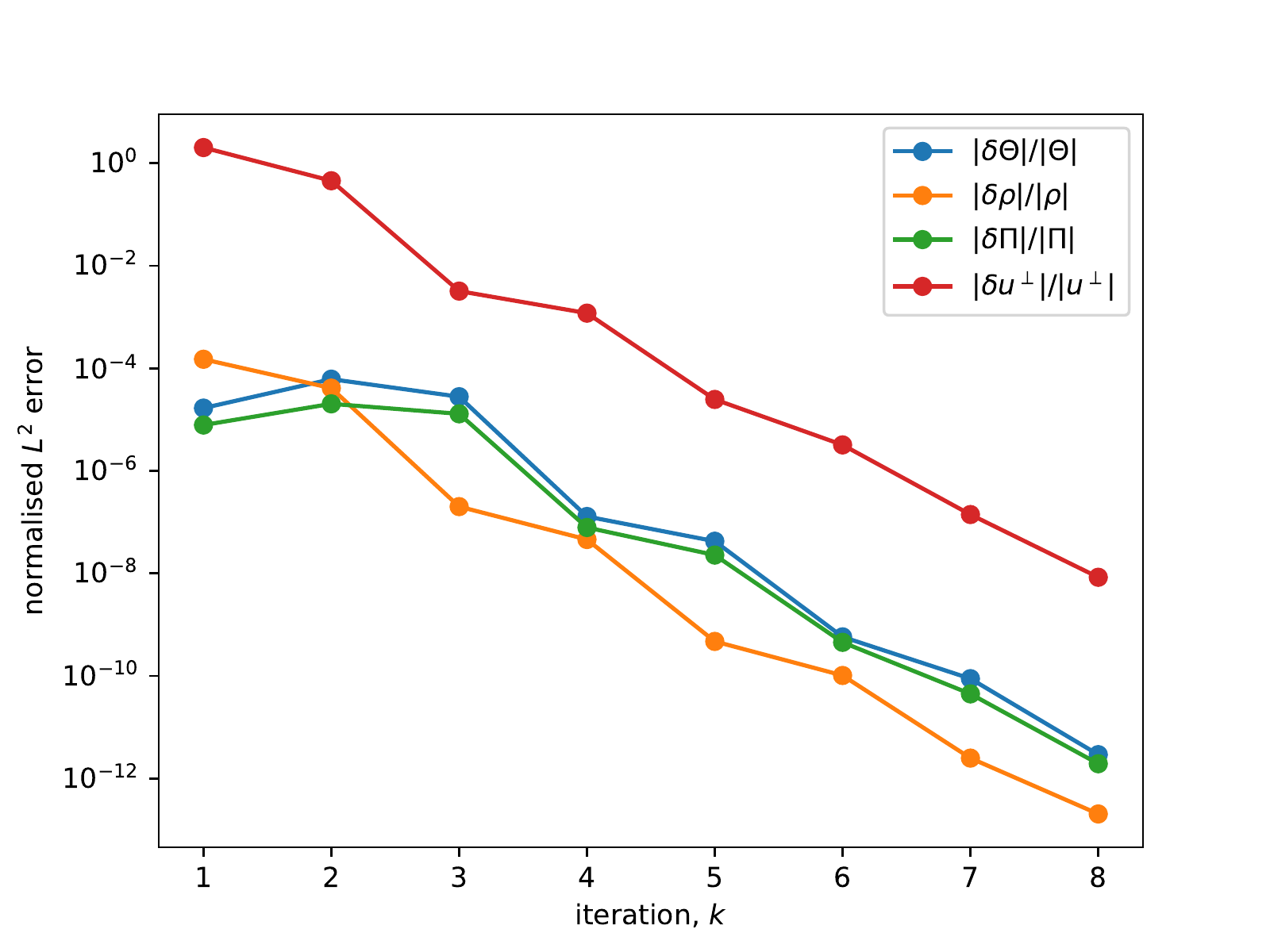}
\caption{
Baroclinic instability test case: 
Left: power associated with energetic exchanges, days 3 to 10. Right: typical 
convergence of the normalised $L^2$ magnitude of the solution updates with Newton iteration
(taken at day 9).}
\label{fig::energetics_2}
\end{center}
\end{figure}

Figure \ref{fig::energetics_2} shows the power exchanges between kinetic, potential and 
internal energy, as given in \eqref{eq::dKdt}-\eqref{eq::dIdt}. The evolution of the baroclinic
instability is evident in the increase of power associated with the internal to horizontal 
kinetic energy exchanges. This figure also shows the typical convergence of the iterative 
solver with Newton iteration (taken at day 9). As can be seen, the convergence is poor
for the first couple of Newton iterations, and the convergence of the velocity lags the other
variables. 
After the first couple of iterations, the convergence 
improves, and reduces the errors between one and two decades at each iteration. Note that
as the Jacobian is only approximate, we do not anticipate the convergence to be strictly
quadratic. Also note that this convergence plot shows the convergence for the column of 
elements with the largest error, which may not necessarily be the same column between 
Newton iterations.

\begin{figure}[!hbtp]
\begin{center}
\includegraphics[width=0.48\textwidth,height=0.36\textwidth]{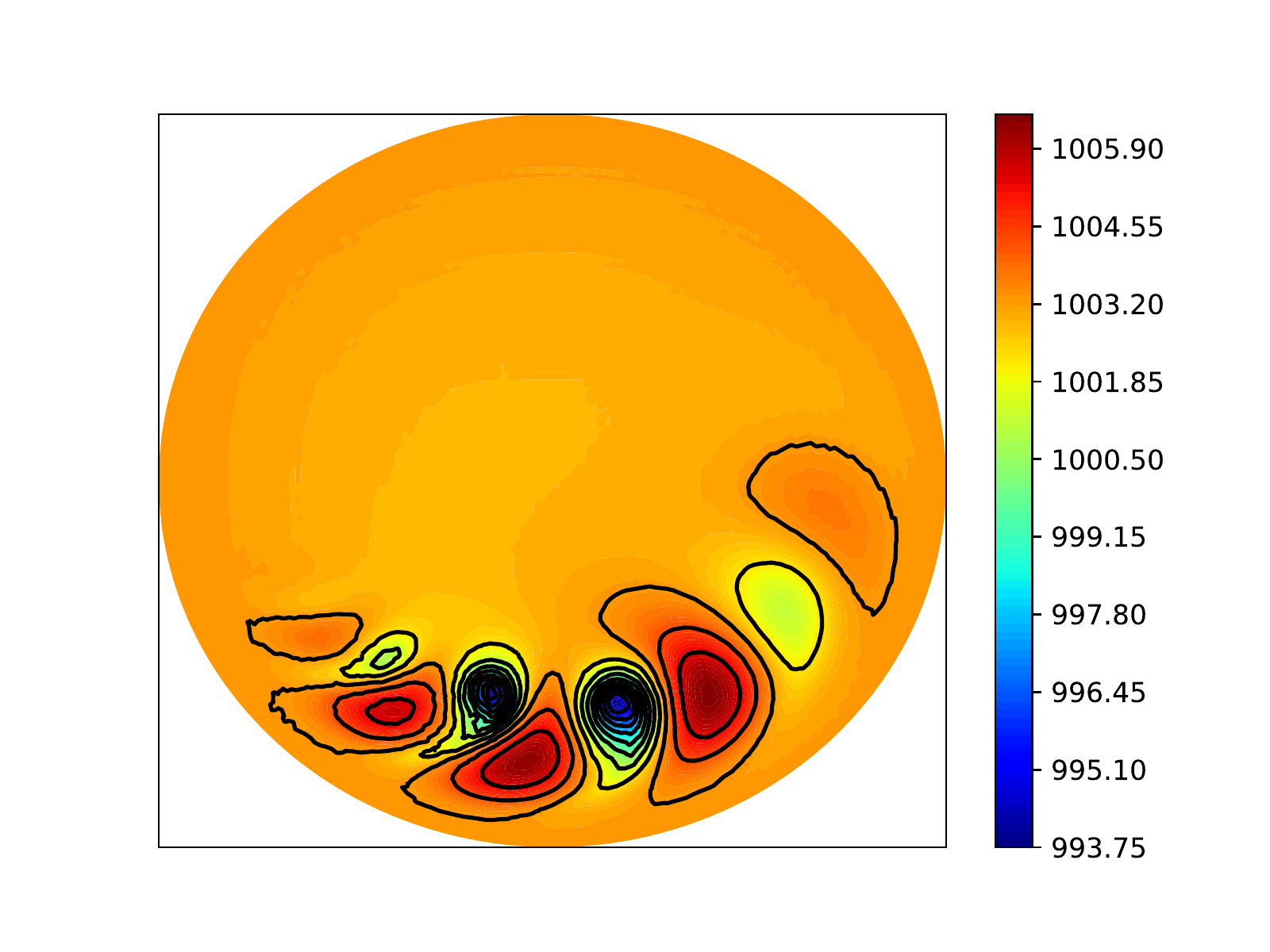}
\includegraphics[width=0.48\textwidth,height=0.36\textwidth]{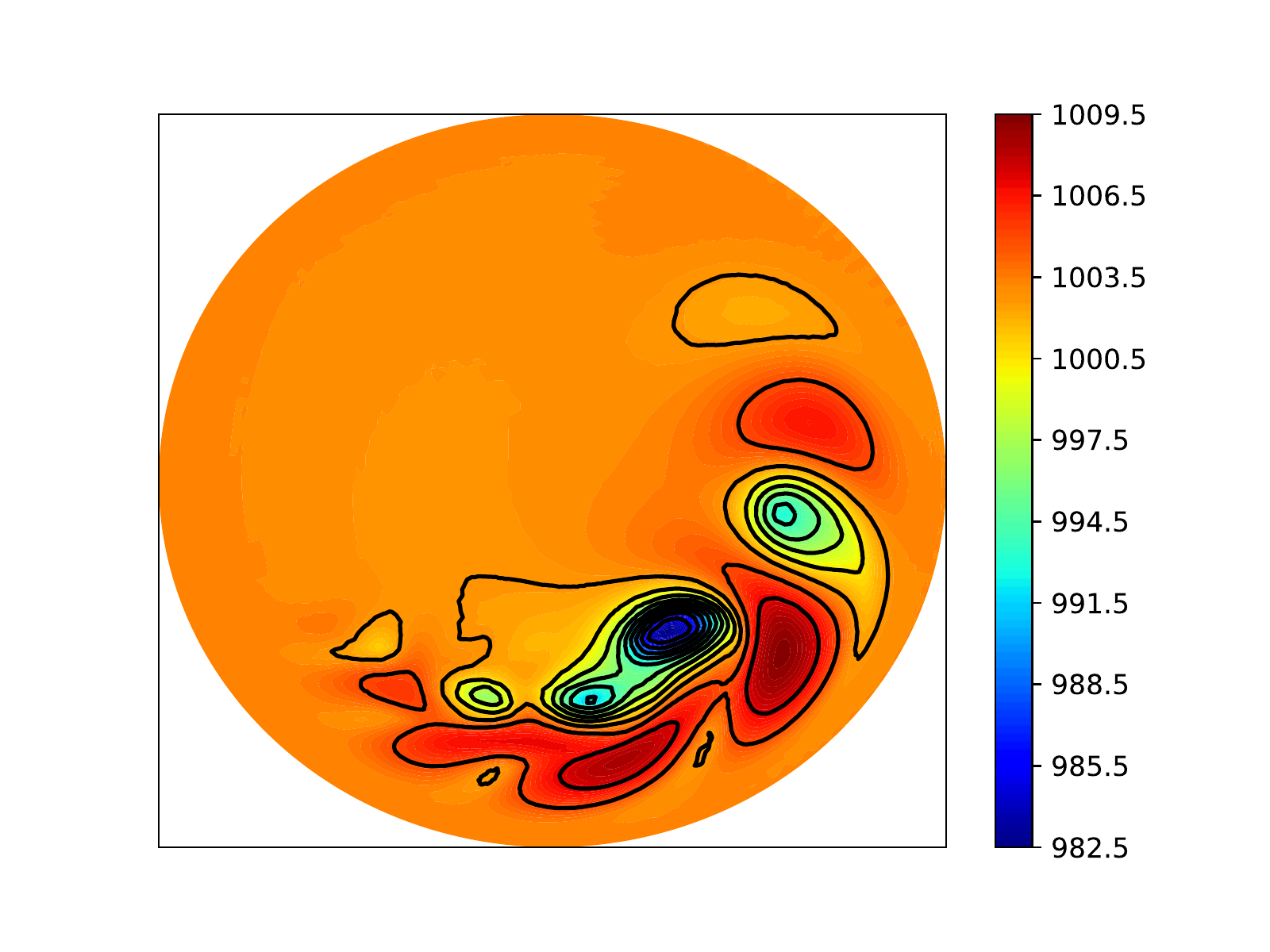}
\caption{
Baroclinic instability test case: 
Bottom level Exner pressure, $\Pi_h$ (in $\mathrm{hPa}$) day 8 (left) and 10 (right).}
\label{fig::pressure}
\end{center}
\end{figure}

\begin{figure}[!hbtp]
\begin{center}
\includegraphics[width=0.48\textwidth,height=0.36\textwidth]{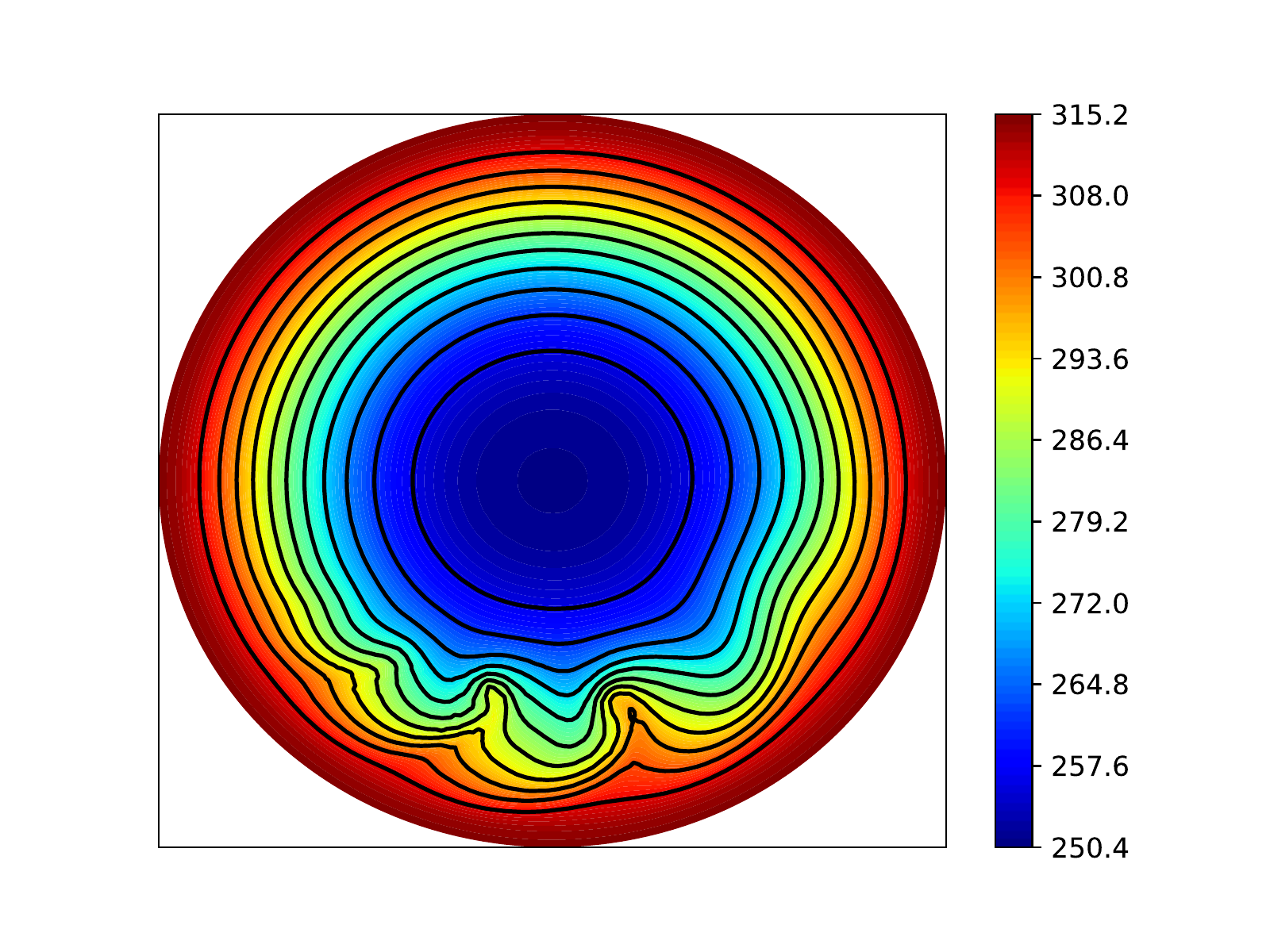}
\includegraphics[width=0.48\textwidth,height=0.36\textwidth]{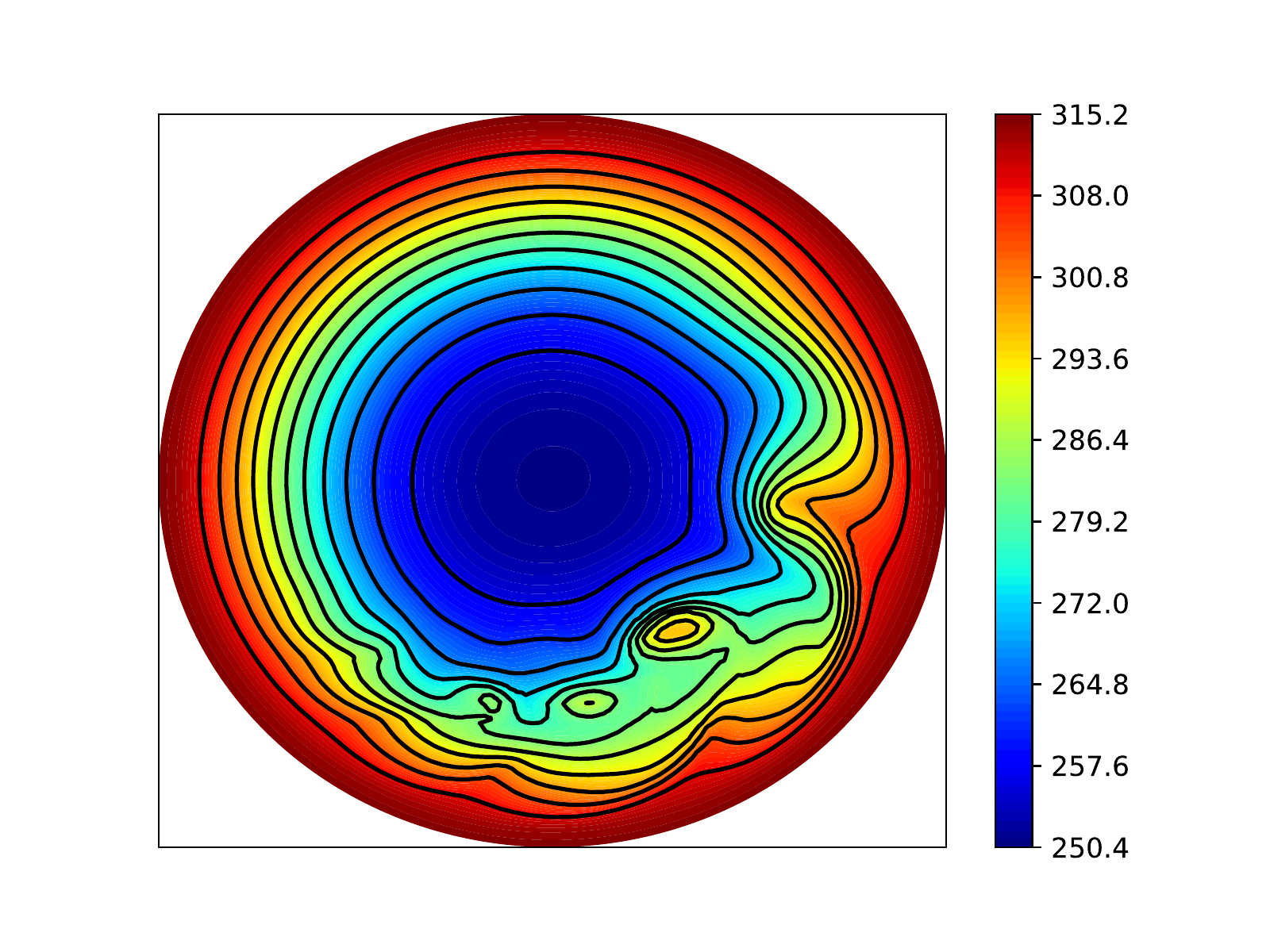}
\caption{
Baroclinic instability test case: 
Potential temperature, $\theta_h$ (in $^{\circ}\mathrm{K}$) at $z \approx 1.5\mathrm{km}$, 
day 8 (left) and 10 (right).}
\label{fig::temperature}
\end{center}
\end{figure}

\begin{figure}[!hbtp]
\begin{center}
\includegraphics[width=0.48\textwidth,height=0.36\textwidth]{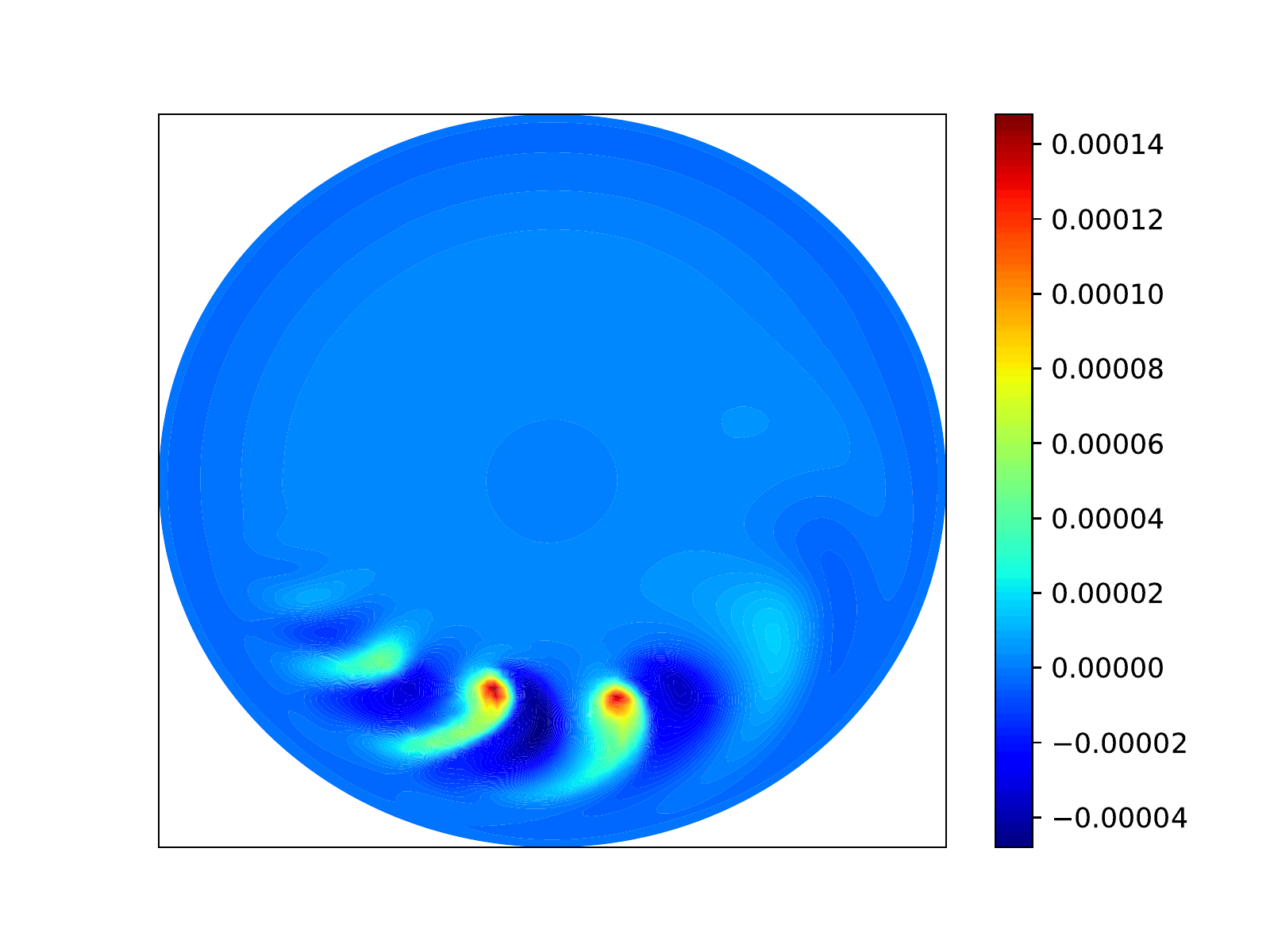}
\includegraphics[width=0.48\textwidth,height=0.36\textwidth]{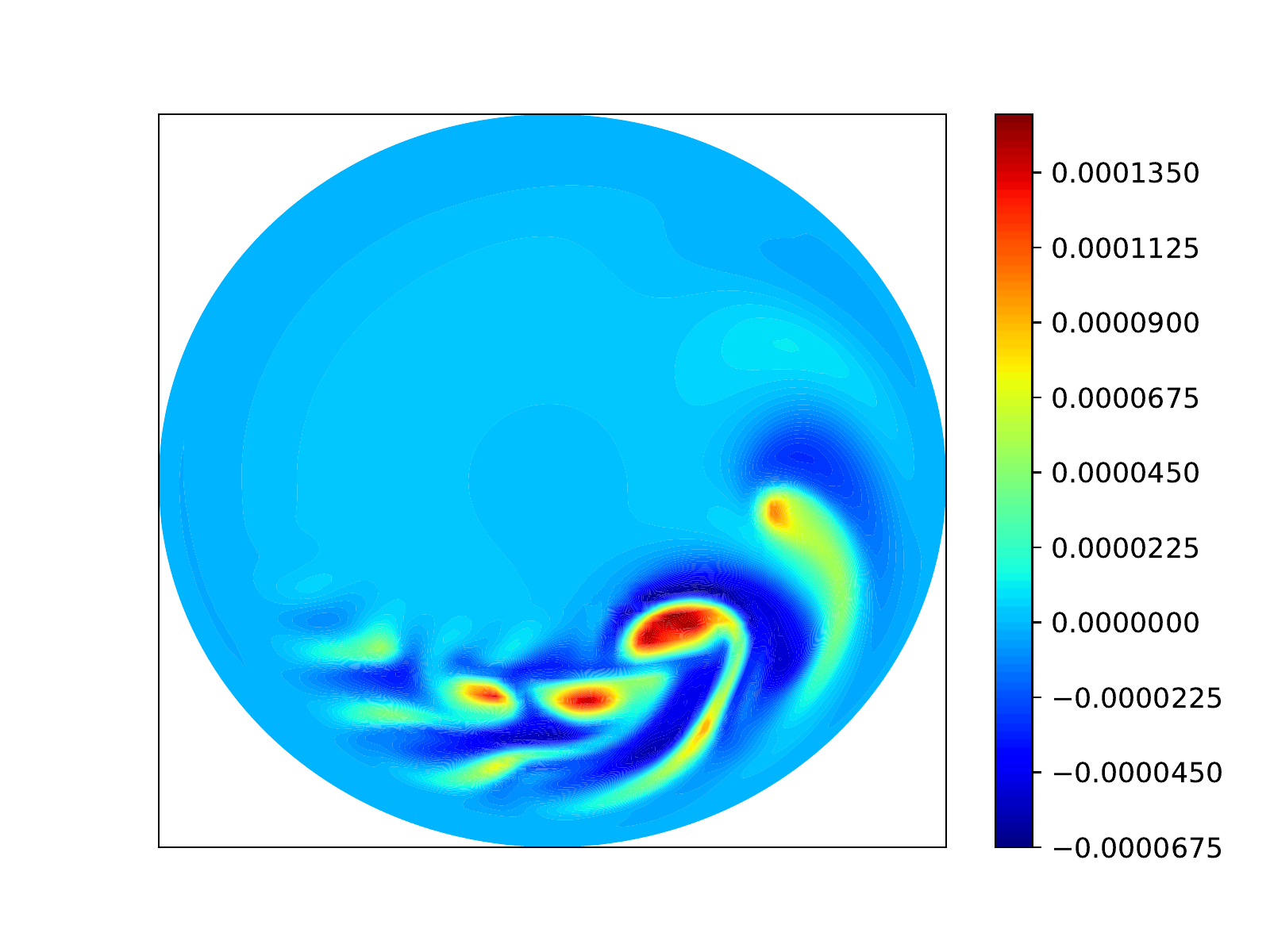}
\caption{
Baroclinic instability test case: 
Vertical component of the relative vorticity, $\omega_h$ (in $\mathrm{s}^{-1}$)
at $z \approx 1.5\mathrm{km}$, day 8 (left) and 10 (right).}
\label{fig::vorticity}
\end{center}
\end{figure}

\begin{figure}[!hbtp]
\begin{center}
\includegraphics[width=0.48\textwidth,height=0.36\textwidth]{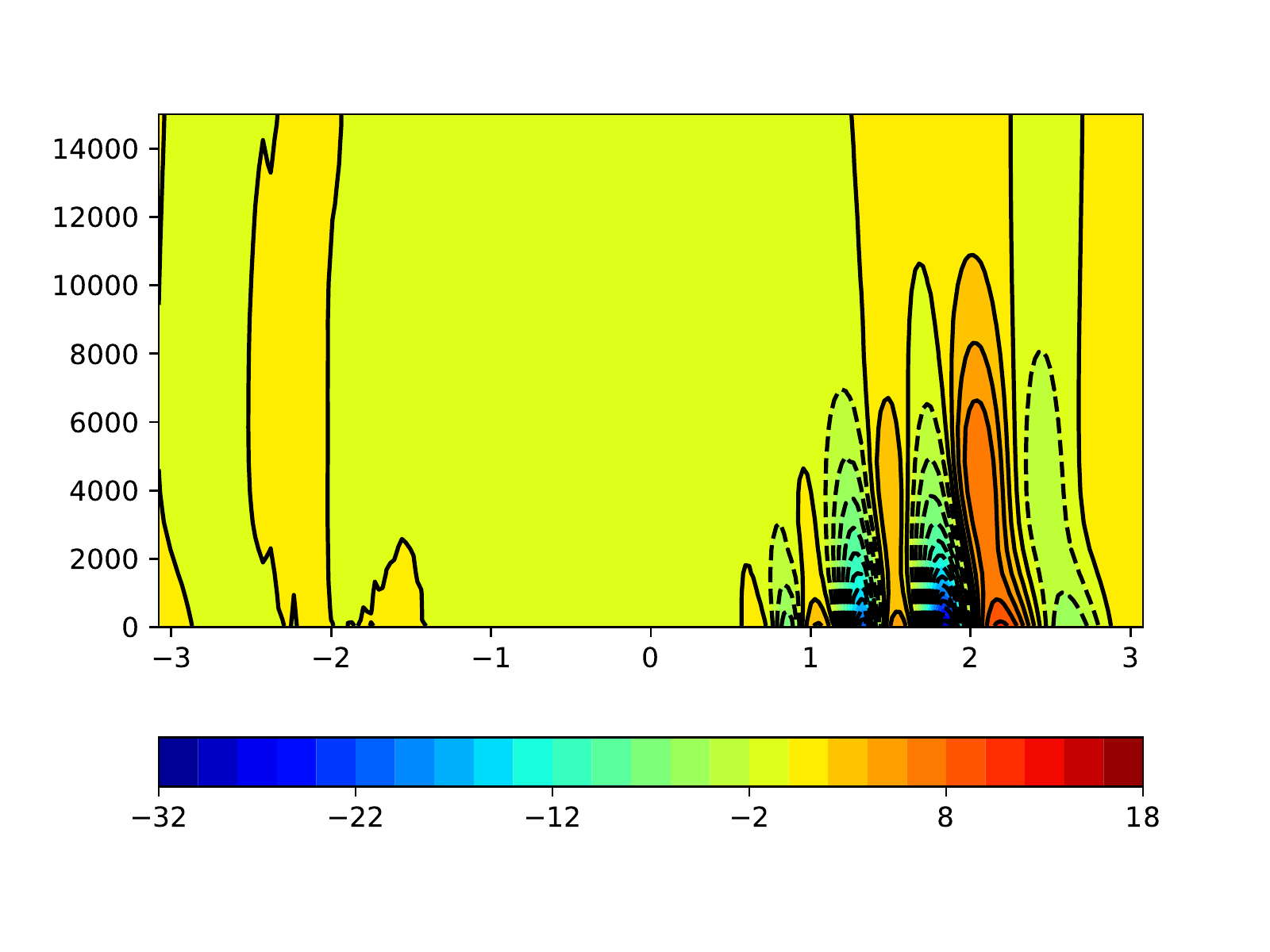}
\includegraphics[width=0.48\textwidth,height=0.36\textwidth]{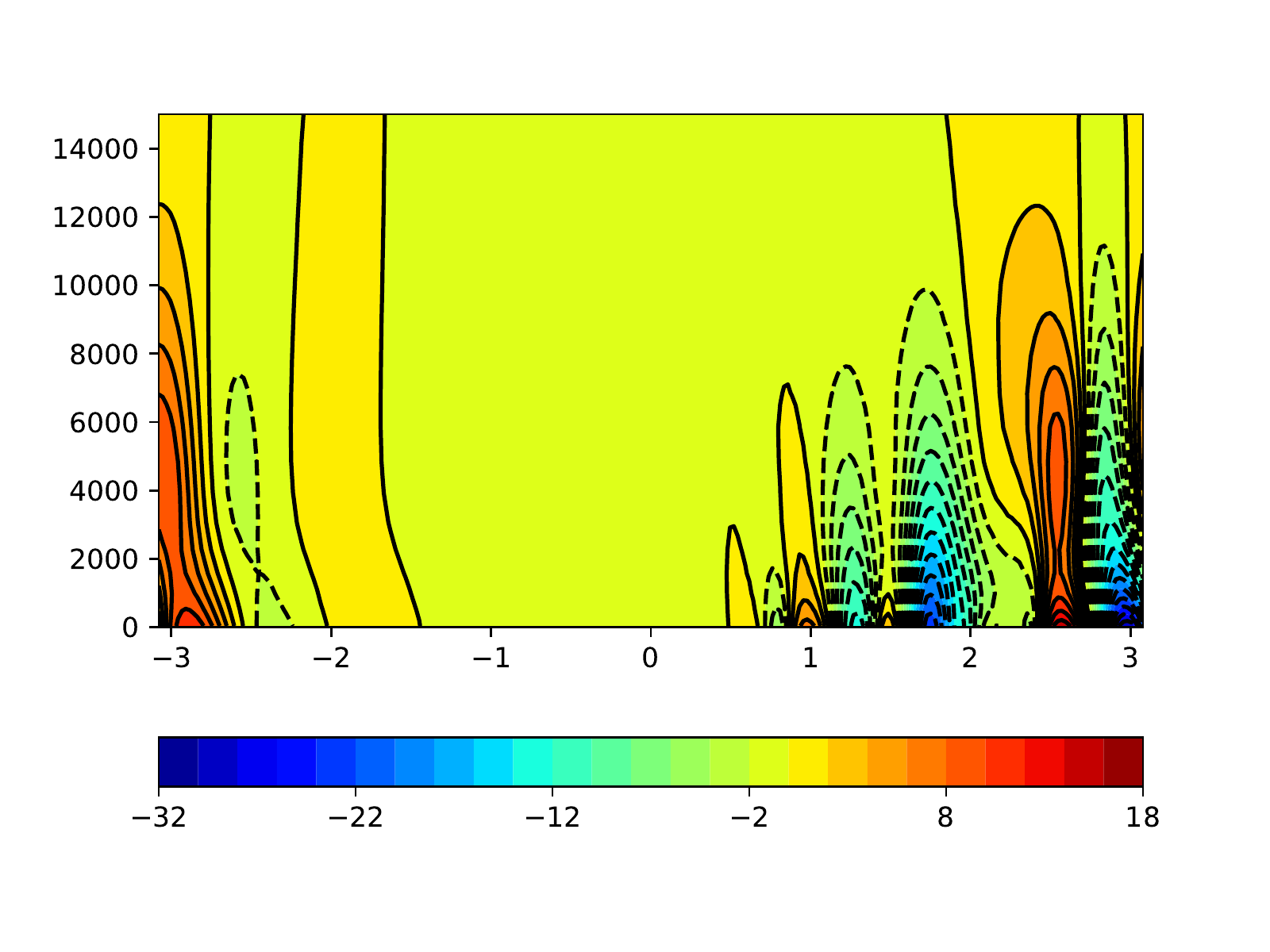}
\caption{Baroclinic instability test case: 
Vertical cross section of the pressure perturbation, $p_h-\bar p_h$ (in $\mathrm{hPa}$)
at $50^{\circ}\mathrm{N}$, day 8 (left) and day 10 (right).}
\label{fig::pressure_cross_section}
\end{center}
\end{figure}

The bottom level Exner pressure $\Pi$, potential temperature, $\theta$, and the vertical 
component of the relative vorticity, $\omega$ are presented for $z\approx 1.5\mathrm{km}$ 
at days $8$ and $10$ in Figs. \ref{fig::pressure}, \ref{fig::temperature} and \ref{fig::vorticity},
and the meridional cross section of the pressure perturbation at $50^{\circ}\mathrm{N}$ is given
in Fig. \ref{fig::pressure_cross_section}. In the cases of the pressure, this is reconstructed 
from the model variables as $p=p_0(\Pi/c_p)^{c_p/R}$.
The pressure perturbation in Fig. \ref{fig::pressure_cross_section} is then derived by removing
the average pressure at the corresponding vertical level at $50^{\circ}\mathrm{S}$. These results
compare well with the previously published test case results \cite{Ullrich14}, and the vertical 
pressure perturbation profile exhibits less noise than the previous version of the implicit solver
\cite{LP20}.

\begin{figure}[!hbtp]
\begin{center}
\includegraphics[width=0.48\textwidth,height=0.36\textwidth]{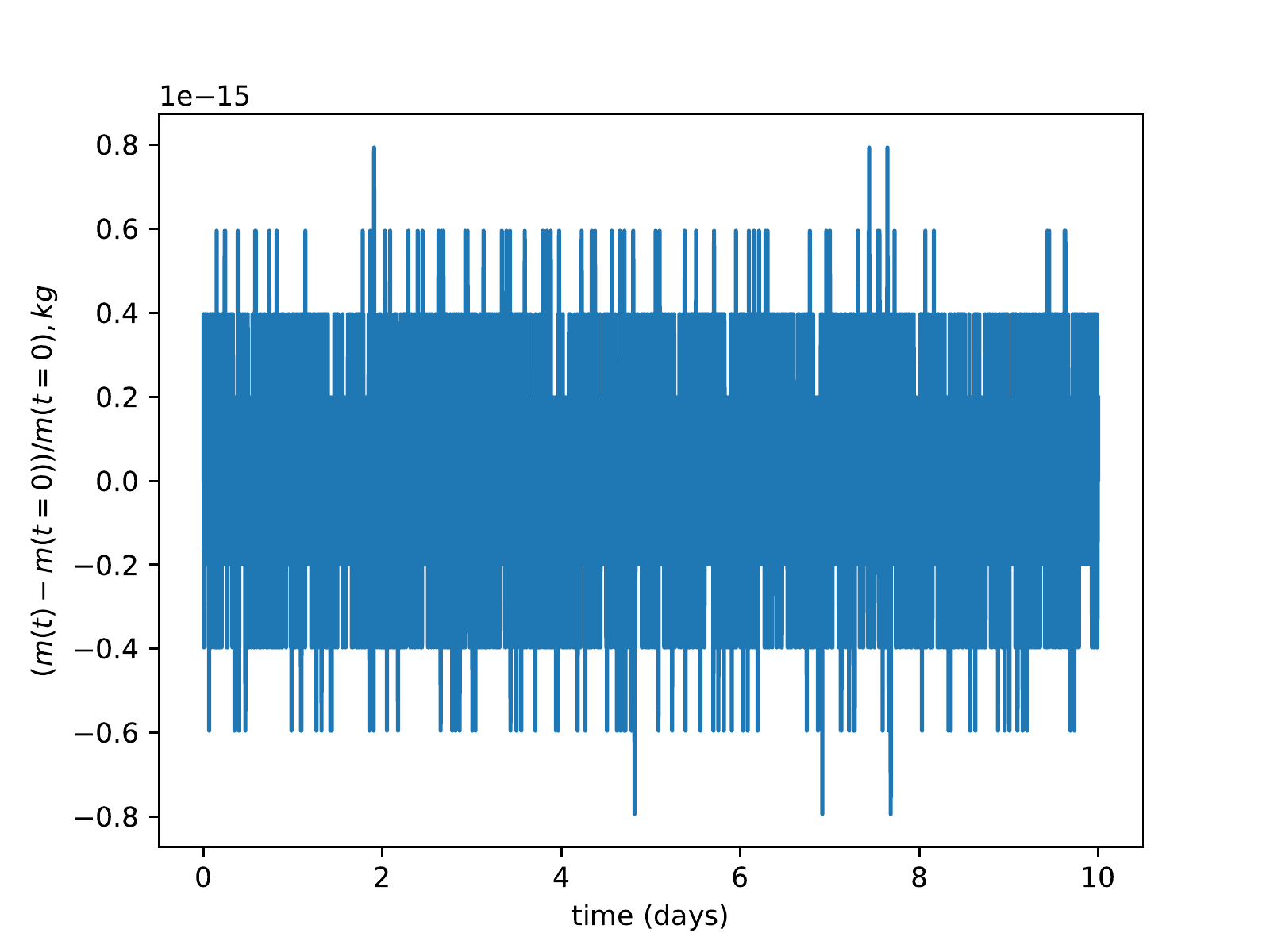}
\includegraphics[width=0.48\textwidth,height=0.36\textwidth]{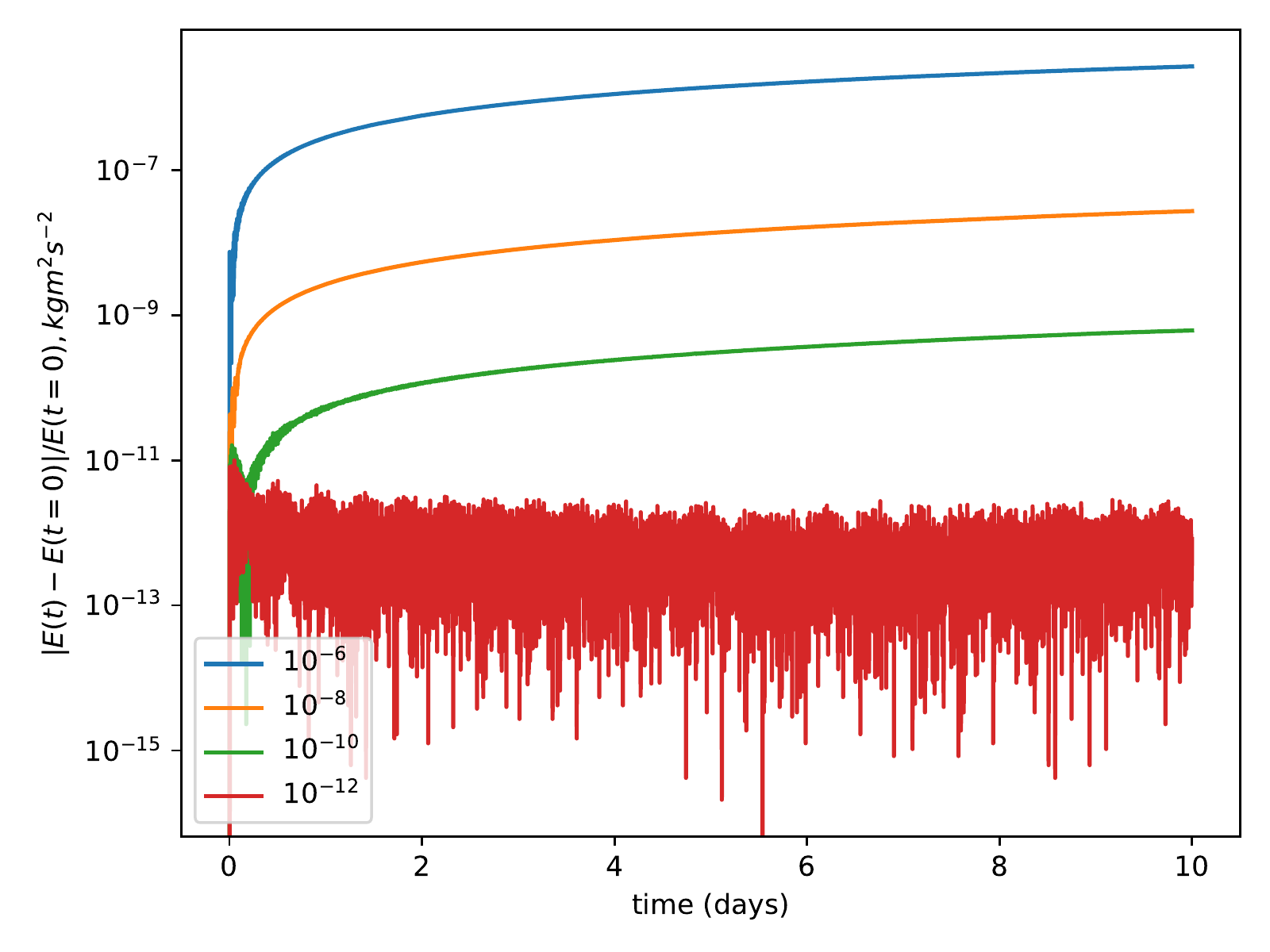}
\caption{Baroclinic instability test case (vertical dynamics only): 
Global conservation errors for mass (left) and total energy as a function of solver tolerance (right).}
\label{fig::conservation}
\end{center}
\end{figure}

In order to validate the conservation properties of the implicit solver the horizontal dynamics
were turned off and the model was re-run with only the implicit vertical solver active (such 
that the horizontal velocity was not updated between time steps), and all the dissipation terms
removed. Note that while the horizontal 
scheme conserves mass \cite{LP20}, it does not conserve energy due to the explicit time integration.
As shown in Fig. \ref{fig::conservation} while the mass is conserved to machine precision, 
the conservation of energy is satisfied only for a solver tolerance below approximately $10.0^{-12}$.
For tolerances greater than this a small drift in energy conservation is observed.

\begin{figure}[!hbtp]
\begin{center}
\includegraphics[width=0.48\textwidth,height=0.36\textwidth]{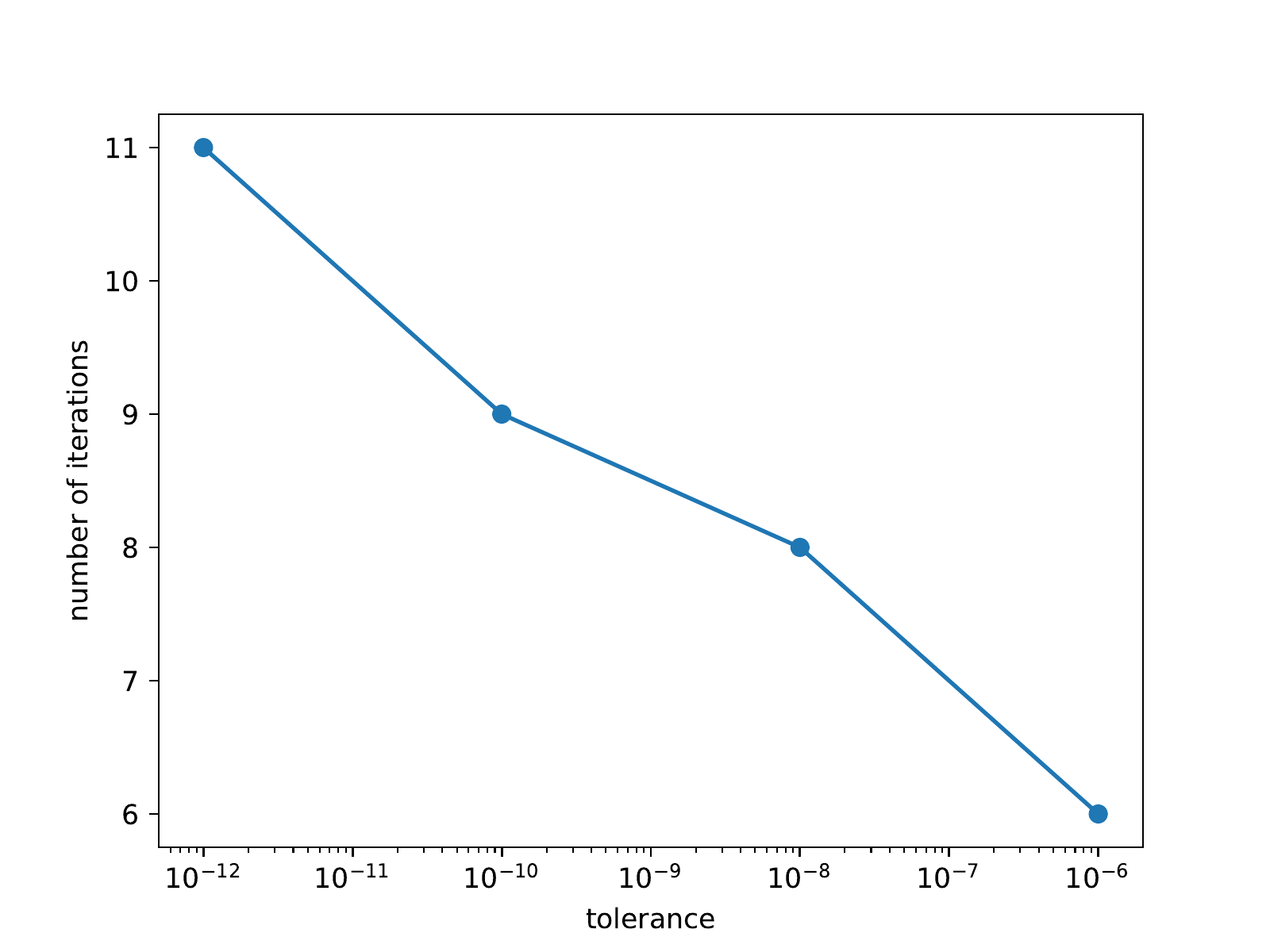}
\caption{Baroclinic instability test case (vertical dynamics only): 
Number of nonlinear iterations to reach convergence as function of solver tolerance.}
\label{fig::conservation_1}
\end{center}
\end{figure}

Figure \ref{fig::conservation_1} shows the number of nonlinear iterations required to achieve 
convergence as a function of solver tolerance (again with the horizontal dynamics disabled). Note that the 
vertical velocity error has been removed from the set of convergence criteria, since the vertical kinetic 
energy is extremely small compared to the amounts of potential and internal energy, and the vertical velocity
convergence error lags that of the other solution variables, as shown in Fig. \ref{fig::energetics_2}. As can be 
observed the number of iterations increases somewhat linearly with solver tolerance at a rate of approximately
one iteration per decade. Note that the Jacobian as detailed in \eqref{eq::jacobian} and 
\eqref{eq::jacobian_blocks} is only approximate, so quadratic convergence is not anticipated.

The energy conservation errors in Fig. \ref{fig::conservation} are compared to those for the standard 
Crank-Nicolson semi-implicit integrator as used in each vertical sub-step of the original TRAP(2,3,2) scheme \cite{Weller13,Lock14}.
The Crank-Nicolson scheme is reverse-engineered from the existing energetically balanced scheme by 
removing all of the cross terms between time level $n$ and Newton iteration $k$ as given in \eqref{eq::residuals}, 
\eqref{eq::variational_derivs_discrete}, and weighting the system at time levels $n$ and $k$ as $0.5$. All
other aspects of the code are un-modified, including the energetically balanced spatial discretisation \cite{LP20},
and in both cases a tolerance of $10^{-12}$ is used (with the vertical velocity again omitted from the convergence 
criteria). 

\begin{figure}[!hbtp]
\begin{center}
\includegraphics[width=0.48\textwidth,height=0.36\textwidth]{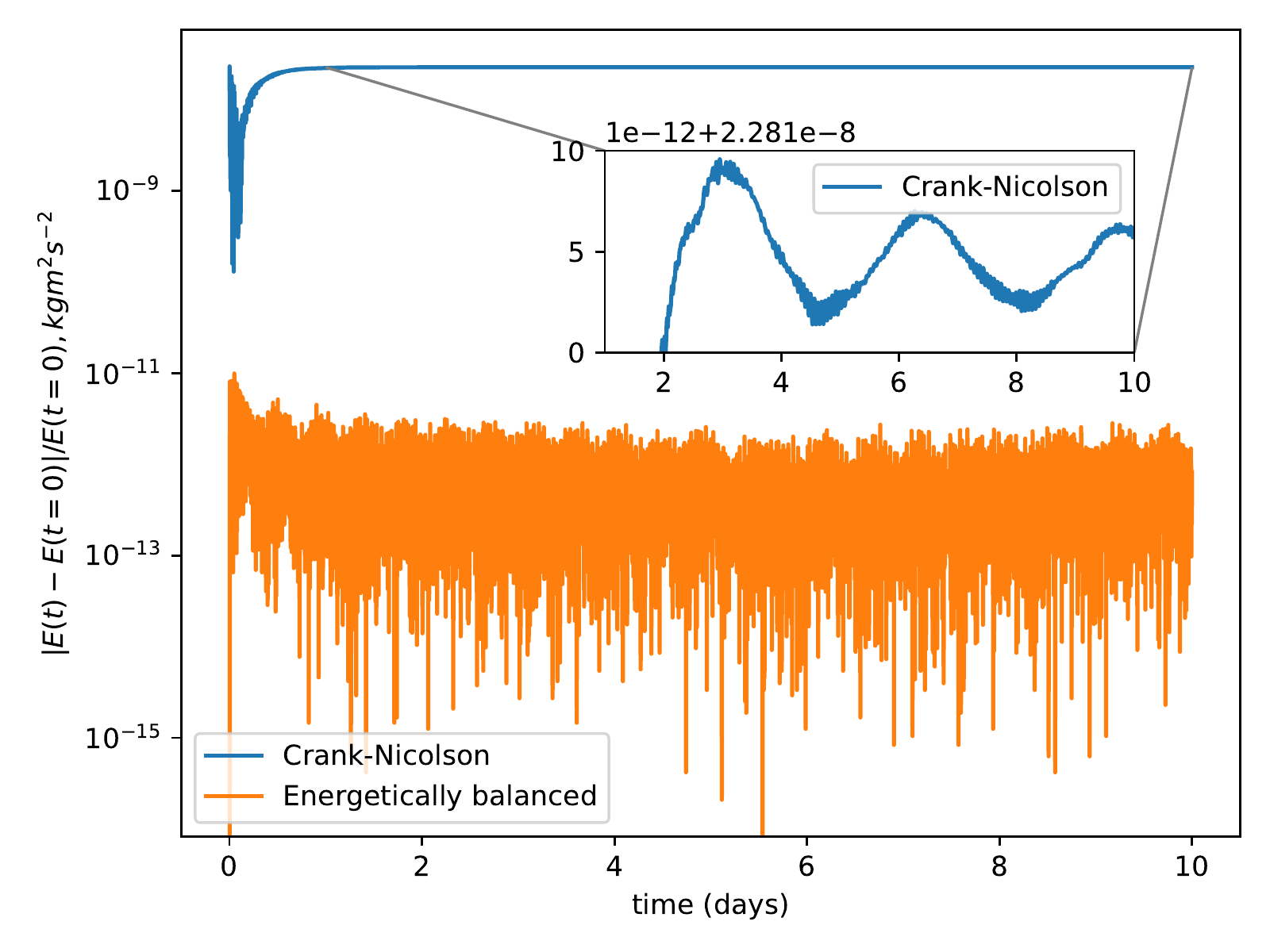}
\includegraphics[width=0.48\textwidth,height=0.36\textwidth]{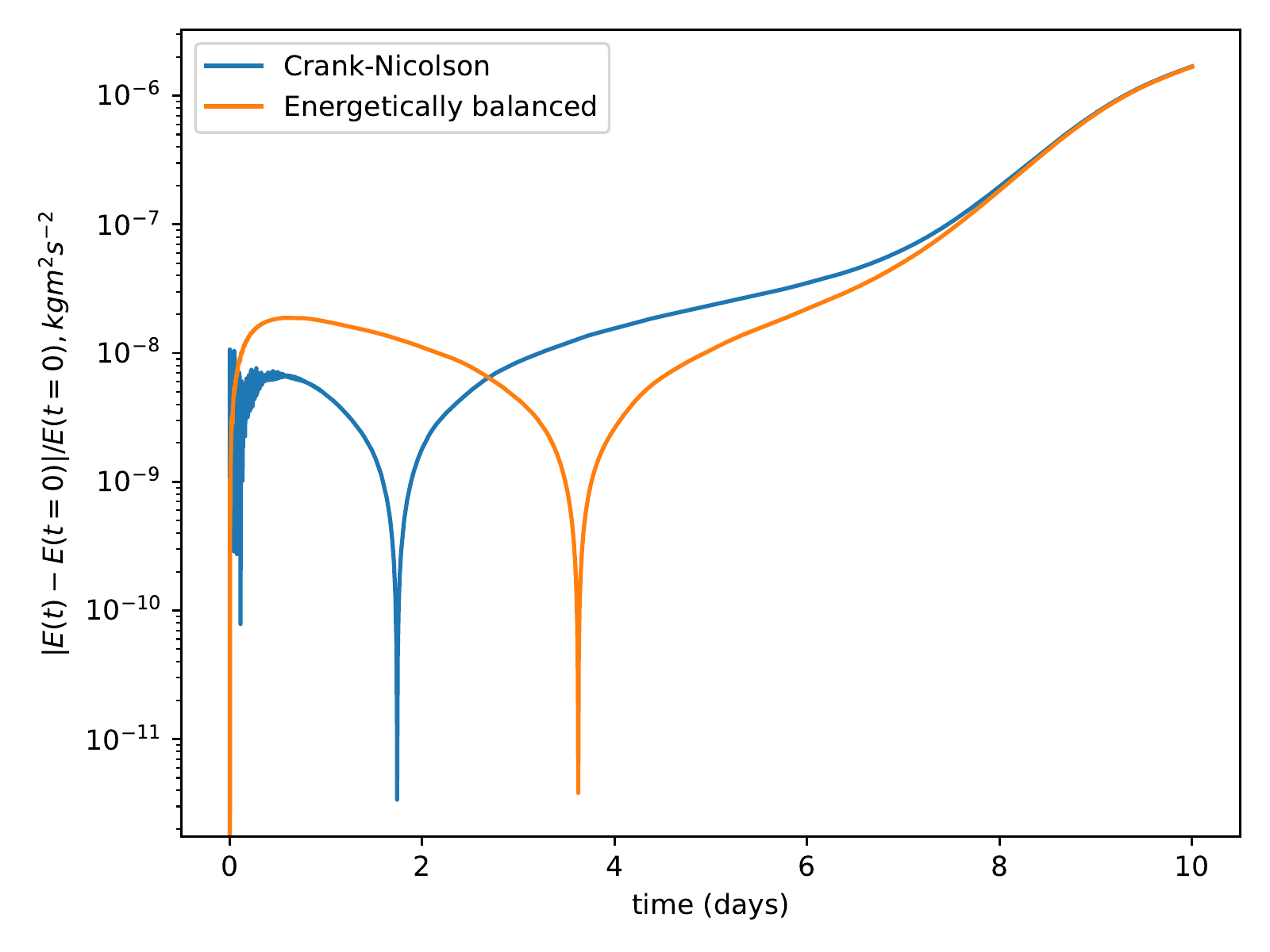}
\caption{Baroclinic instability test case: 
Comparison of global energy conservation errors for the energetically balanced and standard Crank-Nicolson
integrators, left: vertical dynamical only, right: full model (normalised absolute value
of error).}
\label{fig::conservation_2}
\end{center}
\end{figure}

As observed in Fig \ref{fig::conservation_2}, if only the vertical dynamics are considered, then the energy 
conservation errors of the Crank-Nicolson scheme, while still small, are $\mathcal{O}(10^4)$ times greater 
than those of the energetically balanced integrator.
While the Crank-Nicolson scheme is known to be unconditionally stable, such that the energy conservation
errors remain bounded, the exchanges in energy remain unbalanced, and there is a noticeable oscillation in total
energy with a period of approximately 3.5 days. This is despite the fact that both schemes involve balanced 
energy exchanges within the spatial discretisation, and differ only in the evaluation of additional cross 
terms in the temporal discretisation. For the full three
dimensional model, the Crank-Nicolson integrator exhibits greater oscillation in the energy conservation
errors during the initial period of hydrostatic adjustment, however after this initial period the energy 
conservation errors are actually greater for the energetically balanced scheme. This is perhaps due to the
inclusion of the Rayleigh damping term in the top three levels. The energetically balanced scheme is perhaps
more efficient at transmitting waves to the top of the domain, where these motions are damped by Rayligh friction. 
Nevertheless the growth in energy conservation error is faster for the Crank-Nicolson scheme over this period.
Over the course of the simulation both integrators tend towards the same conservation error (which is marginally
lower for the energetically balanced scheme), suggesting that 
for long times horizontal dissipative terms that are the same in both models dominate the energy conservation error.

\subsection{Non-hydrostatic gravity wave}

In order to validate the vertical integrator at high resolutions at which non-hydrostatic dynamics
become significant the model has also been tested for the propagation of
a non-hydrostatic gravity wave driven by a potential temperature perturbation on a planet with
a reduced radius $125$ times smaller than that of the earth. This test was originally proposed as part
of the 2012 DCMIP workshop \cite{DCMIP31}, and specific details of the initial configuration can be found
within the DCMIP test case document on the web site.

As for the baroclinic instability test, the simulation was run with a resolution of
$24\times 24$ elements of degree $p=3$ in each cubed sphere panel. We use 16 evenly spaced
vertical levels over a total height of $10,000.0$m, and a time step of $\Delta t = 0.5$s for
a total simulation time of $3600$s. In order to account for the highly oscillatory nature of 
the dynamics the horizontal biharmonic viscosity was rescaled by a factor of $2.0$ for both 
the momentum and temperature equations for a value of $0.144\Delta x^{3.2}$. 

One particular challenge of this test case is the inexact integration of the horizontal metric terms 
on a planet of reduced radius. These terms include transcendental functions that cannot be integrated
exactly \cite{Guba14,LP18}. While the relative curvature of earth the with respect to the element size is 
small for the case of an earth of regular size, on an earth of reduced size under-integration of the
metric terms may result in substantial aliasing errors that are not present in models based on
polyhedral meshes.

\begin{figure}[!hbtp]
\begin{center}
\includegraphics[width=0.48\textwidth,height=0.36\textwidth]{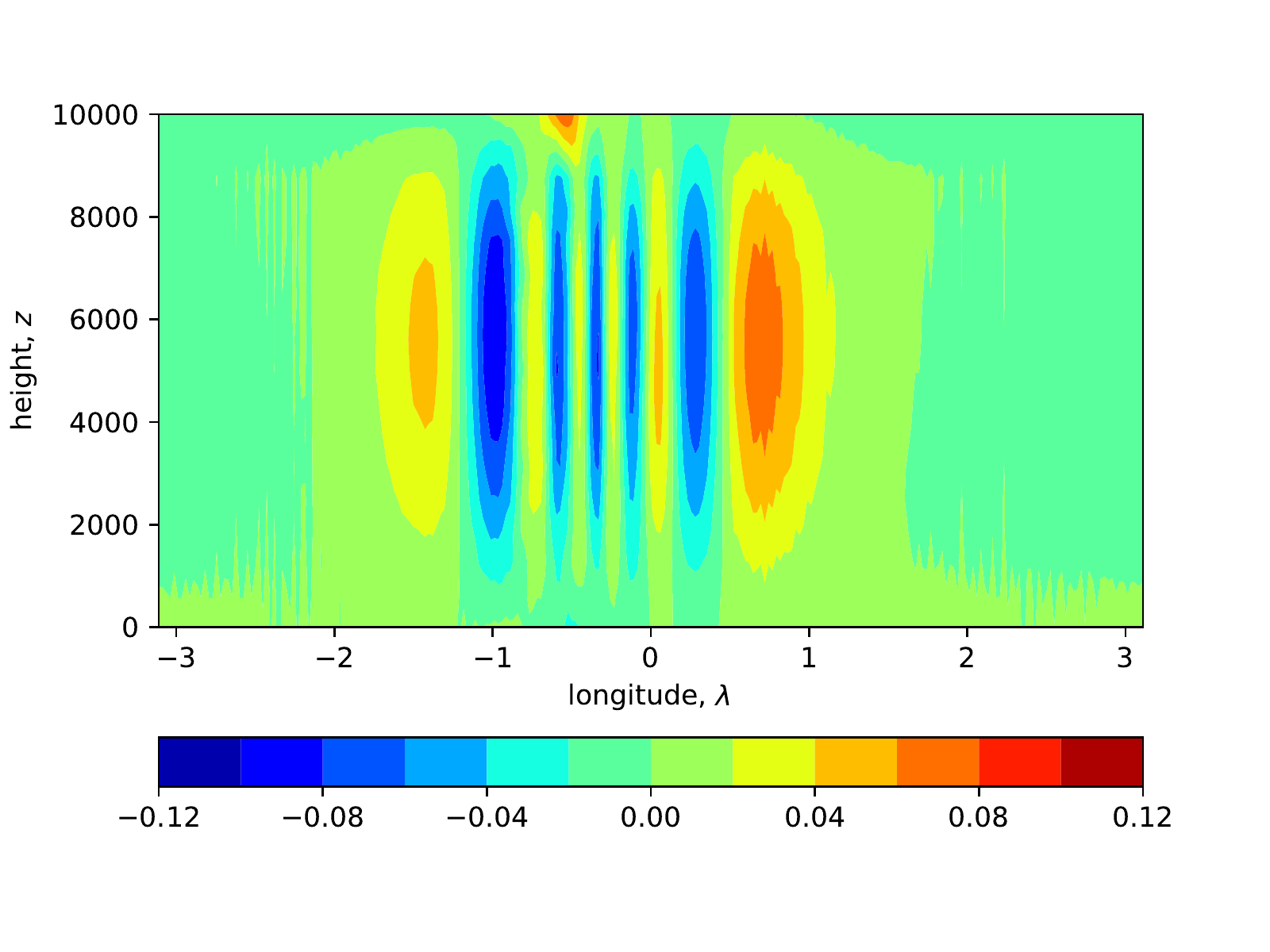}
\includegraphics[width=0.48\textwidth,height=0.36\textwidth]{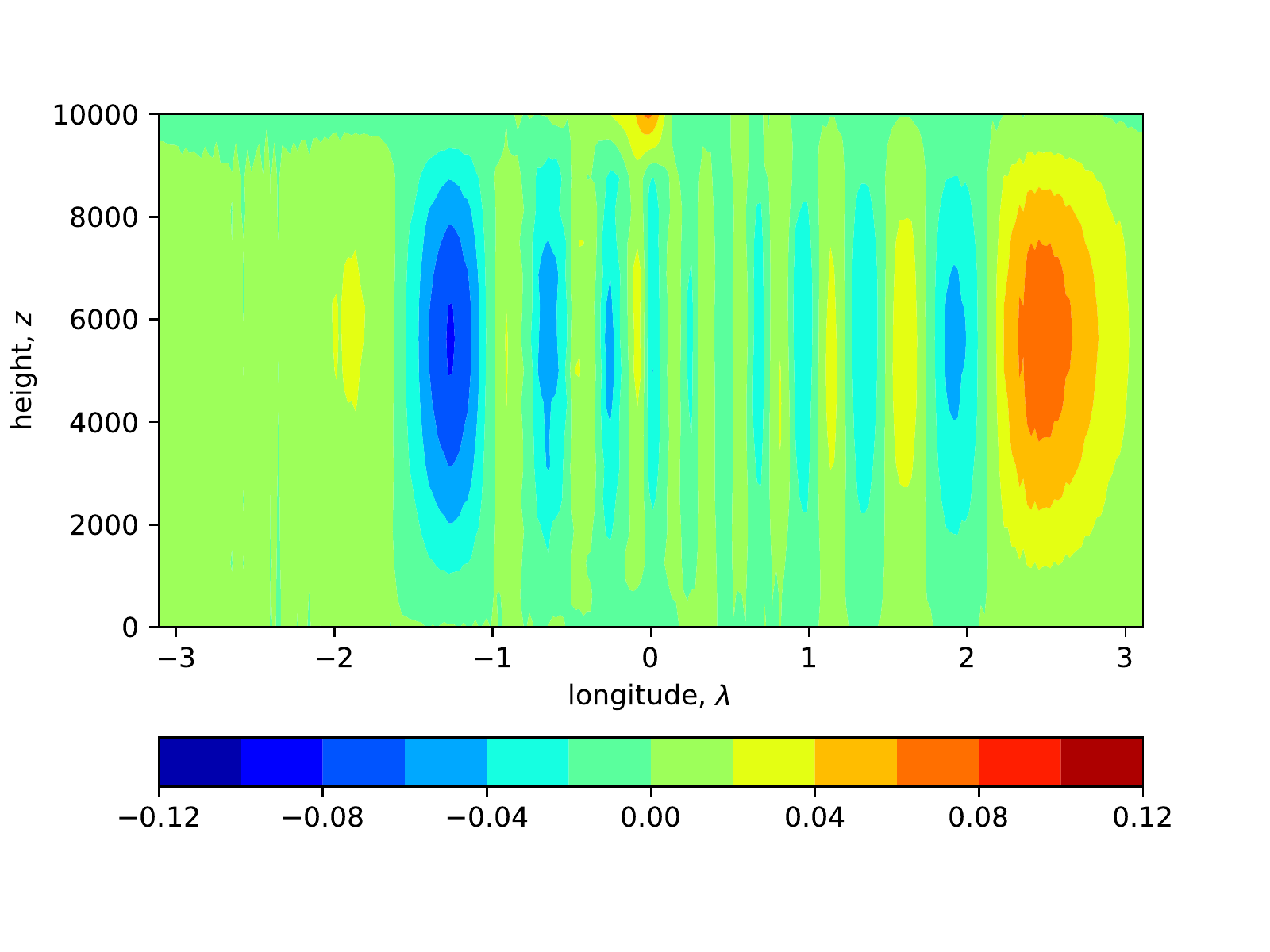}
\caption{Non-hydrostatic gravity wave test case: 
Longitude-height equatorial ($\phi=0^{\circ}$) cross section of the potential temperature
perturbation, $\theta'(\lambda,0,z) = \theta(\lambda,0,z) - \bar{\theta}(z)$
at times $t=30$ minutes (left) and $t=60$ minutes (right).}
\label{fig::grav_wave_theta}
\end{center}
\end{figure}

Figure \ref{fig::grav_wave_theta} shows the longitude-height equatorial ($\phi=0^{\circ}$) cross section of the potential
temperature perturbation, $\theta'(\lambda,0,z) = \theta(\lambda,0,z) - \bar{\theta}(z)$, where $\bar{\theta}(z)$ is the
mean potential temperature at a given height, after 30 minutes and 1 hour. While the structure
and evolution of the perturbation qualitatively match the results presented for other non-hydrostatic
models, the results here are slightly
more oscillatory than those of other models, perhaps reflecting the presence of aliasing errors due to the 
inexact integration of the metric terms on a planet of increased curvature. 
Moreover these results also show the ejection of a small
hot bubble which rises to the top of the domain, which is also perhaps an artefact of the 
inexact spatial integration.

\begin{figure}[!hbtp]
\begin{center}
\includegraphics[width=0.48\textwidth,height=0.36\textwidth]{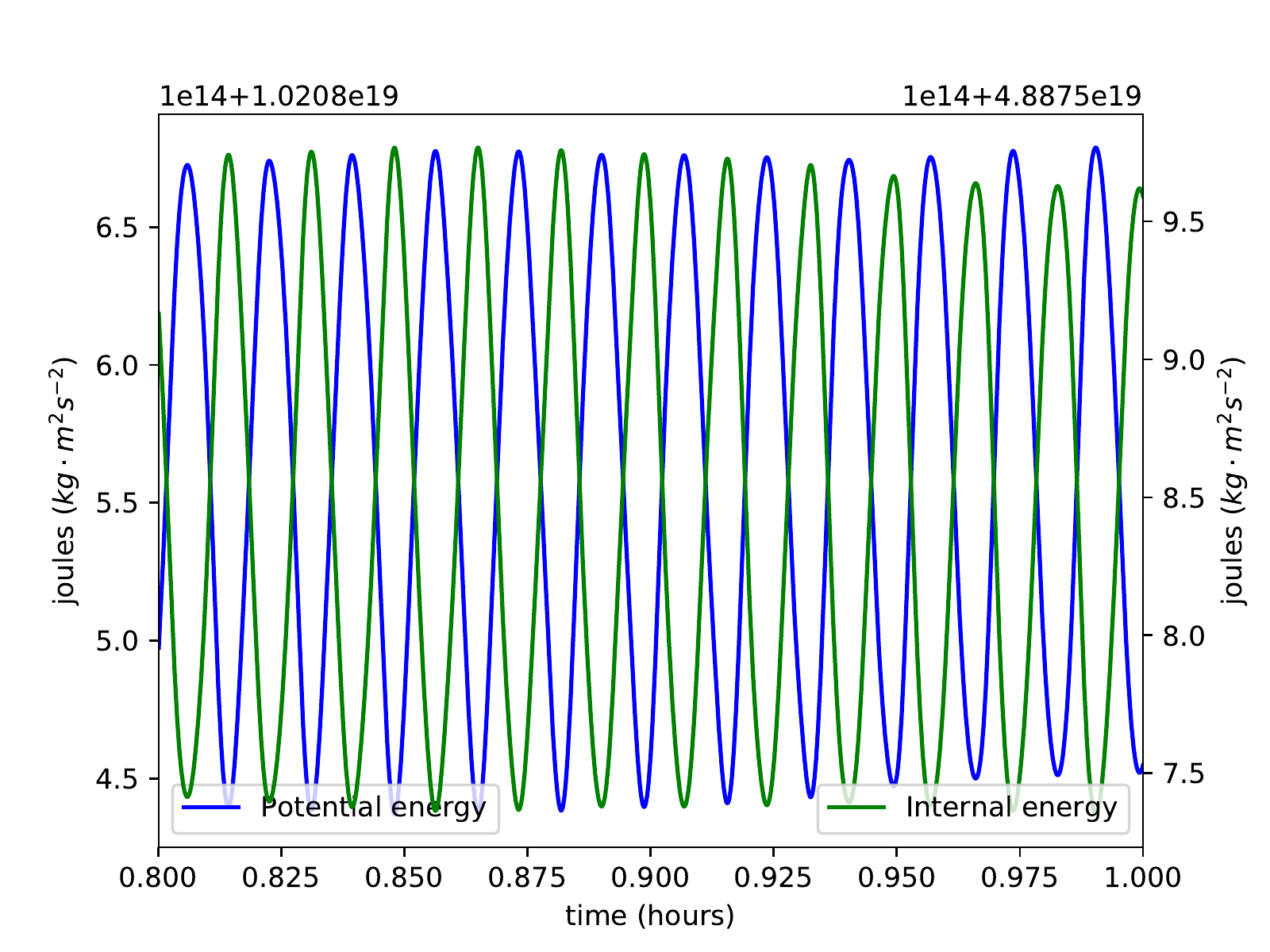}
\includegraphics[width=0.48\textwidth,height=0.36\textwidth]{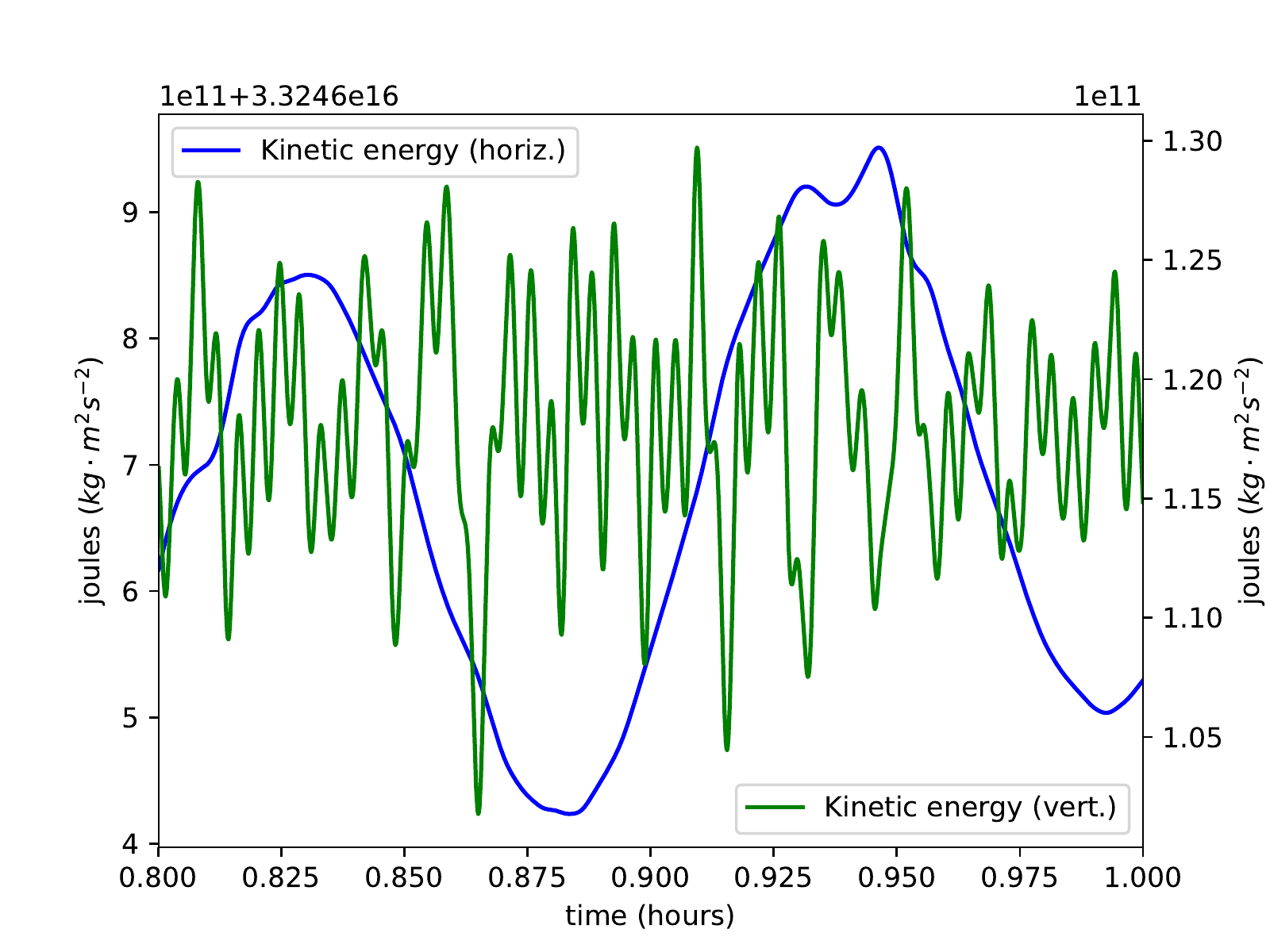}
\caption{Non-hydrostatic gravity wave test case: 
Internal gravity wave; left: time evolution of the potential and internal energies, right: 
time evolution of the kinetic energy (horizontal and vertical).}
\label{fig::energetics_1_gw}
\end{center}
\end{figure}

\begin{figure}[!hbtp]
\begin{center}
\includegraphics[width=0.48\textwidth,height=0.36\textwidth]{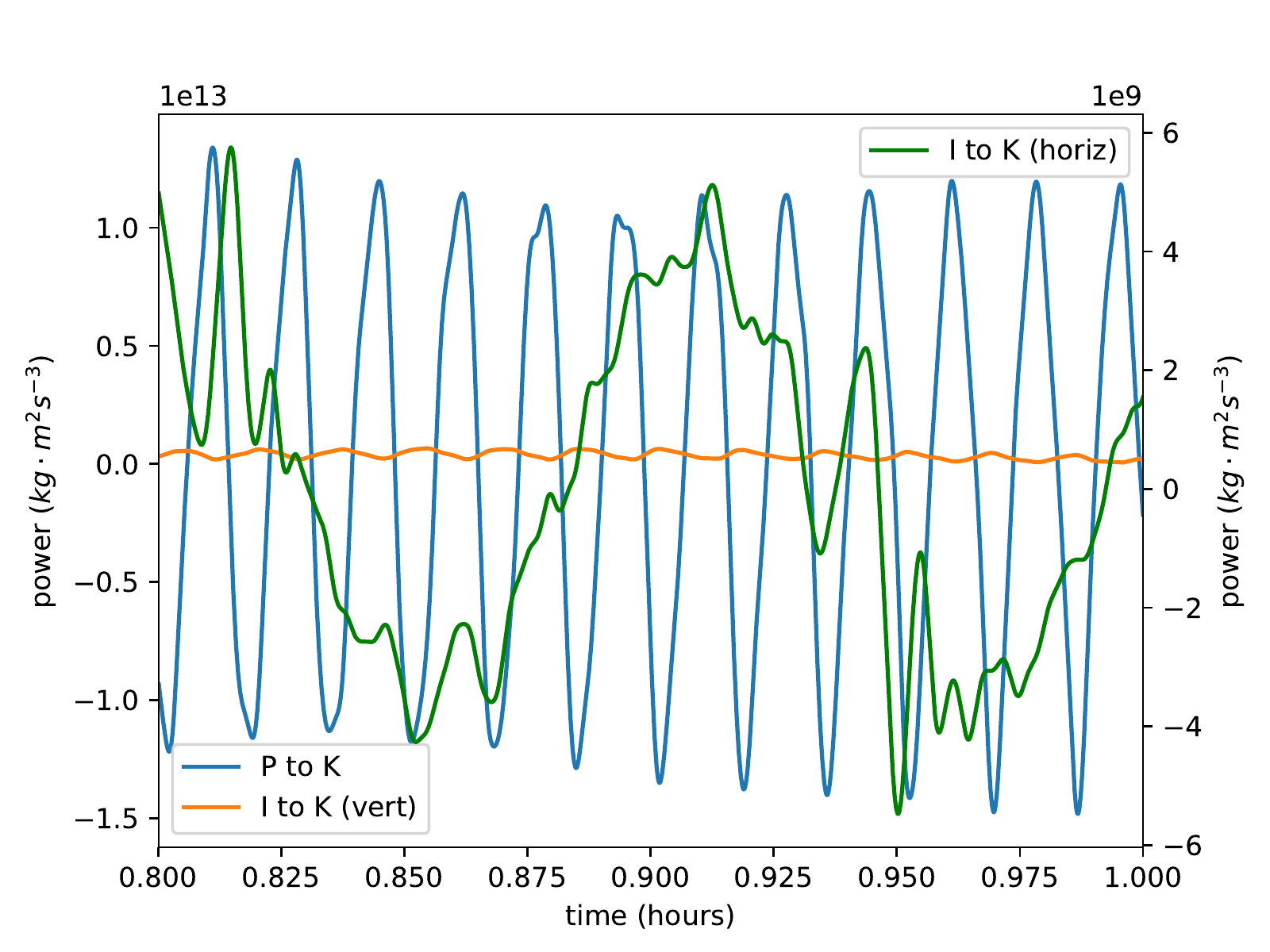}
\includegraphics[width=0.48\textwidth,height=0.36\textwidth]{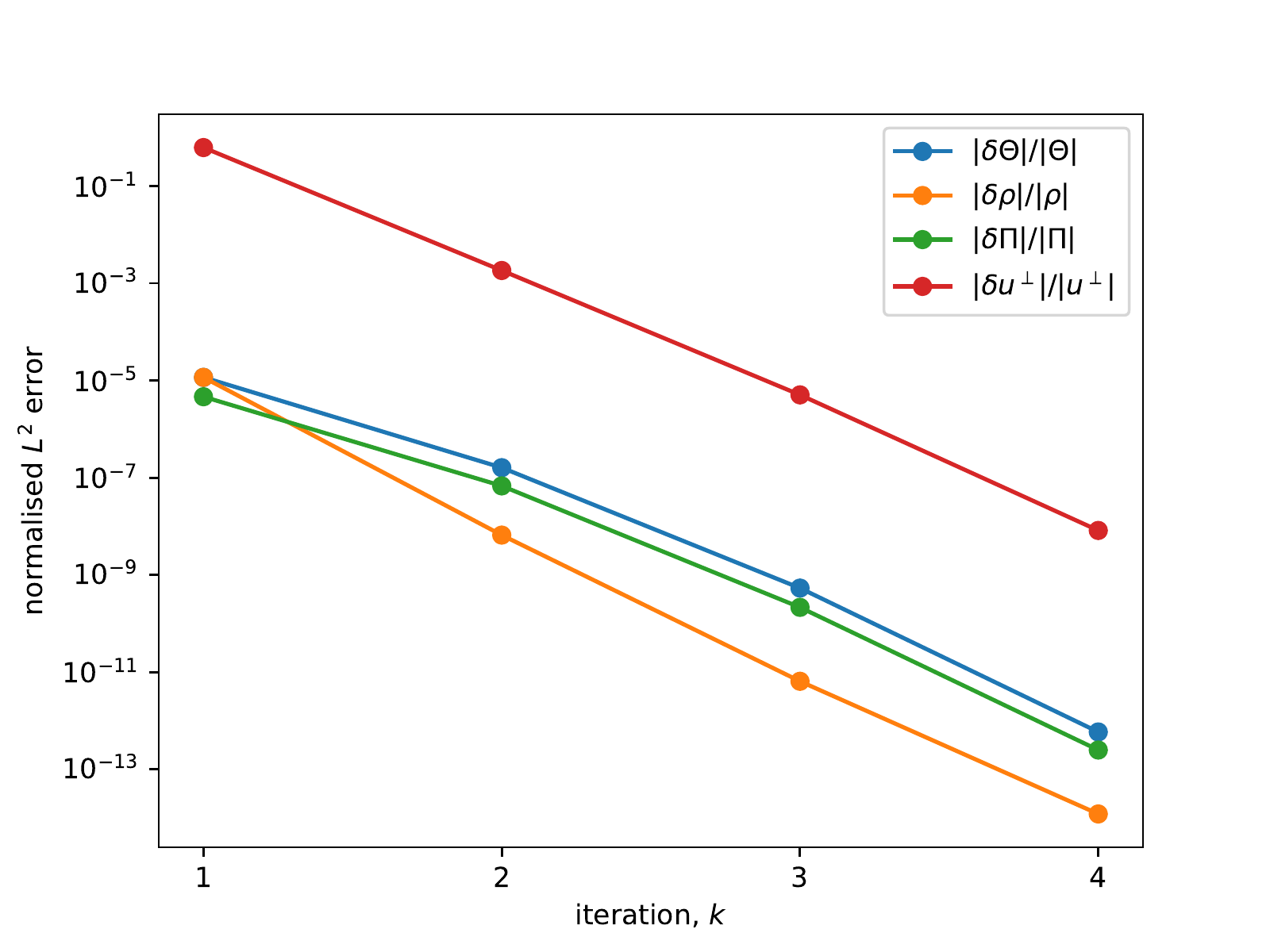}
\caption{Non-hydrostatic gravity wave test case: 
Left: power associated with energetic exchanges, days 3 to 10. Right: typical 
convergence of the normalised $L^2$ magnitude of the solution updates with Newton iteration
(at 1 hour).}
\label{fig::energetics_2_gw}
\end{center}
\end{figure}

The energy partitions are given for the final 12 minutes of the simulation in Fig. \ref{fig::energetics_1_gw}.
The signature of the internal gravity wave is clearly visible in the oscillation between potential
and internal energy. The horizontal dynamics are also observed in the longer time scale oscillation
of kinetic to internal energy exchanges. This is also observed in the power exchanges as given in 
Fig. \ref{fig::energetics_2_gw}.

Figure \ref{fig::energetics_2_gw} also shows nonlinear convergence of the vertical solver at time
$t=1$ hour. Due to the more regular vertical structure than the baroclinic 
test case the convergence of the nonlinear solver is faster and more regular, taking just 4 iterations.

\begin{figure}[!hbtp]
\begin{center}
\includegraphics[width=0.48\textwidth,height=0.36\textwidth]{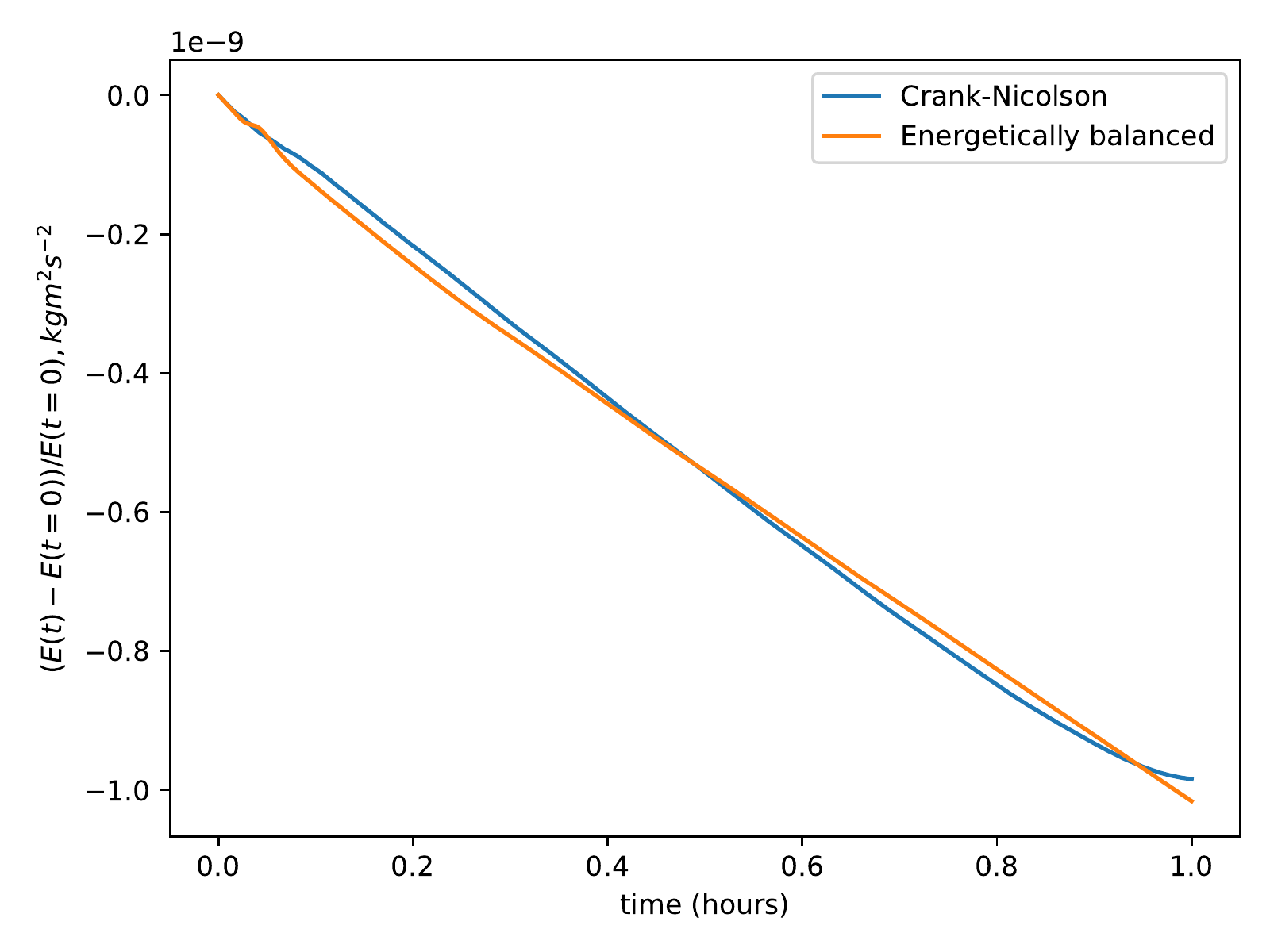}
\includegraphics[width=0.48\textwidth,height=0.36\textwidth]{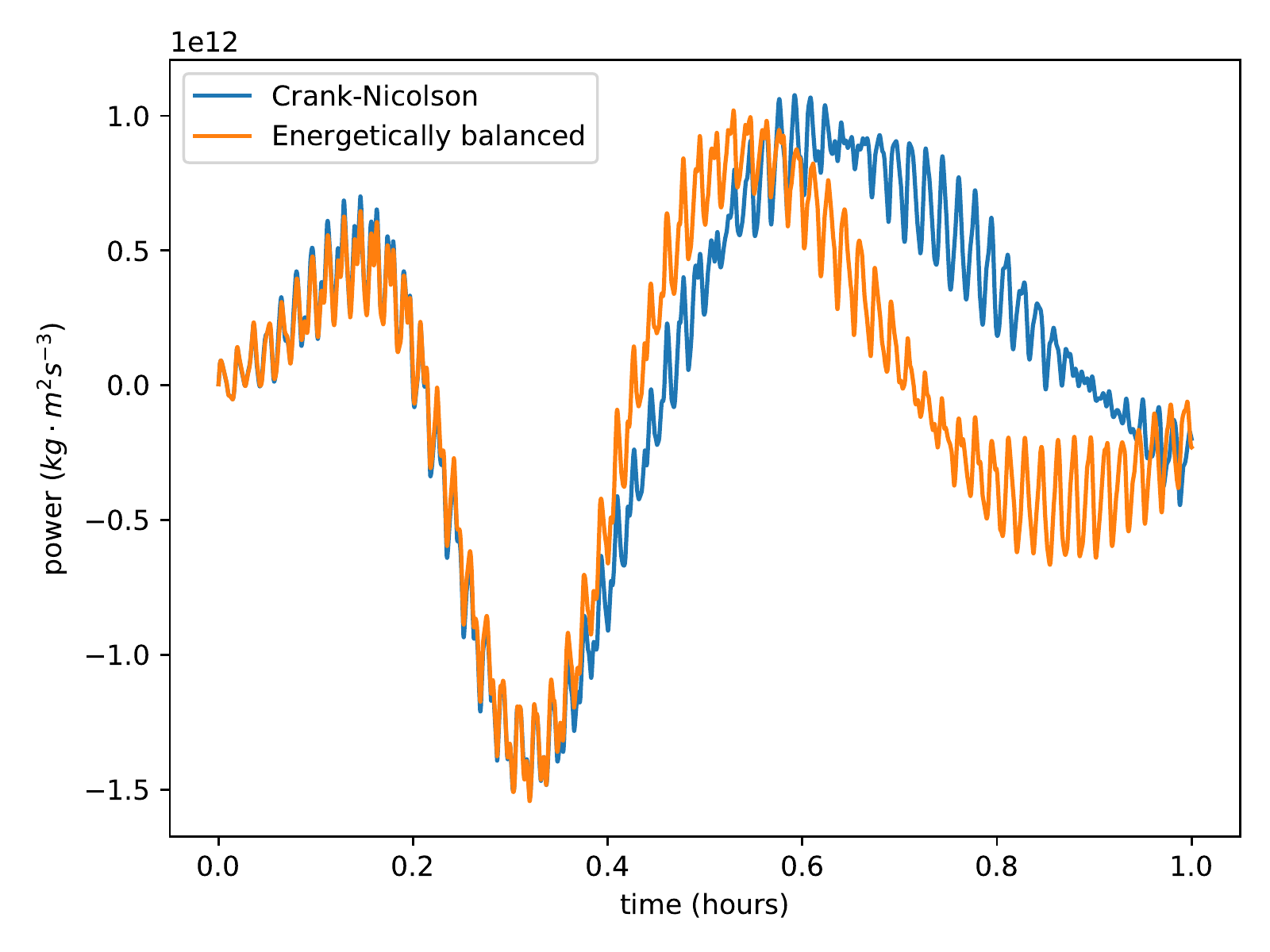}
\caption{
Non-hydrostatic gravity wave test case: 
comparisons of energetically balanced to Crank-Nicolson vertical integrator, energy conservation errors (left)
and vertical internal to kinetic energy exchanges (right).}
\label{fig::energetics_3_gw}
\end{center}
\end{figure}

Figure \ref{fig::energetics_3_gw} shows the energy conservation errors for the non-hydrostatic
gravity wave test case using both the energetically balanced and Crank-Nicolson vertical integrators.
There is little difference between the errors for the two different vertical integrators, which are
due to the horizontal dissipation terms, or
indeed the time variation in kinetic, potential or internal energy (not shown). However over the 
course of the simulation a significant difference develops in the vertical kinetic to internal 
energy exchange. This does not seem to have an appreciable impact on the dynamics over the course
of the simulation, but may possibly become significant over longer simulations.

\subsection{3D rising bubble}

The model is further validated for the case of a 3D rising bubble of potential temperature within a high 
resolution Cartesian geometry.
This involves re-configuring the model geometry as a single panel on a flat surface with 
periodic boundaries, while retaining the original boundary conditions in the vertical dimension. The test is
performed within a high resolution domain of size $1000.0\mathrm{m}\times 1000.0\mathrm{m}\times 1500.0\mathrm{m}$
using $24\times 24$ elements of degree $p=3$ in the horizontal and 150 elements of degree $p=1$ in the vertical and a
time step of $\Delta t = 0.02\mathrm{s}$ in order to resolve the horizontal sound waves. 
The dissipation terms and coefficients are identical to those used in the non-hydrostatic gravity wave test case
(biharmonic viscosity applied to the horizontal momentum and potential temperature advection equations, and 
no dissipation of any kind applied to the vertical dynamics). Notably there is no upwinding
or viscosity within the vertical system, so there is nothing to prevent overshoots or oscillations as the 
gradient of the bubble steepens.

\begin{figure}[!hbtp]
\begin{center}
\includegraphics[width=0.48\textwidth,height=0.36\textwidth]{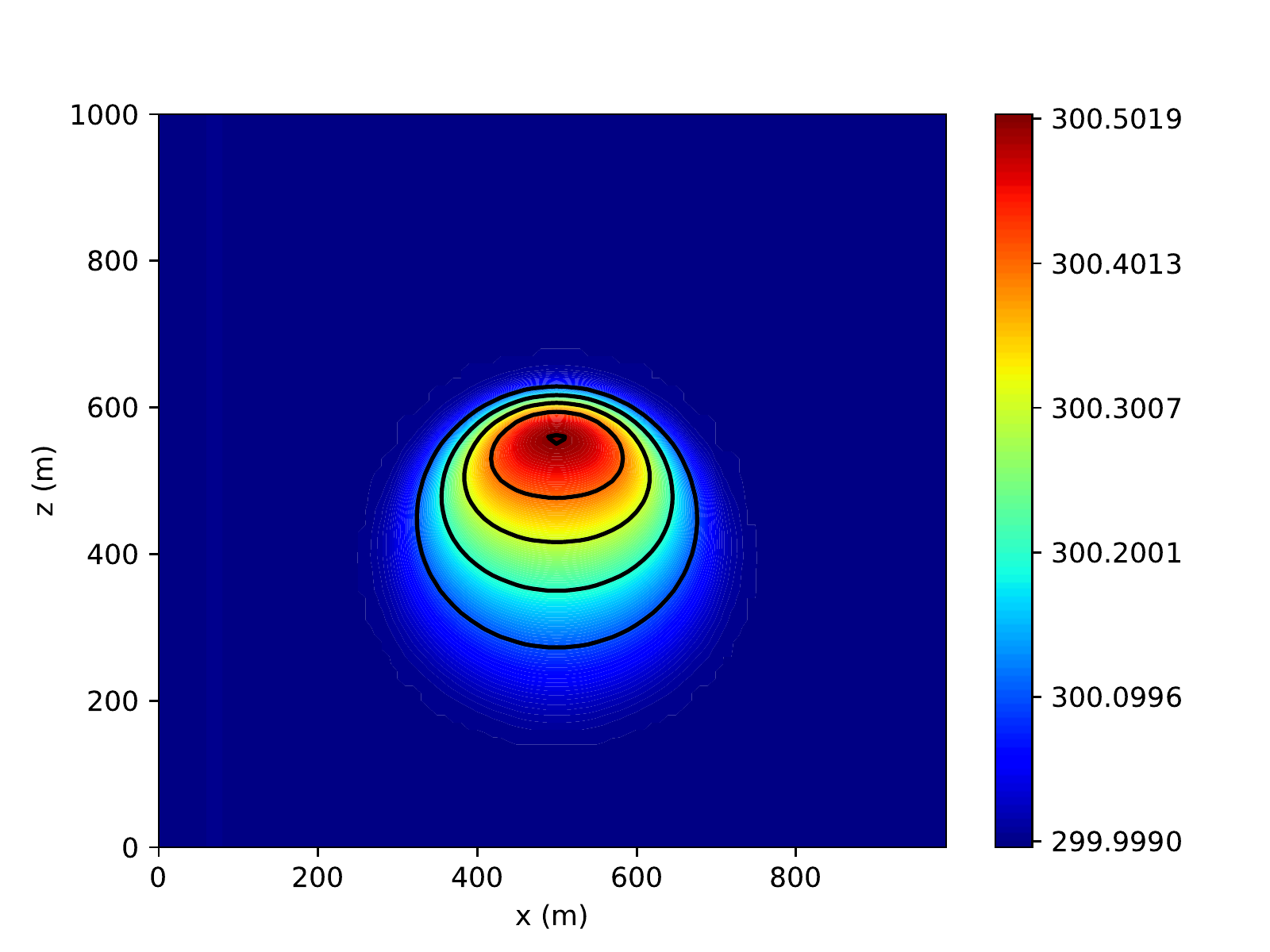}
\includegraphics[width=0.48\textwidth,height=0.36\textwidth]{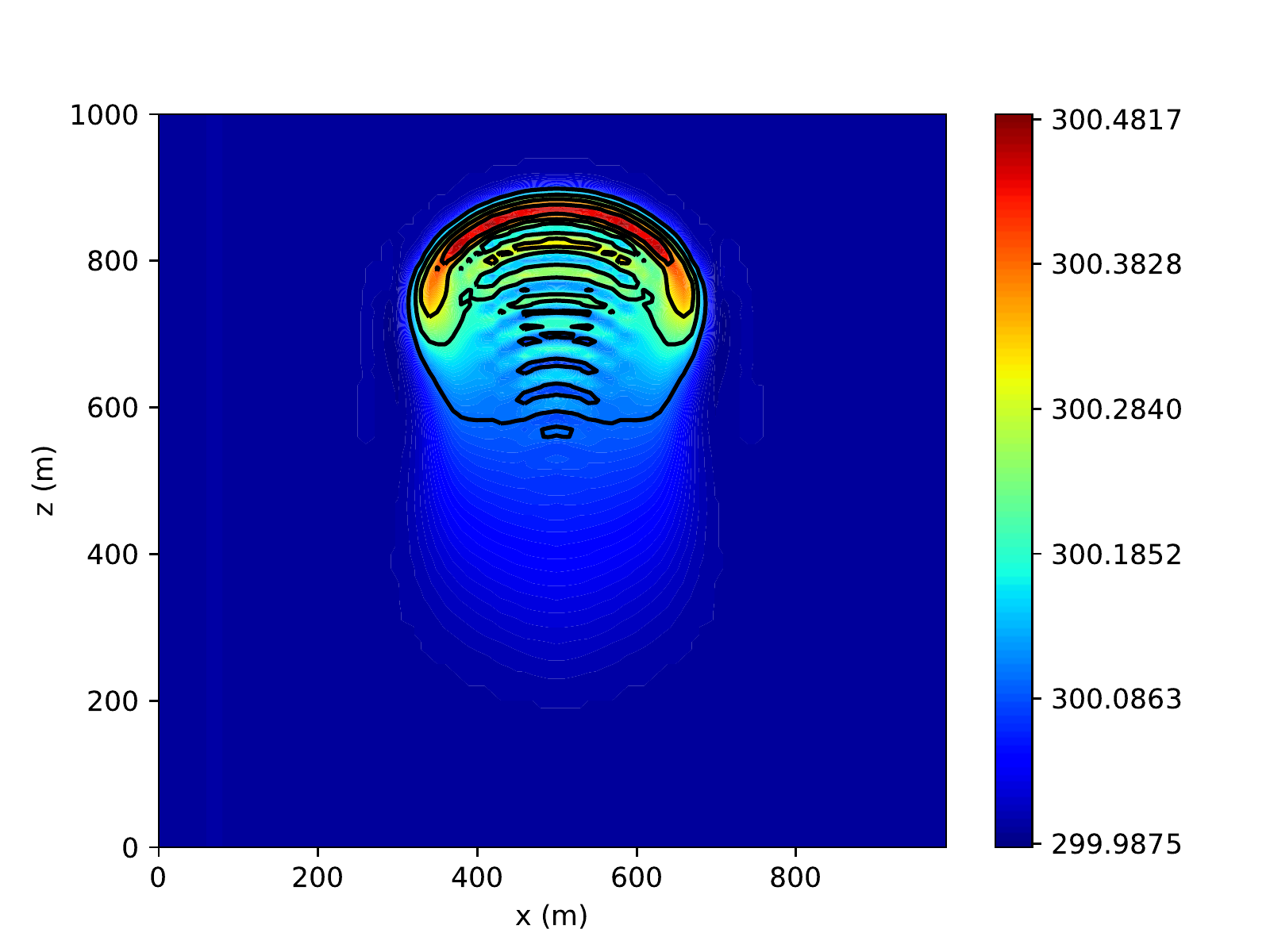}
\caption{
3D rising bubble test case: potential temperature, $\theta_h$, cross section at $y=500.0\mathrm{m}$ at time 
$200.0\mathrm{s}$ (left) and time $400.0\mathrm{s}$ (right).
Contours are from $300.1^{\circ}\mathrm{K}$ to $300.5^{\circ}\mathrm{K}$ in increments of $0.1^{\circ}\mathrm{K}$.}
\label{fig::bubble}
\end{center}
\end{figure}

The bubble is initialised against a constant background state of $\theta_0 = 300.0^{\circ}\mathrm{K}$ as
$\theta' = 0.25(1.0 + \cos(\pi r/r_0))$ for $r\le r_0$, where $r = \sqrt{(x-250.0)^2 + (y-250.0)^2 + (z-350.0)^2}$
and $r_0 = 250.0$ \cite{Melvin19}, and the Exner pressure is initialised from a state of hydrostatic balance
as $\Pi = c_p(1 - gz/c_p\theta_0)$ \cite{AbdiAndGiraldo16}. The initial density is then derived via the equation
of state \eqref{eq::eos}. The position, deformation and amplitude of the bubble compare well to previous results \cite{Giraldo13,Melvin19}, 
as observed in Fig. \ref{fig::bubble}. However the absence of any viscous or upwinding terms from the vertical
discretisation ensures that as the bubble ascends a spurious oscillation develops within its wake, and the amplitude
of the disturbance is not strictly monotone. In order
to suppress this oscillation some form of vertical upwinding should be implemented, as is the case in other
non-hydrostatic models \cite{Melvin19,WCB20}.

\begin{figure}[!hbtp]
\begin{center}
\includegraphics[width=0.48\textwidth,height=0.36\textwidth]{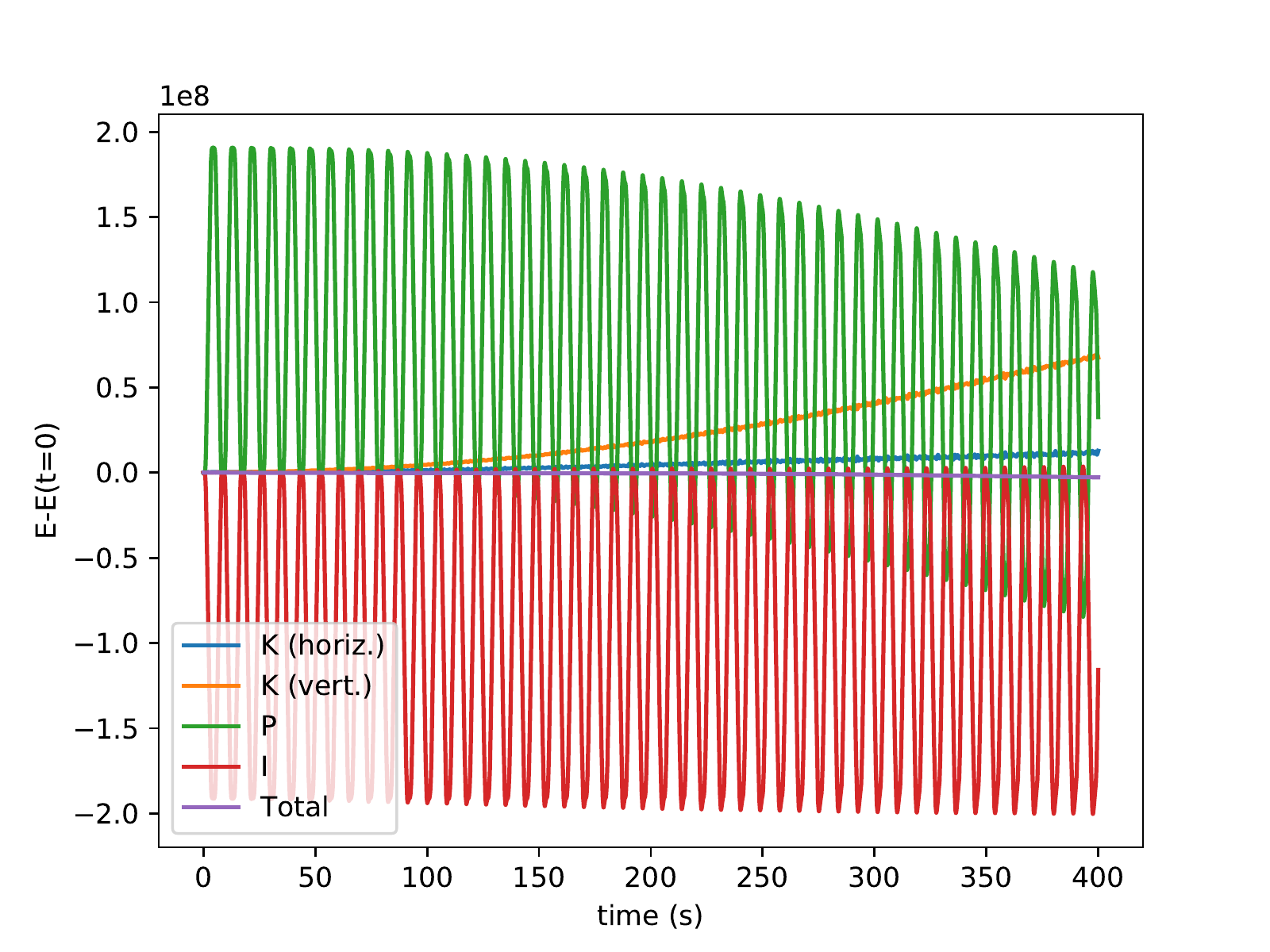}
\includegraphics[width=0.48\textwidth,height=0.36\textwidth]{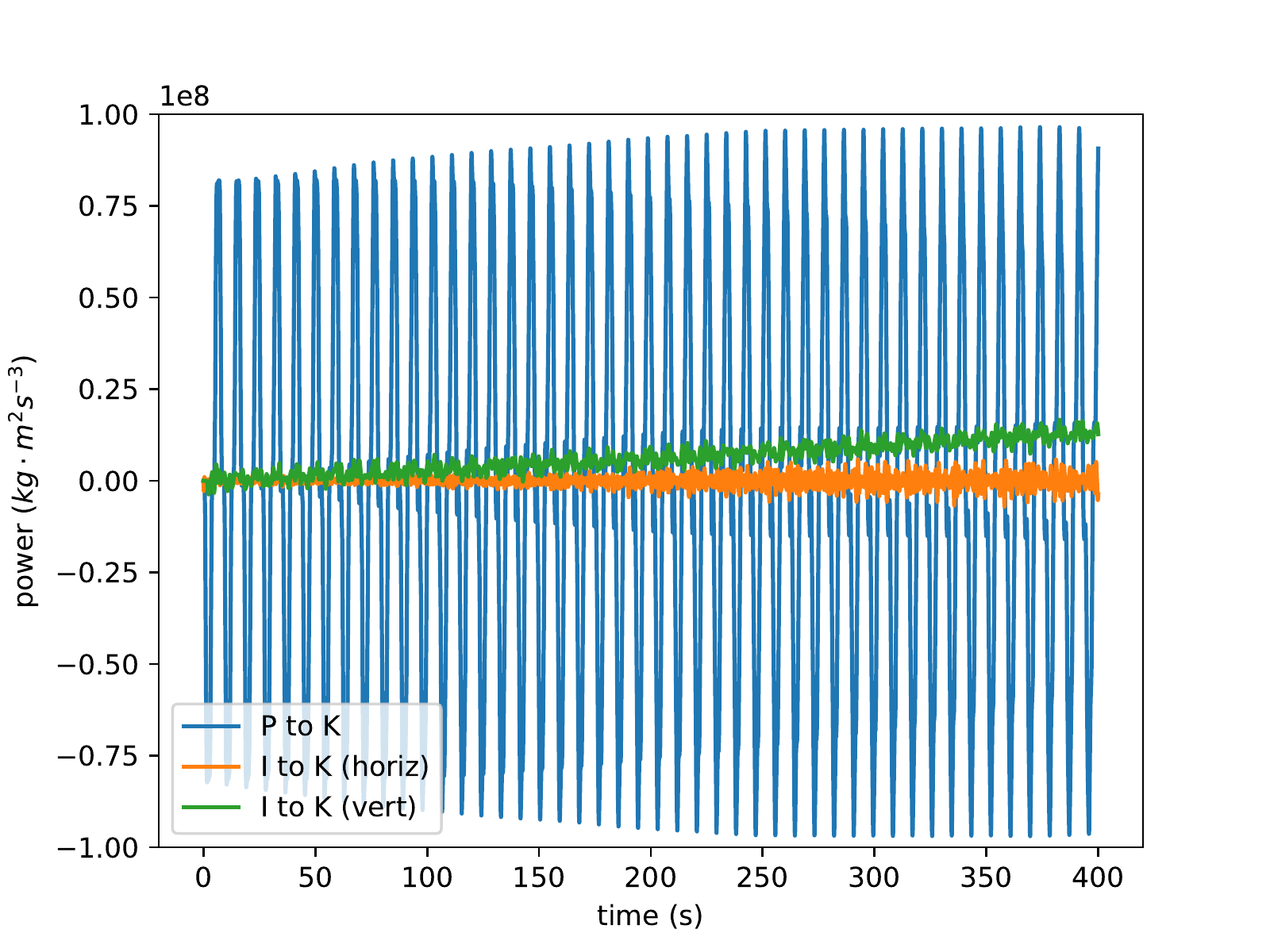}
\caption{3D rising bubble test case: variation in energy from initial state (left) and power exchanges (right).}
\label{fig::bubble_energy}
\end{center}
\end{figure}

Figure \ref{fig::bubble_energy} details the energy variation and power exchanges with time for the rising 
bubble. In addition to the steady decrease in potential energy and increase in vertical kinetic energy as 
the bubble rises, these clearly show the 
signature of an internal gravity wave of period approximately $8\mathrm{s}$ that oscillates within the domain.
The vertical integrator takes just three iterations to converge for the rising bubble test case (with the vertical 
velocity removed from the list of convergence criteria, since the bubble is initially at rest). This is perhaps
because the aspect ratio between the horizontal and vertical scales is unity, such that the requirement that 
the horizontal acoustic modes be resolved in time implies that the vertical modes are resolved also.

\begin{figure}[!hbtp]
\begin{center}
\includegraphics[width=0.48\textwidth,height=0.36\textwidth]{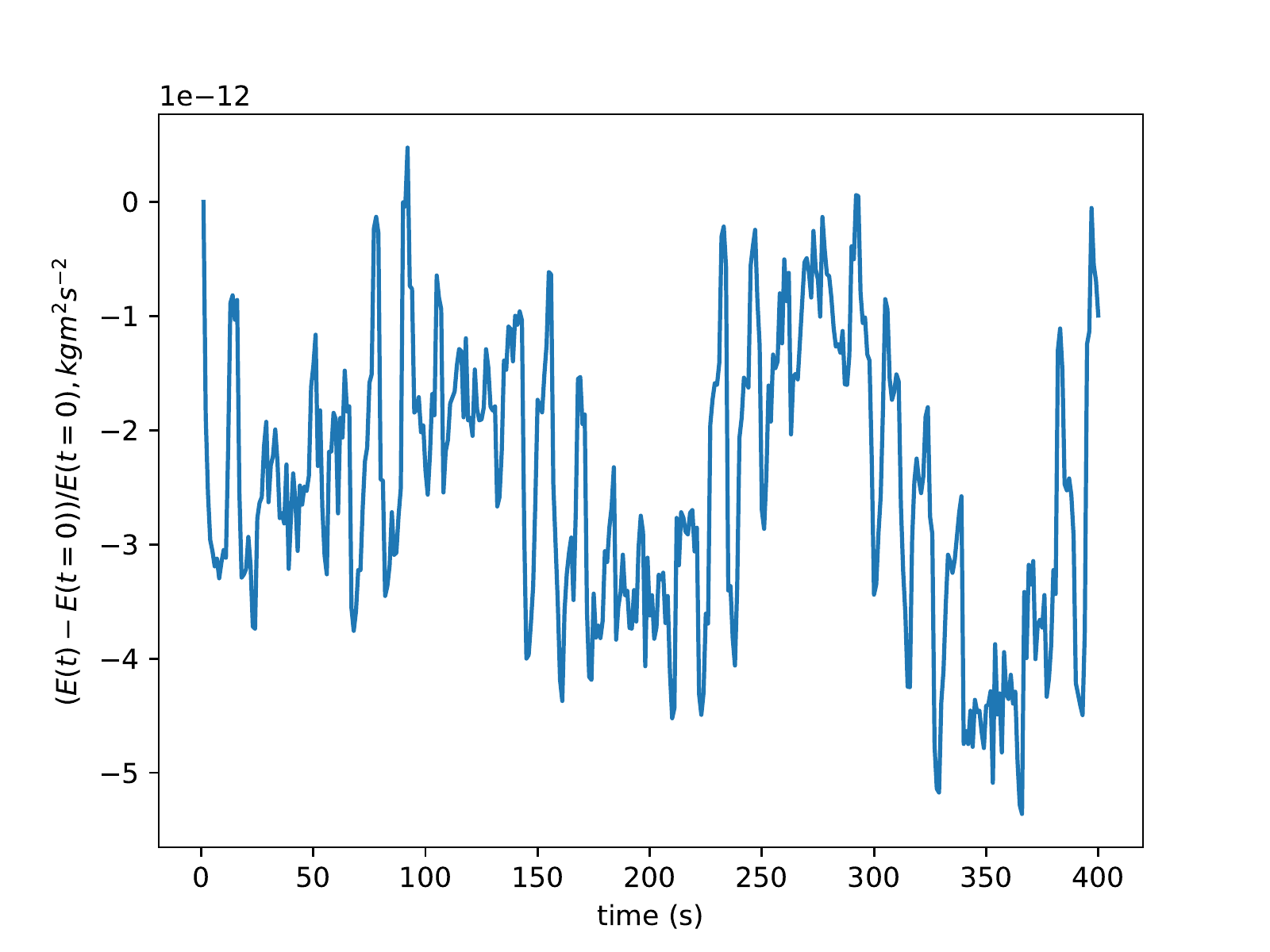}
\caption{3D rising bubble test case (vertical dynamics only): Global total energy conservation error.}
\label{fig::bubble_energy_vert}
\end{center}
\end{figure}

In order to further validate the energy conservation of the new integrator the 3D rising bubble test
case is run with the vertical dynamics only. Since this test is run on a planar geometry, there is no aliasing
error associated with the inexact integration of transcendental metric terms, as is the case for the spherical
geometry. While the dynamics are compressible, the ascension of the bubble occurs primarily in the incompressible
regime, so only the gravity and acoustic modes are discernible without horizontal motions. Without horizontal
acoustic modes, the code can be run with a greatly increased time step, and here we use a step size of 
$\Delta t = 1s$, as well as $4\times 4$ elements of degree $p=3$ at each of the $150$ vertical levels, and no 
horizontal variation in the initial condition. As observed in Fig. \ref{fig::bubble_energy_vert} the energy conservation
errors are comparable to those for the planetary configuration as given in Fig. \ref{fig::conservation},
for which both the metric terms and the vertical structure are more complex.

The differences between the energetics using the energetically balanced and Crank-Nicolson vertical
integrator have also been compared for the 3D rising bubble test case (not shown). Owing to the small time
step required to resolve the horizontal sound waves (which is the same as the time scale for the vertical 
sound waves), and the simple vertical structure of this test case, the differences between the energetics 
for the two vertical integrators are marginal and only barely perceptible.

\section{Conclusions}

An energetically balanced time integrator for vertical atmospheric dynamics is presented. The integrator
balances energetic exchanges associated with vertical motions by preserving the skew-symmetric property
of the non-canonical Hamiltonian form of the compressible Euler equations in space and and the exact
integration of the variational derivatives on the Hamiltonian in time. The
computational efficiency of the integrator is accelerated via a reduction of the four component system
to a single equation for the density weighted potential temperature at each nonlinear iteration. The
integrator is implemented within the context of a horizontally explicit/vertically implicit scheme in 
which the vertical time stepping is centered so as to allow for second order accuracy and exact integration 
in time. The integrator demonstrates robust convergence at both global and non-hydrostatic resolutions.

In future work the issue of the time splitting error associated with the density and density weighted 
potential temperature fluxes will be addressed. These can most likely be reformulated so as to suppress
this error and ensure consistency between time levels. A full three dimensional version of the implicit
integrator will also be investigated, and further examinations between the energetically balanced formulation
and existing schemes will be explored.

\section{Acknowledgments}

David Lee would like to thank Drs. Marcus Thatcher and John McGregor for their continued
support and access to computing resources. 
This project was supported by resources and expertise provided by CSIRO IMT Scientific Computing.


\end{document}